\documentclass[12pt]{article}
\usepackage{graphicx}

\usepackage{amsmath,amssymb,amsthm, amscd, mathrsfs}
\usepackage{ytableau,mathrsfs}
  
  \makeatletter
  \@addtoreset{equation}{section}
  \makeatother

\usepackage[all]{xypic}
\usepackage{tikz,tikz-cd}
\usepackage{longtable}
\allowdisplaybreaks

\newcommand{\Z}{\mathbb{Z}}
\newcommand{\R}{\mathbb{R}}
\newcommand{\C}{\mathbb{C}}
\newcommand{\Q}{\mathbb{Q}}
\newcommand{\A}{\mathbb{A}}
\newcommand{\vphi}{\varphi}
\newcommand{\Cox}{\mrm{Cox}}
\newcommand{\Hom}{\mrm{Hom}}

\newcommand{\mcl}{\mathcal}
\newcommand{\mfk}{\mathfrak}
\newcommand{\mrm}{\mathrm}
\newcommand{\mbf}{\mathbf}

\renewcommand{\tilde}[1]{\widetilde{#1}}
\newtheorem{dfn}{Definition}[section]
\newtheorem{thm}[dfn]{Theorem}
\newtheorem{lem}[dfn]{Lemma}
\newtheorem{prop}[dfn]{Proposition}

\theoremstyle{definition}
\newtheorem{rem}{Remark}[section]
\newtheorem{rems}{Remarks}[section]
\newtheorem{ex}{Example}[section]

\begin{document}

\begin{center} 
{\bf {\LARGE Moduli of $G$-constellations and crepant resolutions II: the Craw-Ishii conjecture}}
\end{center} 
\vspace{0.4cm}

\begin{center}
{\large Ryo Yamagishi}
\end{center} 
\vspace{0.4cm}

\begin{abstract}
For any given finite subgroup $G\subset SL_3(\C)$, we show that every projective crepant resolution $X$ of the quotient variety $\C^3/G$ is isomorphic to the moduli space of $\theta$-stable $G$-constellations for a generic stability condition $\theta$, as conjectured by Craw and Ishii. We also show that generators of the Cox ring of $X$ can be obtained from semi-invariants for representations of the McKay quiver of $G$.
\end{abstract}

\renewcommand{\thefootnote}{\fnsymbol{footnote}} 

\section{Introduction}\label{1}

Let $G$ be a finite subgroup of $SL_n(\C)$. The philosophy of the McKay correspondence is that the geometric properties of a crepant resolution of $\C^n/G$ should reflect the representation-theoretic properties of $G$ (cf. \cite{Re2}). A natural way to see the explicit correspondence is to realize a crepant resolution of $\C^n/G$ as the $G$-Hilbert scheme (see e.g. \cite{Cr}, \cite{CCL}) or, more generally, a moduli space of $G$-constellations. In dimension two, it is well known that the unique crepant resolution is obtained as the $G$-Hilbert scheme. In dimension three, it was also proved by Craw and Ishii \cite{CI} that every projective crepant resolution is obtained as the moduli space $\mcl{M}_{\theta}$ of $G$-constellations for a suitable stability condition $\theta$ when $G$ is abelian. In general $\mcl{M}_{\theta}$ for any finite subgroup $G\subset SL_3(\C)$ is known to be a crepant resolution for any generic stability condition $\theta$ by the result of Bridgeland, King and Reid \cite{BKR}, and it is conjectured that the result by Craw and Ishii still holds for non-abelian $G$ \cite[Conjecture 1.4]{IINdC}. This conjecture has been known to hold in special cases including resolutions whose fibers have dimension at most one \cite{NdCS}, \cite{W}, and iterated $G$-Hilbert schemes \cite{IINdC}. The main result of this article is that this conjecture is affirmative in general:

\begin{thm}(=Theorem \ref{main2})\label{main1}
Let $G\subset SL_3(\C)$ be a finite subgroup. Then every projective crepant resolution of $\C^3/G$ is isomorphic to the moduli space of $\theta$-stable $G$-constellations for a generic stability condition $\theta$. 
\end{thm}

In this paper we also give an explicit construction of a torus-equivariant morphism from the spectrum of the Cox ring of a crepant resolution of $\C^3/G$ to (a certain quotient of) a space of $G$-constellations, as a generalization of the construction in Part I \cite{Y2} of the present paper. In \cite{Y2}, we gave a necessary and sufficient condition for a given (not-necessarily-projective) crepant resolution $X\to \C^n/G$ to admit a moduli description in terms of a natural $G$-constellation family when $G$ is abelian. The key idea to prove this result was to construct an equivariant morphism between certain toric varieties so that it descends to a birational morphism from $X$ to a candidate of a fine moduli space. In this paper we generalize this construction to non-abelian $G$ so that we obtain an embedding $\mcl{W}\hookrightarrow Y_\mcl{W}$ of (a quotient of) a space of $G$-constellations into a toric variety such that it induces an embedding of $\mcl{M}_\theta$ into a toric variety for any generic $\theta$ (Theorem \ref{embedding}). This embedding of $\mcl{M}_\theta$ is so nice that it enables us to understand the behaviour of $\mcl{M}_\theta$ under variation of $\theta$ through the change of the ambient toric variety.  We will demonstrate this embedding and investigate the geometric structure of the moduli space explicitly via concrete examples in Section \ref{3.5}. The embedding will also play a key role to prove Theorem \ref{main1} so that we can overcome one of difficulties which arises when we extend the main result of \cite{CI} to non-abelian cases.

The embedding $\mcl{W}\hookrightarrow Y_\mcl{W}$ will be constructed by using the Cox ring $\Cox(X)$ of $X=\mcl{M}_\theta$, whose generators play a similar role to homogeneous coordinates of a projective space. On the other hand $\mcl{M}_\theta$ is regarded as a moduli space of certain representations of the McKay quiver of $G$, and semi- (or relative) invariants for these representations also play a role of homogeneous coordinates of $\mcl{M}_\theta$. Then the equivariant morphism
$$\mfk{X}:=\mrm{Spec}\,\Cox(X)\to \mcl{W}$$
is obtained by assigning to each semi-invariant a global section of a line bundle on $X$. Moreover, as an application of Theorem \ref{main1}, we will show that generators of $\Cox(X)$ can be computed (at least theoretically) once we know a generating system of semi-invariants for the McKay quiver. More precisely, there is a ring homomorphism $\vphi$ from the ring $\C[\mcl{W}]$ of semi-invariants to the (semi-)invariant ring $\C[\C^3]^{[G,G]}$ such that the ``associated elements" to the images of homogeneous generators of the integral closure of $\C[\mcl{W}]$ under $\vphi$ give a generating system of $\Cox(X)$ (see Proposition \ref{prop:generate} for the precise statement).

The strategy for the proof of Theorem \ref{main1} is based on the one given by Craw and Ishii for abelian cases. To explain this, let us recall that each GIT-chamber $C\subset \Theta$ in the space of stability conditions gives a fine moduli space $X:=\mfk{M}_C$ and its universal family $\mcl{U}_C$. By the result of \cite{BKR}, the Fourier-Mukai transform by $\mcl{U}_C$ gives an equivalence between the bounded derived category of coherent sheaves on $X$ and that of $G$-equivariant sheaves on $\C^3$. This induces an isomorphism
$$\phi_C^*:\Theta\to F^1$$
from the space of stability conditions to the subspace $F^1$ inside the Grothendieck group $K(X)_\R$ generated by sheaves whose supports have codimension at least one. The proof of the theorem will be done by studying the behaviour of $\phi_C^*(C)$ under crossing a wall of {\it Type 0} inside $\Theta$. Here a wall of Type 0 means a codimension-one face $W$ of $\overline{C}$ whose associated morphism $X\to \bar{\mfk{M}}_W$ induced by variation of GIT is an isomorphism. The significant fact to prove Theorem \ref{main1} is that, similarly to abelian cases, crossing such a wall gives a wall-crossing in $F^1$ as well, up to tensoring by a line bundle (Proposition \ref{prop:twist}). The corresponding result \cite[Proposition 7.3]{CI} for abelian $G$ was proven by using a certain rigidity result for which the abelian assumption was essential. In the present paper we bypass this issue by showing that the structure of the unstable locus $D_W$ associated to the wall $W$ is simple enough (Lemma \ref{lem:structure}) to prove Proposition \ref{prop:twist} even if the rigidity does not hold.

Once Proposition \ref{prop:twist} is established, we can find a wall $\tilde{B}$ in $\Theta$ which realizes the small contraction $X\to X_B$ for any given flop $X\to X_B\leftarrow X'$ after crossing walls of Type 0. Using technical Lemma \ref{lem:vertical}, which is deduced from the embedding $\mcl{W}\hookrightarrow Y_\mcl{W}$, we can confirm that crossing the wall $\tilde{B}$ gives the flop $X\dashrightarrow X'$. Since any two projective crepant resolutions are connected by a sequence of flops (cf. \cite{K}), this proves the theorem. In \cite{CI} the fact that $\tilde{B}$ induces a flop was also proven, but again by relying essentially on the abelian assumption. Thus, we can avoid this issue by using the embedding $\mcl{W}\hookrightarrow Y_\mcl{W}$. We remark that, in the construction of the embedding, we essentially use the finite generation of $\Cox(X)$, which is guaranteed by \cite{BCHM} (see Subsection \ref{3.1}).

This paper is organized as follows. In Section \ref{2} we review the construction of a moduli space of $G$-constellations and semi-invariants of quiver representations. In Section \ref{3}, we introduce the Cox ring of a crepant resolution and relate it with the semi-invariants for representations of the McKay quiver of $G$ so that we obtain a nice embedding of a moduli space into  a toric variety (Theorem \ref{embedding}). We also give concrete examples in Subsection \ref{3.5}. In Section \ref{4} we give a proof of the main result (Theorem \ref{main2}). In the appendix, symbols used in \S2-\S4 are listed.

\vspace{5mm}
\noindent{\underline{Conventions \& Notations}\\
In this paper all schemes are algebraic ones defined over $\C$. By a {\it point} of a scheme, we always mean a closed point. For an abelian group $A$, we denote by $A_\R$ its scalar extension $A\otimes_\Z \R$. Similarly, for an homomorphism $f$ of abelian groups, $f_\R$ denotes the scalar extension. For a finite dimensional $\R$-vector space $V$ and a subset $S\subset V$, we denote by $\mrm{Cone}(S)$ the cone generated by $S$. See the appendix for the list of symbols used in \S2-\S4.

\vspace{5mm}
\noindent{\bf Acknowledgements}\\
The author is grateful to Alastair Craw and Akira Ishii for their useful comments to the draft of the present paper. This work was partially supported by the grant (NTU-110VV006) of the Taiwan Ministry of Education, by World Premier International Research Center Initiative (WPI), MEXT, Japan, by JSPS KAKENHI Grant Number JP19K14504, and by Research Project Grant RPG-2021-149 from The Leverhulme Trust.

\section{Moduli spaces of $G$-constellations}\label{2}
\subsection{GIT-construction}\label{2.1} 
In this section we review the construction of a moduli space of $G$-constellations for an arbitrary finite subgroup $G\subset SL(V)$ with $V=\C^n$ (cf. \cite[Section 2]{CI}). See e.g. \cite{CMT} and \cite{Y2} for the special case where $G$ is abelian.

We first recall that a {\it $G$-constellation} $F$ is a $G$-equivariant $\C[V]$-module which is isomorphic to the regular representation $R=\bigoplus_{g\in G}\C g$ as a $\C G$-module where the $G$-action on $\C[V]$ comes from the inclusion $G\subset SL(V)$. Although $G$-constellations are defined as $G$-equivariant coherent sheaves on $V$ in \cite{CI}, we identify them with their global sections since $V$ is affine.

Let $\mrm{Irr}(G)$ be the set of isomorphism classes of irreducible representations of $G$. As is well known, $R$ is isomorphic to $\bigoplus_{\rho\in \mrm{Irr}(G)}\rho^{\oplus \dim \rho}$ as a $\C G$-module. We will construct a moduli space of $G$-constellations by taking a quotient of the affine scheme
\begin{equation}\label{N}
\mcl{N}=\{B\in \mrm{Hom}_{\C G}(V^*\otimes_\C R,R) \mid B\wedge B=0\}
\end{equation}
where $V^*$ is the dual representation of $V$ and
\begin{equation}\label{eq:B wedge B}
B\wedge B\in \mrm{Hom}_{\C[G]}(V^*\otimes_\C V^*\otimes_\C R,R)
\end{equation}
is defined by
$$(B\wedge B)(x\otimes y\otimes a)=B(x\otimes B(y\otimes a))-B(y\otimes B(x\otimes a)).$$

We may regard $\mcl{N}$ as a space parametrizing all $G$-constellations, and two elements of $\mcl{N}$ are isomorphic as $G$-equivariant $\C[V]$-modules if and only if they lie in the same $GL_R$-orbit where $GL_R:=\mrm{Aut}_{\C G}(R)$ is the $G$-equivariant automorphism group of $R$, which acts on $\mcl{N}$ by conjugation. By Schur's lemma, $GL_R$ is written as a direct product $\prod_{\rho\in\mrm{Irr}(G)}GL_{\dim \rho}(\C)$. Note that the diagonal subgroup $\C^*$ acts trivially on $\mcl{N}$ and thus the action of $GL_R$ descends to that of $PGL_R:=GL_R/\C^*$.

Since the orbit space $\mcl{N}/PGL_R$ does not admit a reasonable structure as a scheme, we apply the construction of geometric invariant theory (GIT). Concretely, we consider a $PGL_R$-linearization of the trivial line bundle of $\mcl{N}$ and take the GIT-quotient with respect to it. Such a linearization is given by a character of $PGL_R$, and the group $\chi(PGL_R)$ of characters can be identified with the space of {\it stability conditions}
$$\bar{\Theta}:=\{\theta\in\mrm{Hom}_\mathbb{Z}(R(G),\Z)\mid \theta(R)=0\}$$
where $R(G)=\bigoplus_{\rho\in\mrm{Irr}(G)} \Z \rho$ is the representation ring of $G$. Note that $\bar{\Theta}$ has rank $s-1$ with $s$ being the number of conjugacy classes in $G$. The identification is given by
$$\begin{aligned}
&\chi(PGL_R) & &\to & &\bar{\Theta}\\
\chi=(\bar{g}_{\rho}&\mapsto \mrm{det}(g_{\rho})^{\theta_\rho},\;\forall \rho\in\mrm{Irr}(G)) & &\mapsto & \theta_\chi=(\rho&\mapsto \theta_\rho,\;\forall \rho\in\mrm{Irr}(G)).
\end{aligned}$$

For a character $\chi\in \chi(PGL_R)$ we consider the $\chi${\it-semistable locus} 
$$\mcl{N}^{\chi\text{-ss}}=\{p\in \mcl{N} \mid f(p)\ne0\text{ for some }f\in A_{\chi^k}\text{ and }k\in\mathbb{N}\}$$
where $A_\chi$ is the $\C$-vector space consisting of regular functions $f$ on $\mcl{N}$ such that $g\cdot f=\chi(g)f$ for any $g\in PGL_R$. As an open subset of $\mcl{N}^{\chi\text{-ss}}$, the $\chi${\it-stable locus} $\mcl{N}^{\chi\text{-s}}$ is defined as the set of points $p\in\mcl{N}^{\chi\text{-ss}}$ such that the orbit $PGL_R\cdot p$ is closed in $\mcl{N}^{\chi\text{-ss}}$ and that the stabilizer subgroup of $p$ in $PGL_R$ is finite. We say that a character $\chi$ is {\it generic} if $\mcl{N}^{\chi\text{-ss}}=\mcl{N}^{\chi\text{-s}}$.

For any character $\chi$, we have a $PGL_R$-invariant morphism 
$$\pi_\chi:\mcl{N}^{\chi\text{-ss}}\to \mcl{N}/\!/_\chi PGL_R:=\mrm{Proj}\left(\bigoplus_{k=0}^\infty A_{\chi^k}\right)$$
which gives a one-to-one correspondence between the set of closed $PGL_R$-orbits in $\mcl{N}^{\chi\text{-ss}}$ and the set of closed points in $\mcl{N}/\!/_\chi PGL_R$ (cf. \cite{MFK}). We call $\mcl{N}/\!/_\chi PGL_R$ the {\it GIT-quotient} of $\mcl{N}$ by $PGL_R$ with respect to $\chi$. If $\chi$ is generic, then $\mcl{N}/\!/_\chi PGL_R$ is the orbit space of $\mcl{N}^{\chi\text{-ss}}$, and in such a case $\mcl{N}/\!/_\chi PGL_R$ is called the {\it geometric quotient}.

The stability of a point of $\mcl{N}$ can be rephrased in representation-theoretic terms by the result of King:

\begin{prop}\label{King} \cite[Proposition 3.1]{Kin}
For a point $p\in\mcl{N}$ and a character $\chi\in\chi(PGL_R)$, $p$ is in $\mcl{N}^{\chi\text{-ss}}$ (resp. $\mcl{N}^{\chi\text{-s}}$) if and only if the corresponding $G$-constellation $F_p$ is $\theta_\chi$-semistable (resp. $\theta_\chi$-stable), that is,  $\theta_\chi(M)\ge0$ (resp. $>0$) for any nonzero $G$-equivariant $\C[V]$-submodule $M\subsetneq F_p$.
\end{prop}

We will use the scalar extension $\Theta:=\bar{\Theta}\otimes_\Z \R$
as the space of stability conditions rather than $\bar{\Theta}$ itself in order to consider polyhedral cones inside $\Theta$ (or $\chi(PGL_R)_\R$). For $\theta\in\Theta$ and the corresponding character $\chi_\theta$, we denote $\mcl{N}/\!/_{\chi_\theta} PGL_R$ also by $\mcl{M}_\theta$. From King's result above, we see that $\Theta$ is divided into finitely many open chambers, on which the semistable locus is constant. We call such a chamber a {\it GIT-chamber} of $\Theta$. Since $G$-constellations are regular representations as $\C G$-modules, supporting hyperplanes of the GIT-chambers are of the form $\left\{\sum_{\rho\in \mrm{Irr}(G)} a_\rho \theta_\rho=0\right\}$ with $a_\rho\in\{0,1,\dots,\dim \rho\}$.

Typical examples of $G$-constellations are given by the coordinate rings of  free $G$-orbits, that is, orbits $G\cdot x\subset \C^n$ such that $\sharp(G\cdot x)=\sharp G$. Such $G$-constellations are simple $\C[V]$-modules and hence stable for any $\theta\in\Theta$. As a natural generalization of free $G$-orbits, $G$-constellations obtained as quotients of $\C[V]$ by its ideal are called $G$-{\it clusters}. Since $G$-clusters are characterized as the $G$-constellations generated by the trivial representation, the moduli space of $G$-clusters, called the $G$-{\it Hilbert scheme}, is realized as $\mcl{M}_{\theta+}$ for $\theta_+\in\Theta$ satisfying $\theta_+(\rho)>0$ for any nontrivial $\rho\in \mrm{Irr}(G)$.

By \cite[Proposition 2.2]{CI}, there exists an irreducible component of $\mcl{M}_0$ which is isomorphic to $\C^n/G$ such that  its general points parametrize the coordinate rings of free $G$-orbits in $\C^n$. Accordingly, there is an irreducible component $\mcl{V}\subset\mcl{N}$ (with the reduced structure) such that the restriction of the quotient map $\mcl{N}\to \mcl{M}_0$ to $\mcl{V}$ factors into
\begin{itemize}
\item the morphism $\mcl{V}\to\C^n/G$ which sends a $G$-constellation to its support, and
\item a closed immersion $\C^n/G\to \mcl{M}_0$ onto an irreducible component of $\mcl{M}_0$.\end{itemize}
In Subsection 3.3 we will give a more explicit description of $\mcl{V}$ (cf. Proposition \ref{component V}). Note that the morphism $\mfk{M}_\theta:= \mcl{V}/\!/_{\chi_\theta} PGL_R\to \C^n/G$ induced by the inclusion $\mcl{V}\cap\mcl{N}^{\chi_\theta\text{-ss}}\hookrightarrow \mcl{V}$ is a birational morphism for any $\theta\in\Theta$. The irreducible component $\mfk{M}_\theta$ of $\mcl{M}_\theta$ is sometimes called the {\it coherent component}. For a GIT-chamber $C\subset \Theta$, the moduli space $\mfk{M}_\theta$ is constant for $\theta\in C$ and thus we also denote it by $\mfk{M}_C$.

The main result of \cite{BKR} shows that the $G$-Hilbert scheme for $G\subset SL_3(\C)$ is a crepant resolution of $\C^3/G$. More generally, one can repeat the proof of this result to show that the moduli space $\mcl{M}_\theta$ of $\theta$-semistable $G$-constellations for any generic $\theta\in \Theta$ is a crepant resolution of $\C^3/G$. In particular $\mcl{M}_\theta$ is irreducible and hence $\mcl{M}_\theta=\mfk{M}_\theta$ in this case.

As explained in \cite[\S2.1]{CI}, the moduli space $\mfk{M}_C$ for a GIT-chamber $C$ admits a universal family $\mcl{U}_C$ on $\mfk{M}_C \times \C^n$, and the associated tautological bundle $(p_{\mfk{M}_C})_*\mcl{U}_C$ on $\mfk{M}_C$ admits a decomposition
$$(p_{\mfk{M}_C})_*\mcl{U}_C=\bigoplus_{\rho \in \mrm{Irr}(G)} \mcl{R}_\rho\otimes \rho$$
into isotypic components where $p_{\mfk{M}_C}:\mfk{M}_C \times \C^n\to\mfk{M}_C$ is the projection. Note that each $\mcl{R}_\rho$ is a vector bundle on $\mfk{M}_C$ of rank $\dim \rho$. Since tensoring by a line bundle to $(p_{\mfk{M}_C})_*\mcl{U}_C$ has no effect on the parametrized $G$-constellations, we assume that $\mcl{R}_{\rho_0}\cong \mcl{O}_{\mfk{M}_C}$ for the trivial representation $\rho_0$ whenever we talk about the universal family or the tautological bundle determined by $C$. Later we will identify $G$-constellations with representations of a quiver in Subsection \ref{2.3}. We will also use the associated tautological bundle $\mcl{R}_C$ on $\mfk{M}_C$ as a moduli space of quiver representations rather than use the bundle $(p_{\mfk{M}_C})_*\mcl{U}_C$.

For later use, we introduce the notion of an {\it orbit cone} for an affine scheme acted on by an algebraic group, following \cite[Definition 2.1]{BH}:

\begin{dfn}\label{orbit cone}
Let $H$ be an algebraic group acting on an affine algebraic scheme $Z$. For a point $z\in Z$, the {\bf orbit cone} of $z$ is the polyhedral cone $C_z\subset\chi(H)_\R$ generated by characters $\chi$ such that $z$ is $\chi$-semistable.
\end{dfn}

\begin{rem}
The orbit cone of $z\in Z$ depends only on the $H$-orbit of $z$. If $Z$ is the affine scheme $\mcl{N}$ acted on by $H=PGL_R$ as above, the orbit cone of a point of $\mcl{N}$ corresponding to a $G$-constellation $F$ is identified with a cone in $\Theta$ consisting of stability conditions $\theta$ such that $F$ is $\theta$-semistable. Moreover, for a GIT-chamber $C\subset\Theta$, its closure $\overline{C}$ is equal to the intersection of all the orbit cones for points of $\mcl{N}^{C\text{-ss}}$. 
\end{rem}

\subsection{Semi-invariants for quiver representations}\label{2.2} 

Let us consider the subgroup $SL_R:=\prod_{\rho\in\mrm{Irr}(G)}SL_{\dim \rho}(\C)$ of $GL_R$. Since $SL_R$ is reductive, there are finitely many generators of the invariant ring of the coordinate ring $\C[\mcl{N}]$ of the affine scheme $\mcl{N}$. We call elements of $\C[\mcl{N}]^{SL_R}$ {\it semi-invariants (with respect to the action of $GL_R$ on $\mcl{N}$).} Note that we have a decomposition of this ring into homogeneous components:
$$\C[\mcl{N}]^{SL_R}=\bigoplus_{\chi\in \chi(PGL_R)} A_\chi.$$

As we will see in the next subsection, $\mcl{N}$ is identified with a space of quiver representations (with relations). In general the ring of semi-invariants for representations of a quiver (in characteristic zero) is generated by {\it determinantal semi-invariants} (see Theorem \ref{thm:semi-invariants} below). In this subsection we explain how to obtain such semi-invariants.

Firstly, we briefly recall the concept of representations of a quiver. See e.g. \cite{DeW2},\cite{Kir} for more details. Let $Q$ be a quiver, that is, an oriented graph consisting of a finite set $I$ of vertices and a finite set $A$ of arrows between two vertices. Here we allow $Q$ to have loops and oriented cycles. For any map $\alpha:I\to \Z_{\ge0}$, we define a space $\mrm{Rep}(\alpha)$ of representations of $Q$ with a dimension vector $\alpha$ as
$$\mrm{Rep}(\alpha)=\bigoplus_{a\in A} \mrm{Hom}_\C\left(\C^{\alpha(t_a)},\C^{\alpha(h_a)}\right)$$
where $t_a$ and $h_a$ denote the tail and the head of $a$ respectively. For a representation $H=(H_a)_{a\in A}\in \mrm{Rep}(\alpha)$ and a path
$$p=a_1\cdots a_k\quad(a_i\in A,\,h_{a_i}=t_{a_{i+1}}),$$
we set $H_p=H_{a_k}\circ\cdots \circ H_{a_1}$. For the trivial path $e_v$ at a vertex $v$, we set $H_{e_v}=\mrm{id}_{\C^{\alpha(v)}}$. $\mrm{Rep}(\alpha)$ is acted on by $GL_\alpha:=\prod_{v\in I} GL_{\alpha(v)}$ by conjugation.

Let $v_1,\dots,v_r,w_1,\dots,w_s$ be vertices of $Q$ with possible repetition such that
$$\sum_{i=1}^r \alpha(v_i)=\sum_{j=1}^s \alpha(w_j),$$
and let $p_{i,j}$ be a $\C$-linear combination of paths from $v_i$ to $w_j$. Then the function $f:\mrm{Rep}(\alpha)\to \C$ defined by
\begin{equation}\label{determinant}
H\mapsto \det
\begin{pmatrix}H_{p_{1,1}}&\cdots&H_{p_{1,s}}\\
\vdots& \ddots&\vdots\\
H_{p_{r,1}}&\cdots&H_{p_{r,s}}
\end{pmatrix}
\end{equation}
gives a semi-invariant satisfying
$$(g_v)_{v\in I}\cdot f=\left(\prod_{v\in I} \det g_v^{c_v}\right)f$$
for $(g_v)_{v\in I}\in GL_\alpha$ where
\begin{equation}\label{eq:weight}
c_v=\sharp\{j \mid w_j=v\}-\sharp\{i \mid v_i=v\}.
\end{equation}
A semi-invariant obtained in this way is called a {\it determinantal semi-invariant} of weight $(c_v)_{v\in I}$. Note that weights are identified with characters of $GL_\alpha$.

\begin{thm}(\cite{DeW1},\cite{DZ},\cite{SV})\label{thm:semi-invariants}
The ring $\C[\mrm{Rep}(\alpha)]^{SL_\alpha}$ of semi-invariants on $\mrm{Rep}(\alpha)$ with respect to the $GL_\alpha$-action is generated by determinantal semi-invariants.
\end{thm}

This was first proved by Derksen and Weyman \cite{DeW1} under the assumption that $Q$ has no oriented cycles, and generalized by Domokos and Zubkov \cite{DZ} and by Schofield and van den Bergh \cite{SV}. We will apply this theorem to representations of the McKay quiver introduced in the next subsection.

\subsection{$G$-constellations as representations of the McKay quiver}\label{2.3}

In this subsection we realize the affine scheme $\mcl{N}$ in (\ref{N}) as a space of quiver representations.

The {\it McKay quiver} $Q_G$ of $G$ is the quiver having $\mrm{Irr}(G)$ as the set of vertices and having $a_{\rho,\rho'}:=\dim_\C(\mrm{Hom}_{\C G}(V^*\otimes_\C \rho,\rho'))$ arrows from $\rho$ to $\rho'$ for any $\rho,\rho'\in\mrm{Irr}(G)$. Decomposing $R$ into irreducible summands, we see that $\mrm{Hom}_{\C G}(V^*\otimes_\C R,R)$ is isomorphic to the representation space $\mrm{Rep}(\alpha_G)$ of $Q_G$ for the dimension vector $\alpha_G:\rho\mapsto \dim \rho$ (cf. \cite{IN}). Thus, $\mcl{N}$ can be realized as a $GL_R(=GL_{\alpha_G})$-invariant closed subscheme of $\mrm{Rep}(\alpha_G)$.

Via the identification of a $G$-constellation with a representation of $Q_G$, the moduli space $\mfk{M}_C$ of $G$-constellations constructed in Subsection \ref{2.1} is also regarded as a moduli space of representations of $Q_G$. Then the tautological bundle $\mcl{R}_C$ on $\mfk{M}_C$ as the latter moduli space admits a decomposition $\mcl{R}_C=\bigoplus_{\rho\in\mrm{Irr}(G)} \mcl{R}_\rho$ into the isotypic components. Note that $\mcl{R}_C$ is equipped with a $\C[V]$-action whose restriction to $\C[V]^G$ coincides with the $H^0(\mcl{O}_{\mfk{M}_C})$-action as a vector bundle on $\mfk{M}_C$. Hereafter, by {\it the tautological bundle on} $\mfk{M}_C$, we mean the bundle $\mcl{R}_C$.

\begin{ex}\label{ex:relations}($D_5$-singularity)
We describe the space $\mcl{N}$ explicitly for the case of the $D_5$-singularity. This singularity is given as a quotient of $\C^2$ by a binary dihedral group $G$ of order 12:
$$G=\left\langle
g_1=\begin{pmatrix}
\zeta_6&0\\
0&\zeta_6^{-1}
\end{pmatrix},
g_2=\begin{pmatrix}
0&1\\
-1&0
\end{pmatrix}
\right\rangle\subset SL_2(\C)$$
where $\zeta_6$ is the sixth root of unity. $\mrm{Irr}(G)$ consists of four 1-dimensional representations $\rho_k:G\to\C^*\,(k=0,1,2,3)$ defined by
$$\rho_k(g_j)=\begin{cases}
(-1)^k&\text{ if }j=1\\
i^k&\text{ if }j=2\end{cases},$$
and two 2-dimensional representations $V_1,V_2$; the first one is the inclusion $G\subset SL_2(\C)$ presented above and the other is presented as $$g_1\mapsto\begin{pmatrix}
\zeta^2_6&0\\
0&\zeta^{-2}_6
\end{pmatrix},\quad
g_2\mapsto\begin{pmatrix}
0&1\\
1&0
\end{pmatrix}.$$

Let $e_0,\dots,e_3$ be a basis of $\rho_0,\dots,\rho_3$ respectively. Also, we let $\{v_{11},v_{12}\}$ (resp. $\{v_{21},v_{22}\}$) be a basis of $V_1$ (resp. $V_2$) such that it gives the matrix presentation above. We moreover prepare their copies $\{v'_{11}, v'_{12}\}$ and $\{v'_{21},v'_{22}\}$ for the second direct summands of $V_1$ and $V_2$ in $R$ respectively. When we emphasize that $v_{11}$ and $v_{12}$ are elements of $V^*$, we instead denote them by $x$ and $y$ respectively.

We fix isomorphisms of $\C G$-modules as follows:
{\small
$$\begin{alignedat}{11}
V^*&\otimes \rho_0 & & \cong V_1, &\qquad V^*&\otimes \rho_1 & &\cong & & V_2,\; &\qquad V^*&\otimes \rho_2  &  &\cong & & V_1,\; &\qquad V^*&\otimes \rho_3 & &\cong & & V_2\\
x&\otimes e_0 & & \mapsto v_{11} & x&\otimes e_1 & &\mapsto & & v_{22} & x&\otimes e_2 &  &\mapsto & & v_{11} & x&\otimes e_3 & &\mapsto & & v_{22}\\
y&\otimes e_0& &\mapsto v_{12} & y&\otimes e_1& &\mapsto & i&v_{21} & y&\otimes e_2& &\mapsto -& & v_{12} & y&\otimes e_3& &\mapsto -i & & v_{21}
\end{alignedat}$$
$$\begin{alignedat}{11}
V^*&\otimes V_1 & & \cong & & \rho_0 & &\oplus \rho_2 \oplus & & V_2, &\qquad V^*&\otimes V_2 & & \cong & & \rho_1 & &\oplus  & & \rho_3 \oplus & & V_1,\\
x&\otimes v_{11} & & \mapsto & & & & & & v_{21} & x&\otimes v_{21} &  &\mapsto  & & e_1 & &+ & & e_3 & & \\
x&\otimes v_{12}& & \mapsto & & e_0 & & +e_2& & & x&\otimes v_{22}& & \mapsto & &  & & & & & & v_{12}\\
y&\otimes v_{11}& & \mapsto - & & e_0 & & +e_2& & & y&\otimes v_{21}& & \mapsto & &  & & & & &-& v_{11}\\
y&\otimes v_{12} & & \mapsto & & & & & & v_{22} & y&\otimes v_{22} &  &\mapsto -i & & e_1 & &+i & & e_3. & &
\end{alignedat}$$}
Thus the McKay quiver $Q_G$ is presented as

\begin{figure}[h]
\centerline{
\xymatrix@=36pt{
& \rho_0 \ar@<-0.5ex>[d]_-{A_0} 
& \rho_1 \ar@<-0.5ex>[d]_-{A_1}
& \\
\rho_2 \ar@<-0.5ex>[r]_-{A_2}
& V_1 \ar@<-0.5ex>[l]_-{B_2} \ar@<-0.5ex>[u]_-{B_0}\ar@[r]@<-0.5ex>[r]_-C  
& V_2 \ar@<-0.5ex>[l]_-D \ar@<-0.5ex>[r]_-{B_3}  \ar@<-0.5ex>[u]_-{B_1}
& \rho_3 \ar@<-0.5ex>[l]_-{A_3}.
}
}
\caption{The McKay quiver for $D_5$-singularity}\label{quiver}
\end{figure}

Using the isomorphisms above, matrices associated to the arrows in $Q_G$ define an action of $\C[V]$ on the vector space
{\small$$\C e_0\oplus\cdots \oplus\C e_3 \oplus (\C v_{11}\oplus \C v_{12})\oplus (\C v'_{11}\oplus \C v'_{12})\oplus  (\C v_{21}\oplus \C v_{22})\oplus (\C v'_{21}\oplus \C v'_{22})\cong R.$$}
For instance, a matrix $\begin{pmatrix}a_1\\a_2\end{pmatrix}\,(a_i\in\C)$ corresponding to the arrow $A_1$ in Figure \ref{quiver} lets $V^*$ act on $\C e_1= \rho_1$ as $x\cdot e_1=a_1 v_{22}+a_2 v'_{22}$ and $y\cdot e_1=a_1  iv_{21}+a_2 iv'_{21}$. Recall that the matrices must satisfy the condition $B\wedge B=0$ (see (\ref{eq:B wedge B})). By writing down the conditions $x\cdot (y\cdot v)-y\cdot (x\cdot v)=0$ for each basis $v$ of $R$, one sees that this condition is expressed as the following matrix relations:
{\small
\begin{equation}\label{path}
B_0A_0=B_1A_1=B_2A_2=B_3A_3=A_0B_0-A_2B_2-DC=A_1B_1-A_3B_3-iCD=0.
\end{equation}}
\end{ex}

\begin{rems}
(1) Usually the quiver $Q_G$ itself can be obtained easily from the character table of $G$, without giving explicit $G$-equivariant homomorphisms $V^*\otimes \rho\to \rho'$ as above. However, this computation will be necessary in order to construct embeddings of moduli spaces into toric varieties, as we will see in the next section.\\
(2) The existence of an ideal $I_G$ of the path algebra $\C Q_G$ such that
$$\mcl{N}\cong \{H\in\mrm{Rep}(\alpha_G) \mid H_p=0,\,p\in I_G\}$$
as in the above example is regarded as a consequence of the Morita equivalence between the skew group algebra $\C[V]\ast G$ and a quotient of $\C Q_G$. For another method to calculate $I_G$, see \cite[Theorem 3.2 and \S5]{BSW}.
\end{rems}

\section{Embedding into toric varieties}\label{3}

In this section we explain how we embed moduli spaces $\mfk{M}_\theta$ into toric varieties such that the birational geometry of the embedded varieties is inherited from that of the ambient toric varieties.

\subsection{The Cox ring of a resolution of a quotient singularity}\label{3.1}

In this subsection we treat the Cox ring, denoted by $\Cox(X)$, of a crepant resolution $X$ of a quotient singularity $\C^n/G$ for $G\subset SL_n(\C)$. $\Cox(X)$ is the direct sum $\bigoplus_{L\in \mrm{Pic}(X)} H^0(X,L)$ as an $H^0(\mcl{O}_{V/G})$-module and admits a multiplicative structure inherited from the tensor product of global sections of line bundles (see \cite{ADHL} for generalities regarding Cox rings). Then $\Cox(X)$ is a $\mrm{Pic}(X)$-graded commutative ring and is known to be finitely generated by the fundamental result of Birker-Cascini-Hacon-McKernan \cite{BCHM}. See e.g. the proof of \cite[Lemma 5.3]{BCRSW} as to how \cite[Corollary 1.3.2]{BCHM} can be applied to our situation to deduce the finite generation of the Cox ring. The most important feature of the Cox ring is that $X$ can be recovered from $\Cox(X)$ as a GIT-quotient of $\mrm{Spec}\,\Cox(X)$ by an algebraic torus with respect to a generic stability parameter if $X\to\C^n/G$ is projective (see the next subsection). Note that the isomorphism class of $\Cox(X)$ is independent of the choice of $X$ since crepant resolutions are mutually isomorphic in codimension one. The structure of $\Cox(X)$ has been studied for special classes of $G$ in \cite{D}, \cite{DW}, \cite{G} and so on. See \cite{Y1} for a general treatment. Moreover, algorithms finding generators of $\Cox(X)$ are given in \cite{DK} and \cite{Y1}.

Now we review how to describe the structure of $\Cox(X)$ explicitly. We refer the reader to \cite{Y1} for details. We first recall the notion of {\it age}. For an element $g\in G$, one can choose a basis $x_1,\dots,x_n$ of the dual space $V^*$ of $V=\C^n$ such that $g\cdot x_j=e^{\frac{2\pi i a_j}{r}}x_j$ with an integer $0\le a_j<r$ where $r$ is the order of $G$. The {\it age} of $g$ is defined as $\frac{1}{r}\sum_{j=1}^n a_j$, which is an integer since $g$ is in $SL_n(\C)$. We call an element $g\in G$ {\it junior} if its age is one. Note that the age is invariant under conjugation. We also define a discrete valuation $\nu_g:\C(V)^*\to \Z$ by setting $\nu_g(x_j)=a_j$.

It is shown in \cite{IR} that there is a one-to-one correspondence between the set of conjugacy classes of junior elements in $G$ and the set of irreducible exceptional divisors of the resolution $\pi:X\to V/G$, such that the divisorial valuation $\nu_E:\C(X)^*\to \Z$ along the $\pi$-exceptional divisor $E$ corresponding to $g$ coincides with the restriction of $\frac{1}{r}\nu_g$ to $\C(X)$ (via the birational map $\pi$). Let $E_1,\dots,E_m$ be the set of $\pi$-exceptional irreducible divisors and let $g_k\,(k=1,\dots,m)$ be a representative of the junior conjugacy class corresponding to $E_k$. We also put $\nu_k=\nu_{g_k}$ and $r_k=\sharp\langle g_k\rangle$.

The Cox ring has an embedding into the Laurent polynomial ring
$$\C[V]^{[G,G]}[t_1^{\pm1},\dots,t_m^{\pm1}]$$
with $m$ variables over the invariant subring for the commutator group $[G,G]$ of $G$ \cite[Lemma 4.2]{Y1}. Note that $\C[V]^{[G,G]}$ is naturally graded by the dual abelian group $Ab(G)^\vee=\mrm{Hom}_\Z(Ab(G),\C^*)$ of the abelianization $Ab(G)=G/[G,G]$ so that its degree-zero part is equal to $\C[V]^G$. For $Ab(G)^\vee$-homogeneous generators $f_1,\dots,f_\ell$ of $\C[V]^{[G,G]}$, we associate to them the following set
\begin{equation}\label{gens}
\left\{\tilde{f}_i:=f_i\prod_{k=1}^m t_k^{\nu_k(f_i)}\right\}_{i=1,\dots,\ell}\bigcup \left\{t_k^{-r_k}\right\}_{k=1,\dots,m}
\end{equation}
of $\ell+m$ elements of $\C[V]^{[G,G]}[t_1^{\pm1},\dots,t_m^{\pm1}]$, which lie in $\Cox(X)$. This set is a natural candidate of a generating system of $\Cox(X)$ but it is not a generating system in general. In fact it is a generating system if and only if $f_i$'s satisfy the {\it valuation lifting condition} cf. \cite[Theorem 2.2]{DG}. We do not explain this condition here as we will not use it later.

\begin{rem}\label{rem:cox(V/G)}
The ring $\C[V]^{[G,G]}$ treated here should be understood as the Cox ring $\Cox(V/G)$ of $V/G$ (\cite[Theorem 3.1]{AG}), and the $Ab(G)^\vee$-grading also corresponds to the $\mrm{Cl}(V/G)$-grading on $\Cox(V/G)$. Thus each $f_i$ above defines a Weil divisor $D_i$ on $V/G$, and the associated element $\tilde{f}_i$ to $f_i$ corresponds to a section defining the strict transform in $X$ of $D_i$ under the birational morphism $X\to V/G$.
\end{rem}

The $\mrm{Pic}(X)$-grading of $\Cox(X)\subset\C[V]^{[G,G]}[t_1^{\pm1},\dots,t_m^{\pm1}]$ can be read off from the $\Z^m$-grading given by the multi-degrees for the variables $t_1,\dots,t_m$. More precisely, the $\mrm{Pic}(X)$-grading is given so that the generator $t_k^{-r_k}$ lies in the homogeneous component corresponding to the divisor class of $E_k$. Note that the divisor classes of $E_1,\dots,E_m$ are a $\Q$-basis (but usually not a $\Z$-basis) of $\mrm{Pic}(X)$.

\subsection{Review of birational geometry of $X$}\label{3.2}

As mentioned in the previous subsection, any projective crepant resolution $X$ of $V/G$ has a finitely generated Cox ring. Such a variety is called a {\it Mori dream space}. This notion was originally introduced by Hu-Keel \cite{HK} for a projective variety (over a field), but it makes sense in a relative setting (cf. \cite{O}) so that the theory applies to the resolution $X\to V/G$. A remarkable property of a Mori dream space is that its birational contractions are completely controlled by polyhedral cones in the Picard group. In this section we review this briefly.

Hereafter we assume $X\to V/G$ is projective. We first consider the polyhedral cone generated by the classes of divisors in $\mrm{Pic}(X)_\R$ whose base loci are of codimension greater than one. This cone is called the {\it movable cone} and denoted by $\mrm{Mov}(X)$. In our case, the movable cone is a convex cone which is the union of the nef cones of all relative minimal models $X'\to V/G$ (i.e. projective crepant birational morphisms with $X'$ having at worst $\Q$-factorial and terminal singularities). Here, the nef cone $\mrm{Nef}(X')\subset \mrm{Pic}(X')_\R$ is regarded as a cone inside $\mrm{Pic}(X)_\R$ via the birational map $X\dashrightarrow X'$ which is an isomorphism in codimension one. Note that, in dimension three, relative minimal models are nothing but projective crepant resolutions.

The nef cones $\mrm{Nef}(X')\subset \mrm{Pic}(X)_\R$ are full-dimensional and have disjoint interiors. If we take any $L$ in the interior $\mrm{Amp}(X')$ of $\mrm{Nef}(X')$, the complete linear system of a sufficiently divisible power of $L$ gives the birational map $X\dashrightarrow X'$. Set
$$T_X=\mrm{Hom}(\mrm{Pic}(X),\C^*)\cong (\C^*)^m,$$
which naturally acts on $\mfk{X}:=\mrm{Spec}\,\Cox(X)$. This birational map is the same as the one $X\dashrightarrow \mfk{X}/\!/_{\chi} T_X$ between the GIT-quotients with the character $\chi\in \chi(T_X)$ corresponding to $L$ via the identification $\chi(T_X)_\R=\mrm{Pic}(X)_\R$ (see the proof of \cite[Theorem 6.7]{O}). Moreover, the similar results hold for any {\it rational contraction}, that is, a rational map $f:X\dashrightarrow Z$ which is the composite of a sequence of flops $X\dashrightarrow X'$ and a proper birational morphism $f':X'\to Z$ over $V/G$ such that
\begin{itemize}
\item $Z$ is a normal variety projective over $V/G$, and 
\item $f'$ has connected fibers, i.e. $f'_*\mcl{O}_{X'}=\mcl{O}_Z$.
\end{itemize}
More precisely, for such $f$, there exists $L\in\mrm{Mov}(X)$ whose linear system recovers $f$, and moreover $f$ equals the rational map $X\dashrightarrow \mfk{X}/\!/_{\chi} T_X$ induced by variation of GIT for the character $\chi$ corresponding to $L$. Note that general points of a facet (i.e. codimension-one face) of $\mrm{Mov}(X)$ corresponds to a primitive divisorial rational contraction (i.e., a composite of a sequence of flops and a birational morphism contracting a single divisor which drops the Picard number by one). In contrast, small contractions are realized by taking divisors (or the corresponding characters) from the interior of $\mrm{Mov}(X)$. In particular every flop $X\dashrightarrow X'$ corresponds to the common facet $\mrm{Nef}(X)\cap \mrm{Nef}(X')$ of the nef cones of $X$ and $X'$.

\begin{rems}
(1) We can also consider divisors (or corresponding characters) outside the movable cone, but they again give rational contractions (cf. \cite[Corollary 5.8]{O}) and thus considering the movable cone is essential. 

(2) Classifying the GIT-quotients of $\mfk{X}$ (or equivalently, rational contractions of $X$) gives a fan structure on $\chi(T_X)_\R$ (or on $\mrm{Pic}(X)_\R$). This fan is called the {\it GIT-fan} for the $T_X$-action on $\mfk{X}$. When $X$ is toric, this is also called the {\it secondary fan}. See \cite[Ch. 14]{CLS} for a more detailed exposition of this fan.
\end{rems}

\begin{ex}\label{abel1}(Abelian group of type $\frac{1}{10}(1,3,6)$)\\
We illustrate how to compute the movable cone and nef cones for a crepant resolution of $\C^n/G$ for abelian $G$ through a particular example. Here we treat the case when $G\subset SL_3(\C)$ is the cyclic group generated by
$$g=\begin{pmatrix}\zeta_{10}& 0& 0\\0& \zeta_{10}^3& 0\\ 0& 0& \zeta_{10}^6\end{pmatrix}$$
with $\zeta_{10}$ the primitive 10th root of unity. In this case $g,g^2,g^4,g^5,g^7$ are the junior elements, and the Cox ring of any crepant resolution $X$ is generated by the following 8 elements in $\C[x,y,z][t_1^{\pm1},\dots,t_5^{\pm1}]$:
\begin{equation}\label{ex:generator}
x t_1 t_2 t_3^2 t_4 t_5^7,\; y t_1^3 t_2^3 t_3 t_4 t_5,\; z t_1^6 t_2 t_3^2  t_5^2,\; t_1^{-10},\; t_2^{-5},\; t_3^{-5},\; t_4^{-2},\; t_5^{-10}
\end{equation}
where $x,y$ and $z$ are the standard coordinates of $\C^3$. Moreover, $X$ is a toric variety corresponding to a fan $\Sigma$ obtained as a triangulation of the cone $\sigma_0=\mrm{Cone}(e_1,e_2,e_3)\subset\R^3$ such that the set of rays (i.e. 1-dimensional cones) of $\Sigma$ is generated by $e_1,e_2,e_3,v_1,\dots,v_5$ where $e_j$ is the standard basis of $\R^3$ and $v_k$'s are given as
$$v_1=\begin{pmatrix}1\\3\\6\end{pmatrix},\;v_2=\begin{pmatrix}1\\3\\1\end{pmatrix},\; v_3=\begin{pmatrix}2\\1\\2\end{pmatrix},\; v_4=\begin{pmatrix}1\\1\\0\end{pmatrix},\;v_5=\begin{pmatrix}7\\1\\2\end{pmatrix}$$
so that $v_k$ corresponds to $E_k$. See e.g. \cite[\S3]{Y2} for details. One can check that there are six crepant resolutions, which are in fact all projective. These resolutions $X_1,\dots,X_6$ fit into the following diagram of flops:\\
\[\xymatrix@=36pt{
&X_3 \ar@{<-->}[d]_-{} 
&
&
&\\
X_2 \ar@{<-->}[r]_{}
&X_1 \ar@{<-->}[r]_{}
&X_4 \ar@{<-->}[r]_{} 
&X_5 \ar@{<-->}[r]_{} 
&X_6.
}\]
The cross sections by the hyperplane $\{(a,b,c)\in \R^3 \mid a+b+c=1\}$ of the fans $\Sigma_1$ and $\Sigma_6$ of $X_1$ and $X_6$ respectively are given in Figure \ref{fan1}. Among $X_i$'s, the $G$-Hilbert scheme corresponds to $X_1$ (cf. \cite{CR}).
\begin{figure}[h]
\begin{center}
 \includegraphics{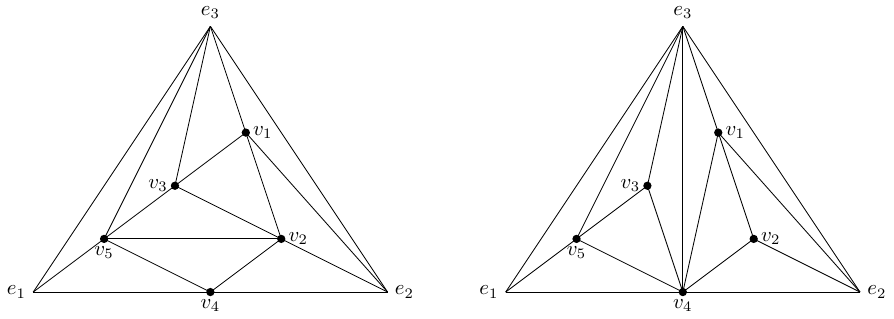}
\caption{The cross sections of the fans $\Sigma_1$ (left) and $\Sigma_6$ (right)}\label{fan1}
\end{center}
\end{figure}

As is well known, each $i$-dimensional cone $\tau\in\Sigma_1$ represents a torus orbit $O_\tau$ in $X_1$ of codimension $i$ (see e.g. \cite[\S3.2]{CLS}). $\tau$ also defines a cone $C_\tau\subset\mrm{Pic}(X_1)_\R$ generated by line bundles having global sections which do not vanish on $O_\tau$. Under the identification of $\mrm{Pic}(X_1)_\R$ with $\chi(T_X)_\R$, $C_\tau$ is the same as the orbit cone (see Definition \ref{orbit cone}) of some (and hence any) lift $\bar{x}\in\mfk{X}$ of some (and hence any) $x\in O_\tau$. $C_\tau$ is explicitly computed as the cone generated by the weights of the generators in (\ref{ex:generator}) whose corresponding rays are not contained in $\tau$. For example, $C_\tau$ for $\tau=\mrm{Cone}(e_3,v_1)\in\Sigma_1$ is generated by the weights of the generators in (\ref{ex:generator}) with the third and fourth ones removed.

The movable cone is obtained as the intersection $\bigcap_{i=1}^5 C_{\tau_i}$ where $\tau_i=\R_{\ge0}v_i$. Similarly, the nef cone of $X_1$ is the intersection $\bigcap_{\tau\in\Sigma_1} C_\tau$. In our case, however, we do not have to compute $C_\tau$ for all $\tau\in \Sigma_1$. $\mrm{Nef}(X_1)$ is given as the intersection of $\mrm{Mov}(X)$ and just one cone $C_\tau$ for $\tau=\mrm{Cone}(v_2,v_4,v_5)$ since $\Sigma_1$ is clearly the only triangulation containing this three-dimensional cone.

The above computations of cones can be done, for instance by using the package ``Polyhedra" of Macaulay2 \cite{GS}. Then we see that $\mrm{Mov}(X_1)$ has 9 rays and 7 facets. The fact that $\mrm{Mov}(X_1)$ is non-simplicial comes from non-uniqueness of primitive rational contractions of exceptional divisors. For example, as we can see from Figure \ref{fan1}, the divisor $E_3\subset X_1$ is mapped to a curve by a primitive contraction while it is mapped to a point by a primitive contraction inside $X_6$. We can also see that $\mrm{Nef}(X_1)$ is a simplicial cone having facets corresponding to 2 primitive divisorial contractions for $E_1,E_3$ and 3 flops along the curves $E_2\cap E_4, E_2\cap E_5$, and $E_4\cap E_5$.

Later we will compare cones in $\mrm{Mov}(X_1)$ and orbit cones in $\Theta$ for moduli spaces of $G$-constellations (see Example \ref{abel2}).
\end{ex}

\subsection{Construction of embeddings}\label{3.3}
In this subsection we give a construction of nice embeddings of moduli spaces $\mfk{M}_\theta$ into toric varieties.

Let us consider the homomorphism $\iota_X^*:\C[y_1,\dots,y_{\ell+m}]\to \Cox(X)$ sending $y_i$'s to homogeneous generators of $\Cox(X)$ of the form in (\ref{gens}) so that we obtain an embedding $\iota_\mfk{X}$ of $\mfk{X}=\mrm{Spec}(\Cox(X))$ into the affine space $\A^{\ell+m}$. We endow $\C[y_1,\dots,y_{\ell+m}]$ with a $\Z^m$-grading so that $\iota_X^*$ preserves the grading. Recall that a character $\chi\in \chi(T_X)_\R = \mrm{Pic}(X)_\R\cong\R^m$ corresponding to an ample divisor on $X$ satisfies $\mfk{X}/\!/_\chi T_X=X$ by the GIT construction (see the previous subsection). We can also consider the GIT-quotient $Y=\A^{\ell+m}/\!/_\chi T_X$ for the same character, which is a normal toric variety such that $\mrm{Pic}(Y)\cong\chi(T_X)\cong\mrm{Pic}(X)$ \cite{Co}, and $\iota_\mfk{X}$ induces an embedding $\iota_X:X\hookrightarrow Y$ between the quotients.

\begin{rem}
If the generators of the Cox ring are minimal, then $\iota_X$ is a {\it neat embedding} in the sense of \cite{Ro}. See \cite[\S2.3]{Ro} for details.
\end{rem}

This embedding has a nice property from the viewpoint of birational geometry. As explained in the previous subsection, every rational contraction $X\dashrightarrow X'$ is realized as the rational map $X\dashrightarrow \mfk{X}/\!/_{\chi'} T_X$ for some $\chi'$. The same character $\chi'$ gives a rational contraction $Y\dashrightarrow Y'$ as well, and $\iota_\mfk{X}$ again induces an embedding $X'\to Y'$. Note that GIT-chambers for $\A^{\ell+m}$ are finer than those for $\mfk{X}$ in general, and in particular $Y\dashrightarrow Y'$ may not be an isomorphism even if $X\dashrightarrow X'$ is an isomorphism. In other words, the inclusion $\iota_X^*(\mrm{Nef}(Y))\subset\mrm{Nef}(X)$ is strict in general. Note, however, that we have
$$\iota_X^*(\mrm{Mov}(Y))=\mrm{Mov}(X)$$
since every exceptional divisor of $Y\to \A^{\ell+m}/\!/_0 T_X$ restricts to an exceptional divisor of $X\to \C^n/G$ by the construction of $Y$.

Next we construct an embedding of a moduli space $\mfk{M}_\theta$ into a toric variety. Recall that $\mfk{M}_\theta$ is obtained as the GIT-quotient of the variety $\mcl{V}$ by the action of $PGL_R$ with respect to the character $\chi_\theta\in\chi(PGL_R)_\R \cong\Theta$. Let $T_\Theta\subset PGL_R$ be the algebraic torus consisting of (the images of) the component-wise scalar matrices. Since the restriction map $\chi(PGL_R)_\R\to \chi(T_\Theta)_\R$ is an isomorphism, via this identification $\mfk{M}_\theta$ is isomorphic to $\mcl{W}/\!/_{\chi_\theta} T_\Theta$ with the affine categorical quotient
$$\mcl{W}=\mcl{V}/\!/SL_R.$$
Note that the points of $\mcl{W}$ are in bijection with the closed $SL_R$-orbits in $\mcl{V}$.

In the proof of \cite[Proposition 2.2]{CI}, it is shown that $\mcl{M}_0=\mcl{N}/\!/PGL_R$ admits an irreducible component which is isomorphic to $\C^n/G$ (see also Subsection \ref{2.1}) by constructing a morphism $V\to \mcl{M}_0$ induced from a $G$-constellation family $\mcl{F}$ over $V=\C^n$. We will explicitly describe $V\to \mcl{M}_0$ by lifting it to a morphism $V\to \mcl{V}$. To do this, let
$$\mu:\C[V]\to\C[G]\otimes_\C \C[V]$$
be the ring homomorphism corresponding to the action of $G\subset SL(V)$. We let $G$ act on $\C[G]$ as $(g\cdot f)(h)=f(g^{-1}h)$ for any $g,h\in G$ and $f\in \C[G]$. In terms of algebras, the $\C[V]$-action on each fiber $\mcl{F}_p:=\C[G]\,(p\in V)$ of the family $\mcl{F}$ is given by the homomorphism
$$\C[V]\otimes \C[G]\to\C[G];\quad x\otimes f\mapsto (\mu(x)f\;\mrm{mod}\; \mfk{m}_p)$$
where $\mfk{m}_p\subset \C[V]$ is the maximal ideal for $p$. This is indeed $G$-equivariant since we have
$$(g\cdot(x\cdot f))(h)=(x\cdot f)(g^{-1}h)=(\mu(x)f\;\mrm{mod}\; \mfk{m}_p)(g^{-1}h)=x(g^{-1}h\cdot p)f(g^{-1}h)$$
while
$$\begin{aligned}
((g\cdot x)\cdot(g\cdot f))(h)&=(\mu(g\cdot x)(g\cdot f)\;\mrm{mod}\; \mfk{m}_p)(h)=((g\cdot x)(h\cdot p))((g\cdot f)(h))\\
&=x(g^{-1}h\cdot p)f(g^{-1}h)
\end{aligned}$$
for any $g,h\in G,\, x\in \C[V]$ and $f\in \mcl{F}_p=\C[G]$.

For each $\rho\in\mrm{Irr}(G)$, we fix its $\C$-basis $v_1,\dots,v_{\dim \rho}$. The matrix presentation of the representation $\rho$ with respect to this basis gives a matrix $(f_{i,j})_{i,j}$ with entries in $\C[G]$, that is, $g\cdot v_j=\sum_{i=1}^{\dim \rho}f_{i,j}(g)v_i$ for all $g\in G$. Let $v_1^*,\dots,v_{\dim \rho}^*$ be the dual basis to $v_1,\dots,v_{\dim \rho}$. Then the assignment $f_{i,j}\mapsto v_i^*$ gives $G$-equivariant isomorphisms from the vector subspaces $V_{\rho,j}:=\bigoplus_{i=1}^{\dim \rho}\C f_{i,j}\subset\C[G]$ to the dual (or contragredient) representations $\rho^*$ of $\rho$. This can be checked by comparing the coefficients of $v_i$'s in
$$\begin{aligned}
&g^{-1}h\cdot v_j=\sum_{i=1}^{\dim \rho} f_{i,j}(g^{-1}h)v_i=\sum_{i=1}^{\dim \rho} (g\cdot f_{i,j})(h)v_i \quad \text{and}\\
&g^{-1}\cdot (h\cdot v_j)=g^{-1}\cdot \left(\sum_{i=1}^{\dim \rho}f_{i,j}(h)v_i\right)=\sum_{i=1}^{\dim \rho}f_{i,j}(h)\left(\sum_{k=1}^{\dim \rho}f_{k,i}(g^{-1})v_k\right)
\end{aligned}$$
with $g,h\in G$. Since $\C[G]$ is a direct sum of $\{V_{\rho, j}\}_{\rho, j}$ by the character theory for finite groups, we obtain an explicit $G$-equivariant isomorphism $\C[G]\to\bigoplus_{\rho\in\mrm{Irr}(G)}\rho^{\oplus \dim\rho}$. We will describe the family $\mcl{F}$ as quiver representations with the matrix presentation with respect to the basis of $\C[G]$ obtained in this way. To this end, we also fix a basis of the $a_{\rho,\rho'}$-dimensional vector space $\mrm{Hom}_{\C G}(V^*\otimes_\C \rho,\rho')$ for each pair $\rho,\rho'\in\mrm{Irr}(G)$.

For each arrow $a\in A$ from $\rho$ to $\rho'$ in $Q_G$, we have an embedding $\rho'\to V^*\otimes\rho$. Regarding the fixed basis for $\rho$ and $\rho'$, we obtain a $(\dim\rho'\times\dim\rho)$-matrix $\tilde{H}_a$ with entries in $V^*$. Then we obtain a set of matrices $\{H_a^p:=(\tilde{H}_a\;\mrm{mod}\; \mfk{m}_p)\}_{a\in A}$, which can be regarded as a set of representations of $Q_G$ (see Subsection \ref{2.3}) parametrized by $V$.

\begin{lem}
The representation of $Q_G$ corresponding to the $G$-constellation $\mcl{F}_p$ is isomorphic to $\{H_a^p\}_{a\in A}$ for any $p\in V$.
\end{lem}

\begin{proof}
We take any arrow $a\in A$ from $\rho$ to $\rho'$ and set $m=\dim \rho$ and $n=\dim \rho'$. Let $(f_{i,j})_{1\le i,j\le m}$ (resp. $(g_{i,j})_{1\le i,j\le n}$) be the matrix presentation of $\rho^*$ (resp. $\rho'^*$) with respect to the fixed basis. We write the $n\times m$-matrix $\tilde{H}_a$ as $(h_{i,j})_{i,j}$ with $h_{i,j}\in V^*$. Regarded as a representation of $Q_G$, $H_a^p$ defines an action of $V^*$ on the vector space $V_{\rho^*,j}=\bigoplus_{i=1}^m\C f_{i,j}\subset\C[G]$, which is isomorphic to $\rho$ as a $\C G$-module, so that
\begin{equation}\label{Ha}
\begin{aligned}
\sum_{i=1}^m h_{1,i}\cdot f_{i,k} &=\sum_{j=1}^n h_{j,k}(p) g_{1,j}\\
\sum_{i=1}^m h_{2,i}\cdot f_{i,k} &=\sum_{j=1}^n h_{j,k}(p) g_{2,j}\\
&\vdots\\
\sum_{i=1}^m h_{n,i}\cdot f_{i,k} &=\sum_{j=1}^n h_{j,k}(p) g_{n,j}
\end{aligned}
\end{equation}
for each $k=1,\dots,m$. We show that this action coincides with the one for $\mcl{F}_p$. By the definition of the $\C[V]$-module structure on $\mcl{F}_p=\C[G]$, the LHS of the $j$-th row of (\ref{Ha}), as an element of $\mcl{F}_p$, is the function
$$g\mapsto h_{j,1}(g\cdot p)f_{1,k}(g)+h_{j,2}(g\cdot p)f_{2,k}(g)+\cdots+h_{j,m}(g\cdot p)f_{m,k}(g)$$
on $G$ for each $j$. Note that the value of this function at the identity $1\in G$ is $h_{j,k}(p)$ and this is also equal to the value of the RHS of the $j$-th row of (\ref{Ha}) (as a function on $G$) at $1\in G$. By the choice of $\tilde{H}_a$, the LHSs and the RHSs of (\ref{Ha}) give $G$-equivariantly isomorphic basis of $\rho'$. Noticing that $f(g)=g^{-1}\cdot f(1)$ for any $f\in\C[G]$ and $g\in G$, this implies that the two functions attain the same values for all $g\in G$, and hence the claim.
\end{proof}

We will give concrete examples of the representations $\{H_a^p\}_{a\in A}$ in Subsection \ref{3.5}.

The generic representation $\{\tilde{H}_a\}_a$ gives a morphism $V\to\mcl{N}$, which factors through $V\to\mcl{V}$ since $V$ is mapped onto $V/G\subset \mcl{M}_0$. Let $f$ be a semi-invariant function on $\mcl{N}$ and hence on $\mcl{V}$. Then the (scheme-theoretic) zero locus of $f$ on $\mcl{V}$ descends to a Weil divisor $D_f$ on $V/G$. Thus, sending $f$ to the section of the divisorial sheaf $\mcl{O}_{V/G}(D_f)$ defines a ring homomorphism
$$\vphi:\C[\mcl{N}]^{SL_R}\to \Cox(V/G)\cong \C[V]^{[G,G]}$$
(cf. Remark \ref{rem:cox(V/G)}). More explicitly, for a determinantal semi-invariant $f:H\mapsto \det (H_{p_{i,j}})_{i,j}$ defined by paths $p_{i,j}$ as in (\ref{determinant}), $\vphi(f)$ is given as $\det \mcl{H}\in\C[V]$ where $\mcl{H}=(\tilde{H}_{p_{i,j}})_{i,j}$ is defined similarly to $\tilde{H}_a$. Note that if $f$ has weight $(\theta_\rho) _{\rho\in\mrm{Irr}(G)}$, then $\vphi(f)$ is acted on by $G$ by the character
$$\det\left(\bigoplus_{\rho\in\mrm{Irr}(G)} \rho^{\oplus \theta_\rho}\right):G\to \C^*.$$
In particular this naturally gives a homomorphism $\bar{\Theta}\to Ab(G)^\vee$.

The homomorphism $\vphi$ becomes more meaningful when we restrict our attension from $\mcl{N}$ to its irreducible component $\mcl{V}$. Recall that we have the following commutative diagram of affine schemes and their GIT quotinets
\[
\xymatrix{
\mcl{N} \ar[r] & \mcl{N}/\!/SL_R \ar[r]_{/\!/T_\Theta}  \ar@{-->}@/^1pc/[rr]^{/\!/_{\chi_\theta} T_\Theta} & \mcl{M}_0 & \mcl{M}_\theta \ar[l] \\
\mcl{V} \ar[u] \ar[r]^{/\!/SL_R} & \mcl{W} \ar[u] \ar[r]^{/\!/T_\Theta} \ar@{-->}@/_1pc/[rr]_{/\!/_{\chi_\theta} T_\Theta} & V/G \ar[u] & \mfk{M}_\theta \ar[l] \ar[u]
}\]
for any $\theta\in\Theta$ where the vertical arrows are all closed immersions of irreducible components. Using $\vphi$, we can give a description of the coordinate ring of $\mcl{W}=\mcl{V}/\!/SL_R$ similarly to the description of $\Cox(X)$ in (\ref{gens}). For this, let us consider the $\bar{\Theta}$-graded algebra $\bigoplus_{\theta\in\bar{\Theta}}\C[V]^{[G,G]} t_{\theta}$ where $t_{\theta}$ is a formal variable indicating the graded component for $\theta\in\bar{\Theta}$.

\begin{prop}\label{component V}
The algebra $\C[\mcl{W}]=\C[\mcl{V}]^{SL_R}$ is isomorphic to the graded subring $\mathbf{S}_\mcl{W}$ of $\bigoplus_{\theta\in\bar{\Theta}}\C[V]^{[G,G]} t_{\theta}$ generated by $\vphi(h)t_\theta$ for semi-invariants $h\in \C[\mcl{N}]^{SL_R}$ of weight $\theta$.
\end{prop}

\begin{proof}
The $T_\Theta$-equivariant morphism
$$V/[G,G]\times T_\Theta\to \mcl{N}/\!/SL_R$$
corresponding to the homomorphism $\C[\mcl{N}]^{SL_R}\to \bigoplus_{\theta\in \bar{\Theta}} \C[V]^{[G,G]} t_{\theta}$ is induced from the morphism $V\to \mcl{V} \subset \mcl{N}$ (and thus descends to the closed immersion $V/G\to \mcl{M}_0$). Then the claim follows since the spectrum of the integral domain $\mathbf{S}_\mcl{W}$ is isomorphic to the scheme-theoretic image of the equivariant morphism above, which is equal to $\mcl{V}/\!/SL_R=\mcl{W}$.
\end{proof}

From now on we assume $n=3$, i.e. $V=\C^3$. We fix any GIT-chamber $C\subset \Theta$ and let $X:=\mfk{M}_C$ be a crepant resolution of $V/G$. Similar to the construction of $\vphi$, we define a ring homomorphism
$$\vphi_C:\C[\mcl{N}]^{SL_R}\to \Cox(X)$$
so that, for a homogeneous $h\in \C[\mcl{N}]^{SL_R}$, the zero locus of $h$ in $\mcl{N}^{C\text{-ss}}$ descends to a divisor $D_h$ of $X$ which is defined by $\vphi_C(h)\in H^0(\mcl{O}_X(D_h))\subset \Cox(X)$. More explicitly, $\vphi_C(h)$ is given as
\begin{equation}\label{eq:vphi}
\vphi_C(h)=\vphi(h)\prod t_k^{-r_k\cdot \mrm{ord}_{E_k}(h)}
\end{equation}
where $\mrm{ord}_{E_k}(h)$ is the coefficient of the exceptional divisor $E_k$ in $D_h$. Note that $\vphi$ fits into the following commutative diagram}

\[
\begin{tikzcd}
\C[\mcl{N}]^{SL_R} \dar{\vphi_C} \rar[twoheadrightarrow] \arrow{drr}{\vphi} & \C[\mcl{V}]^{SL_R} \rar[]{\cong} & \mathbf{S}_\mcl{W} \dar[swap]{t_\theta\mapsto 1}\\
\Cox(X) \rar[hook] & \C[V]^{[G,G]}[t_1^{\pm1},\dots,t_m^{\pm1}] \rar[swap]{t_k\mapsto 1}& \C[V]^{[G,G]}
\end{tikzcd}
\]
where the upper right isomorphism comes from Proposition \ref{component V}. Since $\mcl{V}$ is the irreducible component of $\mcl{N}$ which dominates $\mfk{M}_C$, $\vphi_C$ factors through $\C[\mcl{V}]^{SL_R}$ and we thus obtain a morphism
$$\vphi_C^*: \mfk{X}\to \mcl{W}$$
which descends to the identification $X\to \mfk{M}_C$.

The ring homomorphism $\vphi$ also induces a homomorphism
\begin{equation}\label{psi}
\psi_C:\Theta\to\mrm{Pic}(X)_\R
\end{equation}
such that the degree of $\vphi_C(f)$ for a homogeneous element $f$ of weight $\theta\in\Theta$ is equal to $\psi_C(\theta)$. This is indeed well-defined since two semi-invariants on $\mcl{V}$ of the same degree differ by multiplication by a $G$-invariant rational function on $V$ and thus define the same class in $\mrm{Pic}(X)$. Later we will see that $\psi_C$ can be computed explicitly using an embedding of $\mcl{W}$ into a toric variety. See Remarks \ref{rem:vphi}(1) and examples in Subsection \ref{3.5}. We will also see in Lemma \ref{lem:commute} that $\psi_C$ may be defined in terms of tautological bundles on the moduli space $X=\mfk{M}_C$.

Now we are ready to define an embedding of $\mcl{W}$ into a toric variety. We first choose $Ab(G)^\vee$-homogeneous generators $f_1,\dots,f_\ell\in\C[V]^{[G,G]}$ such that the associated elements in (\ref{gens}) generate $\Cox(X)$. Let $\iota_\mfk{X}:\mfk{X}\to \A^{\ell+m}$ be the associated embedding constructed in the beginning of this subsection. Set $N_X=\Z^\ell\oplus\mrm{Pic}(X)^\vee$ and let $\sigma_\A\subset (N_X)_\R\cong\R^{\ell+m}$ be the cone defining the affine toric variety $\A^{\ell+m}$ with the $T_X$-action. We also set $N_\mcl{W}=\Z^\ell\oplus \bar{\Theta}^\vee$ and define $\sigma_\mcl{W}\subset (N_\mcl{W})_\R\cong \R^{\ell+s-1}$ as the cone generated by the set of cones $\{(\mrm{id}_{\R^\ell}\oplus \psi_C^*)(\sigma_\A)\}_C$ where $C$ runs through all GIT-chambers $C\subset\Theta$ and the choice of $f_i$'s is consistent regardless of the choice of $C$. Note that $\psi_C$ and its pullback $\psi_C^*:\mrm{Pic}(X)_\R^\vee\to \Theta^\vee$ depend on $C$.

Let $\bar{Y}_\mcl{W}$ be the affine normal toric variety defined by $\sigma_\mcl{W}$. We show that the normalization $\bar{\mcl{W}}$ of $\mcl{W}$ has an embedding into $\bar{Y}_\mcl{W}$. By definition, the coordinate ring of $\bar{Y}_\mcl{W}$ is the semigroup algebra $\C[S_\mcl{W}]$ for the semigroup $S_\mcl{W}=M_\mcl{W}\cap \sigma_\mcl{W}^\vee$ where $M_\mcl{W}= N_\mcl{W}^\vee \cong\Z^{\ell+s-1}$ and $\sigma_\mcl{W}^\vee$ is the dual cone
$$\{f\in (M_\mcl{W})_\R \mid f(v)\ge0\text{ for all }v\in \sigma_\mcl{W}\}$$ to $\sigma_\mcl{W}$.
Since $f_i$'s generate $\C[V]^{[G,G]}$, we can choose homogeneous generators $\{h_j\}_j$ of $\C[\mcl{W}]=\C[\mcl{V}]^{SL_R}$ such that each $\vphi(h_j)$ is of the form $f_1^{a_1}\cdots f_\ell^{a_\ell}$. This presentation may not be unique in general, but we fix it for each $\vphi(h_j)$. Then $\mcl{W}$ admits a closed immersion into a possibly non-normal affine toric variety
\[
Y_\mcl{W}=\mrm{Spec}\,\C[S'_\mcl{W}]
\]
defined by the semigroup $S'_\mcl{W}\subset M_\mcl{W}$ generated by $\{(a_1,\dots,a_\ell,\theta_j)\}_j$ where $\theta_j\in\bar{\Theta}$ is the weight of $h_j$. This also induces a closed immersion of $\bar{\mcl{W}}$ into the normalization of $Y_\mcl{W}$. Note that we have $S'_\mcl{W}\subset S_\mcl{W}$ and thus $\sigma_\mcl{W}$ is contained in the cone dual to $S'_\mcl{W}$.

\begin{lem}\label{integral closure}
We have the following commutative diagram
\begin{equation}\label{CD1}
\begin{CD}
\A^{\ell+m} @>>> \bar{Y}_\mcl{W} @>>> Y_\mcl{W}\\
@A \iota_\mfk{X} AA @AAA @AAA\\
\mfk{X}@>>> \bar{\mcl{W}} @>>> \mcl{W}\\ 
\end{CD}
\end{equation}
where the three vertical arrows are closed immersions and the two right-hand horizontal arrows are the normalization maps.
\end{lem}

\begin{proof}
Given the above, it suffices to show that the semigroup algebra $\C[S_\mcl{W}]$ coincides with the integral closure of $\C[S'_\mcl{W}]$.

For each $f\in S_\mcl{W}$, it defines a rational function $\bar{f}$ on $\bar{\mcl{W}}$ via the above embedding of $\bar{\mcl{W}}$. It suffices to show that $\bar{f}$ is regular on $\bar{\mcl{W}}$. Since $\bar{\mcl{W}}$ is normal it suffices to show that $\bar{f}$ is regular outside a locus of codimension greater than one. By the construction of $\sigma_\mcl{W}$, $\bar{f}$ is regular at general points of any divisor $D\subset \bar{\mcl{W}}$ which descends to a divisor on a resolution $\bar{\mcl{W}}/\!/_\theta T_\Theta$ of $V/G$ for some generic stability condition $\theta$. Then the claim follows by the lemma below.
\end{proof}

\begin{lem}
The complement of the open subset $\bigcup_{\theta\textup{:generic}} \mcl{W}^{\theta\text{-ss}}$ in $\mcl{W}$ is of codimension at least two.
\end{lem}

\begin{proof}
Assume that the complement contains an irreducible divisor $D\subset \mcl{W}$ in order to deduce a contradiction. Let $C\subset \Theta$ be the cone generated by the weights of semi-invariants of $\mcl{V}$ which do not vanish at general points of $D$. By the assumption, $C$ has positive codimension and, equivalently, general points of $D$ has a positive-dimensional stabilizer subgroup in $T_\Theta$. Then the image of $D$ under the quotient map $\mcl{W}^{\theta\text{-ss}}\to \mcl{W}/\!/_\theta T_\Theta$ for $\theta$ in the relative interior of $C$ has dimension greater than or equal to $n=3$. This is contrary to the irreducibility of $\mcl{W}$ since the $T_\Theta$-orbits of general points of $\mcl{W}$, which correspond to free $G$-orbits (see Subsection (\ref{2.1})), are closed and hence $\mcl{W}/\!/_\theta T_\Theta$ has dimension three for any (non-generic) $\theta$.
\end{proof}

For any $\theta\in C$, we obtain the following diagram of maps between quotients induced from (\ref{CD1}):
\begin{equation}\label{CD2}
\begin{CD}
Y @>>> \bar{Y}_\mcl{W}/\!/_{\chi_\theta} T_\Theta @>>> Y_\mcl{W}/\!/_{\chi_\theta} T_\Theta\\
@A \iota_X AA @AAA  @AAA\\
X @>>> \bar{\mcl{W}}/\!/_{\chi_\theta} T_\Theta @>>> \mcl{W}/\!/_{\chi_\theta} T_\Theta\\ 
\end{CD}
\end{equation}
where vertical arrows are again closed immersions and horizontal arrows are all isomorphisms. The results obtained so far are summarized as follows:

\begin{thm}\label{embedding}
The irreducible component $\mcl{W}$ of $\mcl{N}/\!/SL_R$ such that $\mcl{W}/\!/T_\Theta\cong \C^3/G$ admits a $T_\Theta$-equivariant closed immersion $\iota_\mcl{W}: \mcl{W}\to Y_\mcl{W}$ into a possibly non-normal toric variety whose torus contains $T_\Theta$, satisfying the following property:
\begin{itemize}
\item for any generic $\theta\in\Theta$, the induced morphism $\mcl{W}/\!/_{\chi_\theta} T_\Theta\to Y_\mcl{W}/\!/_{\chi_\theta} T_\Theta$ between the GIT-quotients is identified with a closed immersion $\iota_X: X\to Y$ of the crepant resolution $X:=\mfk{M}_\theta$ into a normal toric variety such that the restriction map $\iota_X^*: \mrm{Pic}(Y)_\R\to \mrm{Pic}(X)_\R$ is an isomorphism and $\mrm{Mov}(X)=\iota_X^*(\mrm{Mov}(Y))$.
\end{itemize}
\end{thm}

In the next subsection we will explain how the (semi)stable locus of $\bar{Y}_\mcl{W}$ for a given stability condition is determined in terms of the fan $\sigma_\mcl{W}$.

\begin{rems}\label{rem:vphi}
(1) Let $v_1,\dots,v_m\in \sigma_\mcl{W}$ be the primitive generators of the rays corresponding to the exceptional divisors $E_k$ of $X=\mfk{M}_C\to V/G$. Recall that the order of zeros (or poles) of a rational function $f\in M_\mcl{W}$ along the corresponding divisor to $v_k$ is given by $f(v_k)$, by toric geometry. Thus, the image of a homogeneous semi-invariant $h:=f_1^{a_1}\cdots f_\ell^{a_\ell}t_{\theta}\,(a_i\in \Z)$ under $\vphi_C$ is given by
\begin{equation}\label{phi(h)}
\tilde{f}_1^{a_1}\cdots \tilde{f}_\ell^{a_\ell}\prod_{k=1}^m t_k^{-\mbf{a}(v_k)r_k}
\end{equation}
(see Remark \ref{rem:cox(V/G)} and (\ref{eq:vphi})) where $\mbf{a}=(a_1,\dots,a_\ell,\theta)\in \Z^\ell\oplus \bar{\Theta}$ is identified with an element of $M_\mcl{W}$. In particular, $\psi_C(\theta)\in \mrm{Pic}(X)$ is given by the multi-degree of $t_k$'s in (\ref{phi(h)}). Note that $\psi_C$ depends only on $v_k$'s and therefore $\psi_{C'}$ is equal to $\psi_C$ as long as we take a GIT-chamber $C'$ inside the intersection of the orbit cones of $G$-constellations corresponding to general points of the divisors $E_1,\dots,E_m$.

(2) The homomorphism $\vphi_C:\C[\mcl{N}]^{SL_R}\to \Cox(X)$ is never surjective unless $G$ is trivial since each $t_k^{-r_k}\in\Cox(X)$ is clearly outside of the image. However, we will show that, for homogeneous generators $\{h_i\}$ of $\C[\mcl{N}]^{SL_R}$, the associated elements to the images $\vphi(h_i)$ generate $\Cox(X)$ (Proposition \ref{prop:generate}).
\end{rems}

\subsection{GIT-quotients of a toric variety by a subtorus-action}\label{3.4}

As we saw in the previous subsection, every crepant resolution of the form $\mfk{M}_\theta$ is embedded in a quotient of the toric variety $\bar{Y}_\mcl{W}$ (or $Y_\mcl{W}$) by the algebraic torus $T_\Theta$. In this subsection we show how to compute the $\theta$-semistable locus of $\bar{Y}_\mcl{W}$ for $\theta\in \Theta$ in terms of the fan of $\bar{Y}_\mcl{W}$. This also determines the semistable locus of $\bar{\mcl{W}}$ since we have
$$\bar{Y}_\mcl{W}^{\theta\text{-ss}}\cap \bar{\mcl{W}} =\bar{\mcl{W}}^{\theta\text{-ss}}$$
by the construction of $\bar{Y}_\mcl{W}$, If $\theta$ is generic, $\bar{\mcl{W}}^{\theta\text{-ss}}$ is also identified with $\mcl{W}^{\theta\text{-ss}}$ since the geometric quotient of $\mcl{W}^{\theta\text{-ss}}$ is smooth.

To consider quotients of $\bar{Y}_\mcl{W}$, first note that $T_\Theta$ is regarded as a subtorus of the {\it big torus} $T_\mcl{W}:=\mrm{Hom}_\Z(N_\mcl{W},\C^*)\subset Y_\mcl{W}$. Quotients of a toric variety by an action of a subtorus of the big torus have been studied in \cite{KSZ},\cite{AH},\cite{Hu} and so on. When we consider GIT-quotients (or more generally good quotients), the quotient map naturally induces a one-to-one correspondence between the sets of the maximal cones for the toric variety and its quotient (see e.g. \cite[Proposition 3.2]{AH}). To state this more precisely for our situation, recall that the toric variety $\bar{Y}_\mcl{W}$ is defined by the cone $\sigma_\mcl{W}$ in $(N_\mcl{W})_\R$. For any $\theta\in\Theta$, its semistable locus $\bar{Y}_\mcl{W}^{\theta\text{-ss}}$ is $T_\mcl{W}$-invariant since functions on $\bar{\mcl{W}}$ defining the semistable locus can be taken as homogeneous ones with respect to $T_\mcl{W}$. Thus, $\bar{Y}_\mcl{W}^{\theta\text{-ss}}$ is also a toric variety defined by a subfan $\Sigma_\theta$ of $\Sigma_\mcl{W}$ where $\Sigma_\mcl{W}$ is the fan consisting of all faces of $\sigma_\mcl{W}$. Then the GIT-quotient $Y_\theta:=\bar{Y}_\mcl{W}/\!/_\theta T_{\Theta}$ is a toric variety defined by a fan $\bar{\Sigma}_\theta$ in $\R^\ell$, and the quotient map $\bar{Y}_\mcl{W}^{\theta\text{-ss}}\to Y_\theta$ is a toric morphism induced by the projection
$$p: (N_\mcl{W})_\R\to(N_\mcl{W}/\bar{\Theta}^\vee)_\R=\R^\ell$$
satisfying the following condition:
\begin{itemize}\label{condition*}
\item[(\ref{condition*})] For any cone $\sigma\in\bar{\Sigma}_\theta$, the cone $p^{-1}(\sigma)\cap|\Sigma_\theta|$ belongs to $\Sigma_\theta$ where $|\Sigma_\theta|\subset (N_\mcl{W})_\R$ is the support of the fan $\Sigma_\theta$. Moreover, if $\theta$ is generic, $p$ induces an isomorphism of $\Sigma_\theta$ and $\bar{\Sigma}_\theta$.
\end{itemize}

Now we give an explicit method to compute $\Sigma_\theta$ and hence $\bar{Y}_\mcl{W}^{\theta\text{-ss}}$. To this end, it is sufficient to determine a general point of $p^{-1}(v)\cap|\Sigma_\theta|$ for each $v\in \R^\ell$ by the condition (\ref{condition*}) above. Let $q:(N_\mcl{W})_\R\cong \R^\ell\times\Theta^\vee \to \Theta^\vee$ be the second projection to the dual vector space of $\Theta\cong\R^{s-1}$ (while $p$ is regarded as the first projection). Then we set $P_v=q(p^{-1}(v)\cap|\Sigma_\mcl{W}|)\subset\Theta^\vee$, which is a rational convex polyhedron since $\sigma_\mcl{W}$ is rational and strongly convex. Note that $\theta$ is a linear function on $P_v$. The following result is regarded as a generalization of \cite[Theorem 7.2]{CMT} to the case where $G$ is non-abelian, see Example \ref{3.5.1}.

\begin{prop}\label{prop:minimize}
For any $\theta\in\Theta$, a face $F$ of $\sigma_\mcl{W}$ belongs to $\Sigma_\theta$ if and only if $q(F)\cap P_{p(\tilde{v})}$ equals the set of points of $P_{p(\tilde{v})}$ for which the function $\theta:P_{p(\tilde{v})}\to \R$ is minimized for some (and hence any) point $\tilde{v}$ of the relative interior $\mrm{relint}(F)$ of $F$.
\end{prop}

This result can be proven by using \cite[Proposition 2.7]{CM}, on which the proof of \cite[Theorem 7.2]{CMT} relies. Here we give another proof using the Hilbert-Mumford numerical criterion in order to understand the result in terms of one-parameter subgroups. The numerical criterion states that a point $x\in \bar{Y}_\mcl{W}$ is $\theta$-semistable if and only if we have $\lambda\cdot\theta\ge0$ for any 1-parameter subgroup (1-PS) $\lambda:\C^*\to T_\Theta$ such that $\lim_{t\to0}(\lambda(t)\cdot x)$ exists, where $\lambda\cdot\theta$ means the natural paring between a 1-PS and a character of $T_\Theta$.

To apply this criterion to toric varieties, let us recall that faces $F$ of $\sigma_\mcl{W}$ and $T_\mcl{W}$-orbits of $\mcl{W}$ are in one-to-one correspondence which assigns to $F$ the orbit $O(F)$ of $x_F:=\lim_{t\to0}\lambda^{\tilde{v}}(t)$ for the 1-PS $\lambda^{\tilde{v}}:\C^*\to T_\mcl{W}$ corresponding to $\tilde{v}\in\mrm{relint}(F)$. Moreover, the closure of $O(F)$ is equal to $\bigcup_{F\subset F'\in\Sigma_\mcl{W}} O(F')$, and we have $\lim_{t\to0}(\lambda^u(t)\cdot x_F)\in O(F')$ for any $u\in\mrm{relint}(F')$ where $\lambda^u$ is the  1-PS corresponding to $u$. See \cite[Proposition 3.2.2]{CLS} and its proof for these standard facts. We also observe that the limit $\lim_{t\to0}(\lambda^u(t)\cdot x_F)$ for $u\in N_\mcl{W}$ exists if and only if there exist (sufficiently small) $h>0$ and a face $F'$ of $\sigma_\mcl{W}$ containing $F$ such that $u_1+hu$ lies in $F'$ for some (and hence any) $u_1\in\mrm{relint}(F)$.

\vspace{3mm}

{\it Proof of Proposition \ref{prop:minimize}.} \\
According to the observation above, for a 1-PS $\lambda^w$ of $T_\Theta$ corresponding to a vector $w\in\bar{\Theta}^\vee$ and any point $x$ of the orbit $O(F)$, the limit $\lim_{t\to0} (\lambda^w(t)\cdot x)$ exists if and only if there exist $h>0$ and a face $F'$ containing $F$ such that $u_1+hu\in F$ for some $u_1\in\mrm{relint}(F)$.

Let us assume that $q(F)\cap P_{p(\tilde{v})}$ does not minimize $\theta$. Then there is another face $F'$ of $\sigma_\mcl{W}$ such that $q(F')\cap P_{p(\tilde{v})}$ minimizes $\theta$, and $(w_1-q(\tilde{v}))\cdot \theta<0$ for any $w_1\in \mrm{relint}(q(F')\cap P_{p(\tilde{v})})$. Since $q(\tilde{v})+1\cdot (w_1-q(\tilde{v}))=w_1$ is inside $P_{p(\tilde{v})}$, the convexity of $P_{p(\tilde{v})}$ implies that the limit for the 1-PS corresponding to $w_1$ exists and thus $F$ is not in $\Sigma_\theta$ by the numerical criterion.

Conversely, if $q(F)\cap P_{p(\tilde{v})}$ minimizes $\theta$, then $q(F)\cap P_{p(\tilde{v})}$ is a face of $P_{p(\tilde{v})}$ having $\{\theta=(\text{the minimal value})\}\subset \Theta^\vee$ as the supporting hyperplane. Therefore, $(w-q(\tilde{v}))\cdot \theta\ge0$ for any $w\in P_{p(\tilde{v})}$ and $F$ is in $\Sigma_\theta$ by the numerical criterion again.
\qed

\vspace{3mm}

\begin{rem}
Whenever we are given an affine toric variety and a subtorus of the big torus, the similar statement to the above proposition applies to compute the semistable locus of the given toric variety with respect to a character of the subtorus. A description of the fan of the GIT quotient is also given in \cite[\S3]{CM}.
\end{rem}

\subsection{Examples}\label{3.5}
In this subsection we demonstrate the construction of the toric embedding and the homomorphism $\vphi_C$ in the previous subsection via concrete examples. 

\subsubsection{Abelian cases}\label{3.5.1}
We first consider the case where $G\subset SL_3(\C)$ is a finite abelian subgroup. In this case all irreducible representations of $G$ are 1-dimensional, which particularly implies that $SL_R$ is trivial and hence $\mcl{V}=\mcl{W}$. Then there is a natural choice of generators of the Cox ring of a crepant resolution $X$ and generators of the coordinate ring of $\mcl{V}$ as follows.

Once we fix a basis of $V=\C^3$ such that $G$ is a diagonal group, we obtain generators of $\C[V]^{[G,G]}=\C[V]$ as the standard basis $x_1,x_2,x_3$ of $V^*$ and also standard generators of $\Cox(X)$ (cf. \cite[Proposition 3.5]{Y2}. As explained in \cite{CMT} (see also \cite{Y2}), the coordinate ring of $\mcl{V}$ is generated by $3r$ variables $\{x_{i,j}\}_{0\le i\le r-1,\,1\le j\le 3}$ with certain relations where $r$ is the order of $G$. In terms of the quiver representation, the variable $x_{i,j}$ is a homogeneous semi-invariant corresponding to an arrow $\rho_i$ to $\rho_i\otimes \chi_j$ of the McKay quiver of $G$ where $\mrm{Irr}(G)=\{\rho_0,\dots,\rho_{r-1}\}$ and we regard the character $\chi_j:G\to\C^*$ defined by $g\cdot x_j=\chi_j(g)x_j$ as a 1-dimensional representation of $G$. Then $\mcl{V}$ is already a possibly non-normal toric variety \cite[Theorem 3.10]{CMT}, and the ambient toric variety $Y_\mcl{V}$ constructed from these generators is nothing but $\mcl{V}$ itself.

For each orbit $O\subset \mcl{V}$ by the big torus, it has the {\it distinguished point} $v\in O$ satisfying $x_{i,j}(v)=1$ for all $i,j$ such that $x_{i,j}(O)\ne\{0\}$. For any $\theta$, the fan $\Sigma_\theta$ of the $\theta$-semistable locus of $\mcl{V}$ is computed using Proposition \ref{prop:minimize}, and one obtains the distinguished point corresponding to a cone $\sigma\in \Sigma_\theta$ by setting
$$x_{i,j}(v)=\begin{cases}1&\text{ if }\sigma\subset H_{i,j}\\0&\text{ if }\sigma\not\subset H_{i,j}\end{cases}$$
where $H_{i,j}\subset (N_\mcl{V})_\R$ is the hyperplane defined by $x_{i,j}\in M_\mcl{V}$. This is essentially the same procedure as the one given in \cite[Theorem 7.2]{CMT}, which enables us to obtain the distinguished point corresponding to the cone $\sigma_w\in \Sigma_\theta$ such that $p(\sigma_w)$ contains a given $w\in\R^3$ in its relative interior. Note that one can also consider distinguished points of $\bar{Y}_ \mcl{W}$ for non-abelian $G$ but these points usually do not lie inside $\bar{\mcl{W}}$.

In \cite{Y2}, an explicit description of the homomorphism $\vphi_C:\C[\mcl{V}]\to\Cox(X)$ is given in terms of the notion of a {\it $G$-nat family} introduced by Logvineko \cite{Lo}. Recall that a $G$-nat family on a crepant resolution $X\to \C^3/G$ is a family of $G$-constellations which extends to $X$ the family over the open subset of $\C^3/G$ parametrizing the free $G$-orbits. To give such a family is equivalent to choosing rays of $\sigma_\mcl{V}$ such that their images under $p:(N_\mcl{V})_\R\to\R^3$ are the rays of the toric fan of $X$. The main result of \cite{Y2} for $n=3$ states that $X$ is realized as a fine moduli space of $G$-constellations (not-necessarily of the form $\mfk{M}_\theta$) if and only if  $X$ admits a $G$-nat family $\mcl{F}$ such that the maximal cones of the ``associated" subfan $\Sigma_{X,\mcl{F}}$ of $\Sigma_\mcl{V}$ has the expected dimension $3$ (see \cite[\S4]{Y2} for details). Note that the results in \cite{Y2} work for abelian $G\subset SL_n(\C)$ with arbitrary $n$.

\begin{rem}
If we add redundant homogeneous generators (like $x_j^2+x_{j'}$ if $\chi_j^2=\chi_{j'}$) to the generators $x_1,\dots,x_n$ of $\C[V]$, we obtain a strictly bigger toric variety $\bar{Y}_\mcl{V}$ than $\bar{\mcl{V}}$.
\end{rem}

\begin{ex}\label{abel2}(Continuation of Example \ref{abel1})\\
We consider again the abelian $G$ treated in Example \ref{abel1}. In order to describe the normalization $\bar{\mcl{V}}$ of $\mcl{V}$ concretely, we fix coordinates $\theta_0,\dots,\theta_9$ of $\theta\in\Theta$ so that $\theta_i=\theta(\rho^{\otimes i})\in\R$ where $\rho:G\to \C^*$ is the one-dimensional representation defined by $\rho(g)=\zeta_{10}$. We identify $\Theta$ with $\R^9$ via $\theta\mapsto(\theta_1,\dots,\theta_9)$ from now on. Note that $\sum_{i=0}^9 \theta_i=0$. Then we see that the cone $\sigma_\mcl{V}$ admits 65 rays, whose images under $p:\R^3\times\R^9\to \R^3$ are equal to the rays $\R_{\ge0}e_j\,(j=1,2,3)$ and $\R_{\ge0}v_k\,(k=1,\dots,5)$, as expected.

In order to determine which rays of $\sigma_\mcl{V}$ correspond to exceptional divisors of $X_1$, recall that the $G$-Hilbert scheme is realized as $\mcl{V}/\!/_{\theta_+} T_\Theta$ for any $\theta_+\in\Theta$ satisfying $(\theta_+)_i>0$ for all $i>0$. Note that GIT-chamber $C_1$ for the $G$-Hilbert scheme contains all such $\theta_+$. Using Proposition \ref{prop:minimize}, we see that the rays generated by the following vectors in $\R^3\times\R^9$ are the rays in $\Sigma_{\theta_+}$:
{\small $$\begin{alignedat}{5}
w_1 &= (1,3,6,-1, -& &2, -3,-& &4,-5,-& &6,-7,-& &8,-9)\\
w_2 &= (1,2,3,-1, -& &2, -3,-& &4,-5,-& &1,-2,-& &3,-4)\\
w_3 &= (2,1,2,-2, -& &4, -1,-& &3,-5,-& &2,-4,-& &6,-3)\\
w_4 &= (1,1,0,-1,   & &0, -1, & &0,-1,  & &0,-1,  & &0,-1)\\
w_5 &= (7,1,2,-7, -& &4, -1,-& &8,-5,-& &2,-9,-& &6,-3).
\end{alignedat}$$}
Alternatively, one can determine the $w_k$'s by using the fact that the universal family for the $G$-Hilbert scheme is the {\it maximal shift family} in the sense of \cite[\S3.5]{Lo}. More explicitly, the $(3+j)$-th entry of $w_k$ is given as
\begin{equation}\label{eq:maximal}
-\mrm{min}\{a d_{x,k}+b d_{y,k}+c d_{z,k} \mid a,b,c \in \Z_{\ge0} \text{ s.t. } g\cdot x^a y^b z^c=\rho_j(g)x^a y^b z^c,\; \forall g\in G\}
\end{equation}
where $d_{x,k}, d_{y,k}, d_{z,k}\in \Z_{\ge0}$ are the exponents of $t_k$ for the generators  in (\ref{ex:generator}) associated to $x,y,z$ respectively.

Now we can explicitly describe $\vphi_{C_1}$ as explained in Remarks \ref{rem:vphi}(1). If we embed $\mrm{Pic}(X_1)$ into $\Z^5$ as we have been doing so far (cf. (\ref{ex:generator})), then the map $(\psi_{C_1})_\R:\R^9\to\R^5$ is given by the matrix
{\small $$\begin{pmatrix}1&2&3&4&5&6&7&8&9\\
1&2&3&4&5&1&2&3&4\\
2&4&1&3&5&2&4&6&3\\
1&0&1&0&1&0&1&0&1\\
7&4&1&8&5&2&9&6&3
\end{pmatrix}.$$}
Note that the $(k,j)$-entry of this matrix equals the negative of the $(3+j)$-th entry of $w_k$. Indeed, a monomial $x^a y^b z^c$ which attains the minimum of (\ref{eq:maximal}) naturally defines an semi-invariant which has weight $(0, \dots, 0, \stackrel{j\text{-th}}{1}, 0,\dots, 0)\in \R^9$ and does not vanish along $E_k$. One can check that we have $(\psi_{C_1})_\R(\overline{C}_1)=\mrm{Nef}(X_1)$.

Now we can explicitly confirm that every crepant resolution of $\C^3/G$ is realized as $\mfk{M}_\theta$ for some generic $\theta$, which is guaranteed by the main result of \cite{CI}. As observed in Remarks \ref{rem:vphi}(1), the map $\psi_C$ for a GIT-chamber $C$ coincides with $\psi_{C_1}$ as long as $C$ is contained in the intersection $C_+:=\bigcap_{k=1}^5 C_{w_k}$ of the orbit cones
$$C_{w_k}=\mrm{Cone}(\{\mrm{wt}(x_{i,j})\in\Theta\mid w_k\subset\{x_{i,j}=0\}\}).$$
One sees that we have $(\psi_{C_1})_\R(C_+)=\bigcup_{i=1}^5 \mrm{Nef}(X_i)$, which implies that the interior of
$$\overline{C}_i:=(\psi_{C_1})_\R^{-1}(\mrm{Nef}(X_i))\cap C_+ \quad(i=1,\dots,5)$$
is a GIT-chamber $C_i$ whose associated moduli space equals $X_i$. The resolution $X_6$ is also achieved as a moduli space by crossing a boundary of $C_+$. More precisely, one can check that the hyperplane $\{\theta_0+\theta_1+\theta_3=0\}$ forms a facet of both $C_+$ and $C_3$ and that the adjacent chamber $C'_3$ to $C_3$ regarding this facet corresponds to  the ray of $\sigma_\mcl{V}$ generated by
$$w'_3=(2,1,2,-2,1,-1,2,0,3,1,-1,2).$$
Then the counterpart $C'_+$ of $C_+$ obtained by replacing $C_3$ with $C'_3$ satisfies
$$(\psi_{C'_3})_\R(C'_+)\supset\mrm{Nef}(X_6).$$
Therefore, we obtain a GIT-chamber $C_6\subset C'_+$ such that $\mfk{M}_{C_6}\cong X_6$.
\end{ex}

\subsubsection{$D_5$-singularity, continuation of Example \ref{ex:relations}} \label{3.5.2}

We next consider the case of $D_5$-singularity again (see Example \ref{ex:relations}). Note that the construction of embeddings in Subsection \ref{3.3} is valid in dimension two since the moduli space for any generic $\theta$ is a crepant (or the minimal) resolution of $\C^2/G$ as well.

We use $x,y$ as the standard coordinates of $V=\C^2$ again. In terms of representations of the quiver $Q_G$ of Figure \ref{quiver}, the $G$-constellation family over $\C^2$ are given by the following generic matrices:
{\small
$$\begin{aligned}
&A_0=\begin{pmatrix}x\\y\end{pmatrix},B_0=\begin{pmatrix}-y&x\end{pmatrix},
A_1=\begin{pmatrix}-i y\\x\end{pmatrix},B_1=\begin{pmatrix}x&i y\end{pmatrix},
A_2=\begin{pmatrix}x\\-y\end{pmatrix},B_2=\begin{pmatrix}y&x\end{pmatrix},\\
&A_3=\begin{pmatrix}i y\\x\end{pmatrix},B_3=\begin{pmatrix}x&-i y\end{pmatrix},
C=\begin{pmatrix}x&0\\0&y\end{pmatrix},D=\begin{pmatrix}-y&0\\0&x\end{pmatrix}
\end{aligned}$$
}for the choice of basis of the vector spaces $\mrm{Hom}_{\C G}(V^*\otimes_\C \rho,\rho')$ in Example \ref{ex:relations}. Here we regard the arrows of $Q_G$ as the corresponding matrices by abuse of notation. Note that these matrices satisfy the condition (\ref{path}). We take $Ab(G)^\vee$-homogeneous elements
$$f_1=x^3+iy^3,\quad f_2=xy, \quad f_3=x^3-iy^3$$
so that $\C[V]^{[G,G]}=\C[f_1,f_2,f_3]$. Note that $g\in G$ acts on $f_k$ via the character $\rho_k$. As is shown in \cite[\S6]{D} and \cite[Example 1]{Y1}, these $f_k$'s give rise to generators of the Cox ring of the minimal resolution $X$ of $V/G$. We will consider the embedding of $\mcl{W}$ into a toric variety with respect to these $f_k$'s. 

Let us consider the semi-invariants obtained as the determinants of the following 42 square (block) matrices:
{\small
$$\begin{aligned}
&C,D, B_2A_0,B_0A_2,B_3A_1,B_1A_3, \begin{pmatrix}A_0&A_2\end{pmatrix},\begin{pmatrix}B_0\\B_2\end{pmatrix}, \begin{pmatrix}A_1&A_3\end{pmatrix},\begin{pmatrix}B_1\\B_3\end{pmatrix}, B_1CA_0, B_3CA_0, \\
&B_1CA_2, B_3CA_2, B_0DA_1, B_2DA_1, B_0DA_3, B_0DA_3, \begin{pmatrix}CA_0&A_1\end{pmatrix},\begin{pmatrix}CA_0&A_3\end{pmatrix},\\
&\begin{pmatrix}CA_2&A_1\end{pmatrix},\begin{pmatrix}CA_2&A_3\end{pmatrix},\begin{pmatrix}A_0&DA_1\end{pmatrix},\begin{pmatrix}A_0&DA_3\end{pmatrix},\begin{pmatrix}A_2&DA_1\end{pmatrix},\begin{pmatrix}A_2&DA_3\end{pmatrix},\\
&\begin{pmatrix}B_0\\B_1C\end{pmatrix}, \begin{pmatrix}B_0\\B_3C\end{pmatrix},\begin{pmatrix}B_2\\B_1C\end{pmatrix},\begin{pmatrix}B_2\\B_3C\end{pmatrix},\begin{pmatrix}B_0D\\B_1\end{pmatrix},\begin{pmatrix}B_0D\\B_3\end{pmatrix},\begin{pmatrix}B_2D\\B_1\end{pmatrix}, \begin{pmatrix}B_2D\\B_3\end{pmatrix},\begin{pmatrix}A_0&D\\0&B_1\end{pmatrix}, \\
&\begin{pmatrix}A_0&D\\0&B_3\end{pmatrix}, \begin{pmatrix}A_2&D\\0&B_1\end{pmatrix}, \begin{pmatrix}A_2&D\\0&B_3\end{pmatrix},\begin{pmatrix}B_0&0\\C&A_1\end{pmatrix}, \begin{pmatrix}B_0&0\\C&A_1\end{pmatrix}, \begin{pmatrix}B_2&0\\C&A_3\end{pmatrix}, \begin{pmatrix}B_2&0\\C&A_3\end{pmatrix}.
\end{aligned}$$
} These semi-invariants are chosen so that Condition (\ref{*}) below is satisfied. (See Remark \ref{rem:candidate}.) We will see that these semi-invariants give a generating system of $\C[\mcl{W}]=\C[\mcl{V}]^{SL_R}$. Although it might be possible to check this directly (with help of a computer) by using algorithms such as the one in \cite[4.1]{DeK}, we adopt another approach.

In order to construct $Y_\mcl{W}$ explicitly, we use coordinates $\theta_0,\dots,\theta_5$ of $\theta\in\Theta$ such that $\theta_i=\theta(\rho_i)$ for $i=0,1,2,3$ and $\theta_i=\theta(V_{i-3})$ for $i=4,5$. We then identify $\Theta$ with $\R^5$ via $\theta\mapsto (\theta_1,\dots,\theta_5)$. Note that we have $\theta_0+\theta_1+\theta_2+\theta_3+2\theta_4+2\theta_5=0$ in this case. Each of the 42 semi-invariants above is mapped to one of $f_i$'s under $\vphi$ up to constant multiplication, and we construct an embedding $\mcl{W}\hookrightarrow Y_\mcl{W}$ with respect to $f_i$'s. For example, the associated vector in $M_\mcl{W}\subset\Z^{3+5}$ to the semi-invariant for the matrix $C$  (or the corresponding arrow in $Q_G$) is given as $(0,1,0,0,0,0,-1,1)$ since $\det C=f_2$ and the weight $\theta\in\Theta$ of this semi-invariant satisfies $\theta_1=\theta_2=\theta_3=0,\,\theta_4=-1,\,\theta_5=1$ (see (\ref{eq:weight})). Similarly, the associated vector to $\begin{pmatrix}B_2&0\\C&A_3\end{pmatrix}$ is given as $(1,0,0,0,1,-1,-1,1)$.

Let $\sigma'\subset \R^{3+5}$ be the dual cone to the cone spanned by the vectors in $M_\mcl{W}$ associated to the above semi-invariants. We will show that $\sigma'$ equals the cone $\sigma_\mcl{W}\subset \R^8$ defining $\bar{Y}_\mcl{W}$. Let $\mcl{W}'$ be the spectrum of the subalgebra of $\C[\mcl{W}]$ generated by the 42 semi-invariants above. Then we obtain  an embedding of $\mcl{W}'$ into a toric variety $Y'_\mcl{W}$ whose normalization is defined by $\sigma'$, similar to the construction of $\mcl{W}\hookrightarrow Y_\mcl{W}$. In fact one can check that $Y'_\mcl{W}$ is normal by comparing $\sigma'^\vee\cap M_\mcl{W}$ and the semi-group generated by the elements of $M_\mcl{W}$ associated to the 42 semi-invariants. This computation can be done e.g. by using the package ``Normaliz" of Macaulay2 \cite{GS}.

We see, by calculation, that $\sigma'$ admits 165 rays. The primitive generators in $\Z^8$ of these rays include the first three coordinate vectors $\mbf{e}_1,\mbf{e}_2,\mbf{e}_3$, and the other 162 vectors are divided into five types by looking at their first three coordinates:
\begin{equation}\label{rays}
v_1=(1,2,1),\; v_2=(1,1,1),\; v_3=(3,2,3),\; v_4=(3,2,5), \text{ or }v_5=(5,2,3).
\end{equation}
These five types correspond to the prime exceptional divisors $E_1,\dots,E_5$ of the minimal resolution $X$ defined by the conjugacy classes of $g_1,g_1^2,g_1^3,g_2,g_1g_2$ respectively (cf. \cite[Example 1]{Y1}). More precisely, the rays for $E_1,\dots,E_5$ in the fan of the ambient toric variety $Y$ of $X$ (with respect to $f_i$'s) are generated by the vectors in (\ref{rays}). Note that $\sigma'$ particularly satisfies the following condition:
\begin{equation}\label{*}
\text{the polyhedral cone }p^{-1}(\R_{\ge0}v_k)\cap\sigma'\text{ is generated by rays of }\sigma'\text{ for each }k.
\end{equation}
This will be important when one computes a generating system of semi-invariants efficiently (see Remark \ref{rem:candidate}).

The fact that the rays of $\sigma'$ are projected to those of the fan of $X$ implies that $Y'_\mcl{W}/\!/_\theta T_\Theta$ is isomorphic to the toric variety $Y$ for any generic $\theta\in\Theta$. Therefore, we have $\sigma_\mcl{W}=\sigma'$, and the natural toric morphism $\bar{Y}_\mcl{W}\to{Y}'_\mcl{W}$ is an isomorphism. This also implies that the 42 semi-invariants generate $\C[\mcl{W}]$.
 
Since the $G$-Hilbert scheme is given by a stability condition $\theta_+$ satisfying $(\theta_+)_i>0$ for all $i>0$, we see that generators of the rays of $\sigma_\mcl{W}$ corresponding to the exceptional divisors $E_1,\dots,E_5$ of the $G$-Hilbert scheme are given by
{\small $$\begin{aligned}
&w_1=(1,2,1,-1,-2,-1,-2,-2),\quad &w_2=(1,1,1,-1,-1,-1,-2,-2), \\
&w_3=(3,2,3,-3,-2,-3,-4,-6),\quad &w_4=(3,2,5,-3,-2,-5,-4,-6),\\
&w_5=(5,2,3,-5,-2,-3,-4,-6)\end{aligned}$$
}respectively.  From these data, one can also explicitly compute the map $(\psi_C)_\R:\R^5\to\R^5$ associated to the GIT-chamber $C$ containing $\theta_+$.

\begin{rem}
In \cite[\S2]{Le}, the $G$-Hilbert schemes $X$ for binary dihedral groups are studied, and in particular the author introduced distinguished $G$-clusters which give rise to open affine covering of $X$. For example, the distinguished $G$-cluster $F$ of {\it type B} is defined by an homogeneous ideal of $\C[x,y]$ and is presented as the following $\C$-vector subspace:
$$\C\oplus(\C x\oplus\C y)\oplus\C xy\oplus(\C x^2\oplus\C y^2)\oplus\C x^3 \oplus \C y^3\oplus(\C x^4\oplus\C y^4)\oplus(\C x^5\oplus\C y^5).$$
Its structure as a $\C[V]$-module is also presented as in the following diagram:
\[\xymatrix@=24pt{
\rho_0 \ar[d]_-{} 
&
&\rho_1 \ar[dr]_{}
&
&\\
V_1 \ar[d]_{} \ar[r]_{}
&V_2 \ar[ur]_{} \ar[dr]_{}
& 
&V_2 \ar[r]_{} 
&V_1\\
\rho_2
&
&\rho_3 \ar[ur]_-{} 
&
&
}\]
where each arrow $\rho\to \rho'$ represents a nontrivial action of $V^*$ on $\rho$ into $\rho'$ similarly to the McKay quiver. One sees that, among the 42 semi-invariants, the ones $h_1,\dots,h_9$ for the matrices
{\footnotesize$$B_1CA_0, B_2A_0,B_3CA_0,\begin{pmatrix}CA_0&A_1\end{pmatrix},\begin{pmatrix}CA_0&A_3\end{pmatrix},\begin{pmatrix}A_0&DA_1\end{pmatrix},\begin{pmatrix}A_0&DA_3\end{pmatrix},\begin{pmatrix}A_0&D\\0&B_1\end{pmatrix}, \begin{pmatrix}A_0&D\\0&B_3\end{pmatrix}$$
}do not vanish at $F$, and thus the face of $\sigma_\mcl{W}$ corresponding to $F$ is given as the intersection $\sigma_\mcl{W}\cap\{w_{h_1}=\dots =w_{h_9}=0\}$ where $w_{h_i}\in M_\mcl{W}$ is the associated element to $h_i$. This is equal to $\mrm{Cone}(w_1,w_2 )$, and thus $F$ corresponds to the the intersection point of the two curves $E_1$ and $E_2$.
\end{rem}

\subsubsection{Trihedral group of order 21}\label{3.5.3}
Here we deal with the following three-dimensional example:
{\small $$G=\left\langle
g_1=\begin{pmatrix}
\zeta_7&0&0\\
0&\zeta_7^2&0\\
0&0&\zeta_7^4
\end{pmatrix},
g_2=\begin{pmatrix}
0&0&1\\
1&0&0\\
0&1&0\end{pmatrix}
\right\rangle\subset SL_3(\C)$$
}where $\zeta_7$ is the seventh root of unity. Then $\mrm{Irr}(G)$ consists of three 1-dimensional representations $\rho_k:G\to\C^*\,(k=0,1,2)$ defined by $\rho_k(g_1)=1,\,\rho_k(g_2)=\omega^k$ where $\omega$ is the third root of unity, and two 3-dimensional representations $V_1,V_2$, the first one is the inclusion $G\subset SL_3(\C)$ presented above and the other is presented as
{\small $$g_1\mapsto\begin{pmatrix}
\zeta_7^3&0&0\\
0&\zeta_7^6&0\\
0&0&\zeta_7^5
\end{pmatrix},\;
g_2\mapsto\begin{pmatrix}
0&0&1\\
1&0&0\\
0&1&0
\end{pmatrix}.$$
}We again fix the basis obtained in this way for the subsequent computations.

The McKay quiver $Q_G$ for this case is as follows:

{\large
\[\xymatrix@=60pt{
\rho_0 \ar[d]_-{A_0} 
& \rho_1 \ar[ld]_(.6){A_1}
&\rho_2 \ar[lld]_(.15){A_2}\\
V_1\ar@(l,d)[]_{L_1} \ar[rr]^(.5){C_1} \ar@/^1pc/[rr]^{C_2}
&
&V_2\ar@(r,d)[]^{L_2} \ar[lu]_(.35){B_1} \ar[u]_-{B_2} \ar@/^1pc/[ll]^{D} \ar[llu]_(.85){B_0}
}\]
}

With the standard coordinates $x,y,z$ of $V=\C^3$, the generic matrices are given as
{\small $$\begin{aligned}
&A_0=\begin{pmatrix}x\\y\\z\end{pmatrix},B_0=\begin{pmatrix}z&x&y\end{pmatrix},
A_1=\begin{pmatrix}x\\\omega y\\\omega^2 z\end{pmatrix},B_1=\begin{pmatrix}\omega z&x&\omega^2 y\end{pmatrix},\\
&A_2=\begin{pmatrix}x\\\omega^2 y\\\omega z\end{pmatrix}, B_2=\begin{pmatrix}\omega^2 z&x&\omega y\end{pmatrix}, C_1=\begin{pmatrix}0&x&0\\0&0&y\\z&0&0\end{pmatrix}, C_2=\begin{pmatrix}y&0&0\\0&z&0\\0&0&x\end{pmatrix}, \\
&D=\begin{pmatrix}0&y&0\\0&0&z\\x&0&0\end{pmatrix}, L_1=\begin{pmatrix}0&0&z\\x&0&0\\0&y&0\end{pmatrix},L_2=\begin{pmatrix}0&z&0\\0&0&x\\y&0&0\end{pmatrix}
\end{aligned}$$
}for some choice of basis $\mrm{Hom}_{\C G}(\rho\otimes_\C V^*,\rho')$ with $\rho,\rho'\in\mrm{Irr}(G)$.

In this case, it seems hard to obtain generators of semi-invariants directly using the algorithm in \cite[Ch. 4]{DeK} since the input data are large. Instead we give a candidate of a system of generators and then confirm the correctness as we did in the previous example.

As computed in \cite[\S5.1]{DG}, there are $Ab(G)^\vee$-homogeneous polynomials $f_1,\dots,f_{13}\in\C[V]^{[G,G]}$ such that the associated elements
$$\tilde{f}_1,\dots,\tilde{f}_{13}, t_1^{-7}, t_2^{-3}, t_3^{-3}\in \C[V]^{[G,G]}[t_1^{\pm1}, t_2^{\pm1}, t_3^{\pm1}]$$
generate the Cox ring of a crepant  resolution $X$ of $V/G$, where we set the variables $t_1, t_2, t_3$ so that they correspond to junior elements $g_1, g_2, g_2^2$ respectively (see Subsection \ref{3.1}).

Similarly to the previous examples, we realize the cone $\sigma_\mcl{W}$ inside $\R^{13+5-1}$ by ordering the elements of $\mrm{Irr}(G)$ as $\rho_0,\rho_1,\rho_2,V_1,V_2$. Regarding the embedding $\iota_X:X\hookrightarrow Y$ with respect to $f_1,\dots, f_{13}$, the generators of the rays in $\R^{13}$ corresponding to the exceptional divisors $E_1, E_2, E_3$ are given as
$$\begin{aligned}
&\mbf{v}_1=(1,1,1,1,1,1,1,1,1,1,1,1,1),\\
&\mbf{v}_2=(0,0,2,1,0,2,1,3,2,1,3,2,4),\\
&\mbf{v}_3=(0,0,1,2,0,1,2,3,1,2,3,4,2)
\end{aligned}$$
respectively.
Then we define $\sigma'\subset\R^{17}$ as the cone generated by the coordinate vectors $\mbf{e}_1, \dots, \mbf{e}_{13}$ and the following 16 vectors:
$$\begin{aligned}
&\mbf{v}_{1,1}=(\mbf{v}_1, -1,-1,-1,-2),\; \mbf{v}_{1,2}=(\mbf{v}_1, 0,-1,-1,-2),\;\mbf{v}_{1,3}=(\mbf{v}_1, -1,0,-1,-2),\\
&\mbf{v}_{1,4}=(\mbf{v}_1, 0,0,-1,-2),\;\mbf{v}_{1,5}=(\mbf{v}_1, 0,0,-1,0),\;\mbf{v}_{1,6}=(\mbf{v}_1, 0,0,0,1),\\
&\mbf{v}_{1,7}=(\mbf{v}_1, 0,0,2,1),\; \mbf{v}_{1,8}=(\mbf{v}_1, 0,1,2,1),\;\mbf{v}_{1,9}=(\mbf{v}_1, 1,0,2,1),\; \mbf{v}_{1,10}=(\mbf{v}_1, 1,1,2,1),\\
&\mbf{v}_{2,1}=(\mbf{v}_2,-1,-2,-3,-3),\;\mbf{v}_{2,2}=(\mbf{v}_2,-1,1,0,0),\;\mbf{v}_{2,3}=(\mbf{v}_2,2,1,3,3),\\
&\mbf{v}_{3,1}=(\mbf{v}_3,-2,-1,-3,-3),\;\mbf{v}_{3,2}=(\mbf{v}_3,1,-1,0,0),\;\mbf{v}_{3,3}=(\mbf{v}_3,1,2,3,3).
\end{aligned}$$

By computing $\sigma'$ explicitly (with help of Macaulay2 \cite{GS}), one sees that it has 400 facets and accordingly we obtain 400 primitive elements of $M_\mcl{W}$ defining these facets. In fact we can take determinantal semi-invariants corresponding to these primitive elements. One can also check that the toric variety $Y'_\mcl{W}$ constructed from these elements is already normal again. Since the rays of $\sigma'$ are projected to those of the fan of $X$, we have that $\sigma'=\sigma_\mcl{W}$ again and that the associated 400 semi-invariants generate $\C[\mcl{W}]$. 

\begin{rem}\label{rem:candidate}
In general, a candidate $\sigma'$ of $\sigma_\mcl{W}$ can be given by computing determinantal semi-invariants $h$ (from ones with $\vphi(h)$ having smaller degrees) until the dual cone
$$\sigma'':=\bigcap_h\{v\in\R^\ell\oplus \Theta^\vee\mid w_h(v)\ge0\}$$
satisfies the condition (\ref{*}). This may be achieved by computing relatively small number of semi-invariants (compared to the number of generators of $\C[\mcl{W}]$). Then the candidate $\sigma'$ is defined as the cone generated by the rays of $\sigma''$ which are projected to the rays of the fan of $X$.

 In the above case, such $\sigma''$ is obtained from determinantal semi-invariants coming from $3\times3$-matrices only, while the 400 generators of $\C[\mcl{W}]$ include ones coming from $6\times6$-matrices. One may also verify that $\sigma_\mcl{W}=\sigma'$ by checking the existence of $G$-constellations corresponding to rays of $\sigma'$, as will  be done below.
\end{rem}

Once we know the cone $\sigma_\mcl{W}$, one can investigate the structures of $G$-constellations presented by faces of $\sigma_\mcl{W}$. As a demonstration, we take the ray generated by $\mbf{v}_{1,1}$. If we take the generators $\mbf{v}_{1,1},\dots,\mbf{v}_{1,10}$ of $\sigma'_\mcl{W}$ whose initial coordinates vector is $\mbf{v}_1$, the orbit cone $C_{\mbf{v}_{1,j}}\subset\Theta$ of the corresponding divisor to $\mbf{v}_{1,j}$ is given as the cone
$$\{\theta\in\Theta \mid q(\mbf{v}_{1,j})\cdot\theta=\underset{j=1,\dots,10}{\mrm{min}}q(\mbf{v}_{1,j})\cdot\theta\}$$
by Proposition \ref{prop:minimize}. Note that $\{C_{\mbf{v}_{1,j}}\}_{j=1,\dots,10}$ give a subdivision of $\Theta$ into polyhedral cones. The orbit cone $C_{\mbf{v}_{1,1}}\subset\Theta$ is defined by the following 5 inequalities:
$$\theta_1\ge0,\;\theta_2\ge0,\;\theta_1+\theta_2+2\theta_4\ge0,\;\theta_1+\theta_2+\theta_3+3\theta_4\ge0,\;\theta_1+\theta_2 +3\theta_3  +3\theta_4\ge0.$$
Notice that the coefficients of these inequalities are obtained as the direction vectors $\overrightarrow{\mbf{v}_{1,1}\mbf{v}_{1,j}}$ for $j$ such that the line segment $\overline{\mbf{v}_{1,1}\mbf{v}_{1,j}}$ forms an edge of the convex hull $\mrm{Conv}(\mbf{v}_{1,1},\dots,\mbf{v}_{1,10})$.  One can infer the structure of the $G$-constellation $F$ of a general point of $E_1$ from these inequalities. In this case the structure of $F$ is presented as the following diagram for some decomposition $F=\rho_0\oplus\rho_1\oplus\rho_2\oplus V_1^{\oplus3}\oplus V_2^{\oplus3}$:

\[\xymatrix@=48pt{
\rho_0 \ar[r]_-{} 
&V_1 \ar[r]_{} \ar[d]_{}
&V_2 \ar[d]_{}
&V_2 \ar[r]_{} \ar[dr]_{}
&\rho_1 \\
&V_1 \ar[r]_{} 
&V_1 \ar[r]_{} \ar[ur]_{} 
&V_2 \ar[r]_{} \ar[ur]_{} 
&\rho_2.
}\]
As another example, the orbit cone $C_{\mbf{v}_{1,6}}$ for $\mbf{v}_{1,6}=(\mbf{v}_1, 0,0,0,1)$ is defined by
$$\theta_3\ge0,\;\theta_1+2\theta_3\ge0,\;\theta_2+2\theta_3\ge0,\;\theta_1+\theta_2+2\theta_3\ge0,$$
$$\theta_3+\theta_4\le0,\;\theta_1+\theta_3+3\theta_4\le0,\; \theta_2+\theta_3 +3\theta_4\le0,\;\theta_1+\theta_2+\theta_3+3\theta_4\le0$$
and the structure of the $G$-constellation $F$ for $\mbf{v}_{1,6}$ is presented as 

\[\xymatrix@=36pt{
V_2 \ar[r]_{} \ar[d]_{} \ar[r]_{}
&V_2 \ar[d]_{}
&\rho_0 \ar[r]_{} \ar[dr]_{} 
&V_1 \\
V_1 \ar[r]_{}
&V_2 \ar[r]_{} \ar[ur]_{} \ar[dr]_{}
&\rho_1 \ar[r]_{} \ar[ur]_{}
&V_1. \\
&
&\rho_2 \ar[ur]_{} \ar[uur]_{}
&
}\]
As a more explicit description, $F$ can be realized as a quotient of the ideal of $\C[V]$ generated by $xy,yz$ and $zx$. Note that $\C xy\oplus\C yz\oplus\C zx$ is isomorphic to $V_2$ as a $\C G$-module.

Finally, we compute the homomorphism $\psi_C$ for certain chambers $C$ and see that every projective crepant resolution of $V/G$ is obtained as a moduli space of $G$-constellations. Let $C_{i,j}\subset \Theta$ be the orbit cones corresponding to the rays of $\sigma_\mcl{W}$ generated by $\mbf{v}_{i,j}$. There are 42 choices of a triple $(i,j,k)$ such that $\overline{C}_{1,i}\cap \overline{C}_{2,j}\cap \overline{C}_{3,k}$ has the full dimension 4. Among them, $C_+:=\overline{C}_{1,1}\cap \overline{C}_{2,1}\cap \overline{C}_{3,1}$ is the one containing the GIT-chamber $C_{++}$ for the $G$-Hilbert scheme. $C_{++}$ is in fact defined by the inequalities $\theta_i>0\,(i=1,2,3,4)$.

Similarly to the abelian case in Example \ref{abel2}, one can compute $\psi_{C_{++}}$ and obtain $\psi_{C_{++}}(C_+)=\mrm{Mov}(X)$. Thus, every projective crepant resolution of $V/G$ is realized as $\mfk{M}_\theta$ for some generic $\theta\in C_+$. In this case the movable cone $\mrm{Mov}(X)$ is a three-dimensional simplicial cone divided into four nef cones (cf. \cite[Proposition 5.2]{DG}). The nef cone for the $G$-Hilbert scheme is given as $\psi_{C_{++}}(\overline{C}_{++})$ and we see that it is the central one which can pass to the three remaining (corner) ones. Note that not all of $\overline{C}_{1,i}\cap \overline{C}_{2,j}\cap \overline{C}_{3,k}$ of dimension 4 realize all the resolutions. For example, $\psi_{C'}(C)$ is strictly smaller than $\mrm{Mov}(X)$ for the cone $C=\overline{C}_{1,2}\cap \overline{C}_{2,1}\cap \overline{C}_{3,1}$ adjacent to $C_+$ and any GIT-chamber $C'\subset C$.

\section{Proof of the main result}\label{4}

In this section we give a proof of the following theorem:

\begin{thm}\label{main2}(=Theorem \ref{main1})
Let $G\subset SL_3(\C)$ be a finite subgroup and $X\to\C^3/G$ a projective crepant resolution. Then there is a generic stability condition $\theta\in\Theta$ such that $X\cong \mfk{M}_\theta$.
\end{thm}

To prove Theorem \ref{main2}, we will use a similar argument to the one in \cite[\S8]{CI}. As explained in Introduction, the strategy is to show that, starting from any GIT-chamber $C\subset \Theta$ and the associated moduli space $\mfk{M}_C$, one can reach a wall in $\Theta$ which induces a given flop $\mfk{M}_C \dashrightarrow X'$ by crossing walls of a certain type. Here we say that the intersection $W=\overline{C}\cap\overline{C'}$ for GIT-chambers $C, C'\subset \Theta$ is a {\it wall} (of $C$) if $W$ is a codimension-one face of $C$ (and hence of $C'$). To carry out the strategy, we should know how the tautological bundle $\mcl{R}_C$ of $\mfk{M}_C$ changes when we cross a wall in $\Theta$. 

For any given GIT-chamber $C\subset \Theta$, we choose any general point $\theta_0$ of a wall $W\subset C$. Then we have $\mcl{W}^{C\text{-ss}}\subset\mcl{W}^{\theta_0\text{-ss}}$ and this inclusion induces a morphism $\alpha:\mfk{M}_C\to\bar{\mfk{M}}_{\theta_0}$ over $\mfk{M}_0=\C^3/G$ where $\bar{}$ denotes the normalization. Similarly to the case of $\mfk{M}_C$, we denote the moduli space $\bar{\mfk{M}}_{\theta_0}$ by $\bar{\mfk{M}}_W$. We also denote the crepant resolution $\mfk{M}_C$ by $X$ as in previous sections.

Recall that each line bundle $L\in\mrm{Pic}(X)$ defines a character $\chi_L\in \chi(T_X)$ and its associated GIT-quotient $X_L:=\mfk{X}/\!/_{\chi_L^k}T_X$ with a sufficiently divisible $k>0$. As observed in Subsection \ref{3.2}, $X_L$ is the same as the image of the rational map obtained by the complete linear system of (some power of) $L$. Note that $X_0=\C^3/G$ and that $X_L$ makes sense for $L\in\mrm{Pic}(X)_\R$. The following lemma shows that one can determine $\alpha:\mfk{M}_{C}\to\bar{\mfk{M}}_W$ by looking at the homomorphism $\psi_C:\Theta\to \mrm{Pic}(X)_\R$ in (\ref{psi}):

\begin{lem}\label{descent}(cf. \cite[Lemma 3.3]{CI})
For any point $\theta\in \overline{C}$, the normalized GIT-quotient $\bar{\mfk{M}}_\theta$ is isomorphic to $X_{\psi_C(\theta)}$.
\end{lem}

\begin{proof}
Recall that $X=\mfk{M}_C$ is regarded as a GIT-quotient of $\mfk{X}=\mrm{Spec}\,\Cox(X)$ via the restriction $\mfk{X}^{\psi_C(C)\text{-ss}}\to \mcl{W}^{C\text{-ss}}$ of the equivariant (affine) morphism $\vphi_C^*: \mfk{X}\to \mcl{W}$ obtained from $\vphi_C$ (see Subsection \ref{3.3}). Since $\theta\in \overline{C}$, the line bundle $L:=\psi_C(\theta)$ induces a restriction $\mfk{X}^{L\text{-ss}}\to \mcl{W}^{\theta\text{-ss}}$ of $\vphi_C^*$ which extends $\mfk{X}^{\psi_C(C)\text{-ss}}\to \mcl{W}^{C\text{-ss}}$, by the construction of $\psi_C$. Then the composition $\mfk{X}^{\psi_C(C)\text{-ss}}\hookrightarrow \mfk{X}^{L\text{-ss}}\to \mcl{W}^{\theta\text{-ss}}$ induces a morphism $X\to \bar{\mfk{M}}_\theta$ and $L$ is just the pullback of the ample line bundle $\mcl{O}_{\bar{\mfk{M}}_\theta}(1)$ (induced by the GIT-construction) under this morphism. Thus, the claim follows since $X\to\bar{\mfk{M}}_\theta$ is identified with the morphism coming from the linear system of $L$. Note in particular that (a multiple of) $L$ is globally generated.
\end{proof}

The following technical lemma will be an important step for the proof of the main theorem.

\begin{lem}\label{lem:vertical}
Let $F$ be a face of an orbit cone $\tilde{C}_x\subset\Theta$ of the $G$-constellation corresponding to a point $x\in\mfk{M}_C$. If the relative interior $\mrm{relint}(F)$ of $F$ intersects $\overline{C}$ and if $\psi_C(F)\subset\mrm{Pic}(X)_\R$ is not of full dimension, then $\psi_C(F)$ is on the boundary of the orbit cone $C_x\subset \mrm{Pic}(X)_\R$ of a point of $\mfk{X}$ which is a lift of $x\in X$.

In particular, in this case, any point $\theta_0$ of $\mrm{relint}(F)\cap\overline{C}$ induces a nontrivial birational contraction $X\to \bar{\mfk{M}}_{\theta_0}$ such that every point $x'\in X$ whose orbit cone in $\Theta$ is equal to $\tilde{C}_x$ lies in the exceptional locus.
\end{lem}

\begin{proof}
We first choose homogeneous generators $h_1,\dots,h_\ell$ of $\C[\mcl{W}]=\C[\mcl{V}]^{SL_R}$ so that we have an associated embedding $\mcl{W}\hookrightarrow Y_\mcl{W}$ and its quotient $X\hookrightarrow Y$ as in Theorem \ref{embedding}. Let $\bar{x}\in\mcl{W}$ be any lift of $x$. By reordering $h_i$'s, we may assume that\\
(1) there exists $\ell_1<\ell$ such that $h_1,\dots,h_{\ell_1}$ do not vanish at $\bar{x}$ and their weights generate the cone $F$ and that\\
(2) for any $h_k$ with $k>\ell_1$, either $h_k$ vanishes at $\bar{x}$ or $\mrm{wt}(h_k)$ is not in $F$.\\
The locally closed subset of $\mcl{W}$ consisting of points having $F$ as their orbit cones is given by
$$\mcl{W}_F:=\mcl{W}\cap \{h_1 h_2\cdots h_{\ell_1}\ne0\}\cap \{h_{\ell_1+1}=\cdots=h_{\ell}=0\}.$$
Note that torus-orbits of the toric varieties $Y,\,\mathbb{A}^{\ell+m}$, and $Y_\mcl{W}$ restrict to stratification of $X, \mfk{X}$, and $\mcl{W}$ respectively, and that $\mcl{W}_F$ is such a stratum of $\mcl{W}$.

By construction, the equivariant map $\vphi_C^*:\mfk{X}\to \mcl{W}$ associated to $C$ is the restriction of the toric morphism $\mathbb{A}^{\ell+m}\to Y_\mcl{W}$ which is obtained by sending the rays $\tilde{\tau}_k$ of the fan of $\mathbb{A}^{\ell+m}$ corresponding to the exceptional divisors $E_k$ to the rays $\tau_{C,k}$ of $\sigma_\mcl{W}$ corresponding to $E_k$ in $\mfk{M}_C$. For $\theta_0\in \mrm{relint}(F) \cap \overline{C}$, the restriction $\mfk{X}^{\psi_C(\theta_0)\text{-ss}}\to \mcl{W}^{\theta_0\text{-ss}}$ of $\vphi_C^*$ descends to an isomorphism $X_{\psi_C(\theta_0)}\to\bar{\mfk{M}}_{\theta_0}$ by Lemma \ref{descent}. Note that $\mcl{W}_F$ is a closed subset of $\mcl{W}^{\theta_0\text{-ss}}$ which intersects the closure of the $T_\Theta$-orbit of $\bar{x}$. In particular, the categorical quotient of $\mcl{W}_F$ by $T_\Theta$ is again a closed subset in $\mfk{M}_{\theta_0}$.

We ``perturb" the map $\vphi_C^*$ to obtain another equivariant map $\vphi'^*:\mfk{X}\to \mcl{W}$ such that its restriction to the semistable loci still descends to an isomorphism $X_{\psi_C(\theta_0)}\to \bar{\mfk{M}}_{\theta_0}$ and that every point of $\mcl{W}_F\subset \mcl{W}^{\theta_0\text{-ss}}$ is inside a $T_\Theta$-orbit of a point in the image of $\vphi'^*$. Such $\vphi'^*$ is explicitly constructed as follows. For each $k=1,\dots,m$, let $\sigma_k$ be the cone (contained in $\sigma_\mcl{W}$) generated by $\tau_{C,k}$ and the rays $\tau$ of $\sigma_\mcl{W}$ (if any) such that $\tau$ is projected to $\tau_k$ under $p:(N_\mcl{W})_\R\to (N_\mcl{W}/\bar{\Theta}^\vee)_\R$ (see Subsection \ref{3.4}) and that the orbit cone of $\tau$ contains $F$ as a face. Note that $\sigma_k$ is equal to $\tau_{C,k}$ if $\tau_{C,k}\not\subset \sigma_F$ and otherwise $\sigma_k$ is the intersection $\sigma_F\cap p^{-1}(\tau_k)$ where $\sigma_F$ is the face of $\sigma_\mcl{W}$ corresponding to the stratum $\mcl{W}_F$. Then $\vphi'^*$ is obtained as the restriction of the toric morphism $\mathbb{A}^{\ell+m}\to Y_\mcl{W}$ obtained by sending each ray $\tilde{\tau}_k$ to a ray generated by a general point of $\sigma_k$. Note that $\vphi'^*$ indeed induces an isomorphism $X_{\psi_C(\theta_0)}\to \bar{\mfk{M}}_{\theta_0}$ since the projected cone $p(\sigma)$ of the image $\sigma$ of each cone of the fan of $\mfk{X}^{\psi_C(\theta_0)\text{-ss}}$ under the induced homomorphism $\R^{\ell+m}\to \R^{\ell+s-1}$ by $\vphi_C^*$ remains the same even if we replace $\vphi_C^*$ by $\vphi'^*$.

The map $\vphi'^*$ gives an associated homomorphism $\psi': \Theta\to \mrm{Pic}(X)_\R$ similarly to $\psi_C$. More concretely, for a homogeneous semi-invariant $f\in\C[\mcl{W}]$ of weight $c\in \Theta$, $\psi'(c)$ is the class of the sum of the strict transform of the image in $V/G$ of $\{f=0\}$ and $\sum_k v_k(f)E_k$ where $v_k$ is the primitive element in $\sigma_k$ defining $\vphi'^*$ (cf. Remarks \ref{rem:cox(V/G)} and \ref{rem:vphi}(1)).

One sees that $\psi'|_F=a\psi_C|_F$ for some constant $a\in \Q$ and in particular $\psi'(F)$ is not full-dimensional by the assumption. Let $\vphi': \C[\mcl{W}]\to \C[\mfk{X}]$ be the ring homomorphism associated to $\vphi'^*$. The equivariant map $\vphi'^*$ restricts to a map from the nonempty (affine) locally closed subset
\[
\begin{aligned}
\mfk{X}_F:&=(\vphi'^*)^{-1}(\mcl{W}_F)\\
&=\mfk{X}\cap \{\vphi'(h_1 h_2\cdots h_{\ell_1})\ne0\}\cap \{\vphi'(h_{\ell_1+1}) =\cdots =\vphi_C(h_{\ell})=0\}
\end{aligned}
\]
of $\mfk{X}$ to $\mcl{W}_F$, which descends to an isomorphism between the (categorical) quotients by $T_X$ and $T_\Theta$.

We assume that $\psi'(F)$ does not lie on the boundary of $C_x$ to deduce a contradiction. Then the quotients of $\mfk{X}_F$ and $\mcl{W}_F$ are isomorphic to the stratum $S_x$ of $X$ containing $x$. If $\mfk{X}_F$ itself is a stratum of $\mfk{X}$, then it must be the set $\mfk{X}_x$ of points of $\mfk{X}$ whose orbit cones are equal to $C_x$. The fact that $\psi'(F)$ is not of full-dimension implies that every point of $\mfk{X}_x$ has a positive-dimensional stabilizer in $T_X$, which is a contradiction since $C_x$ is full-dimensional. If there exists a stratum $\mfk{X}_{x'}\subsetneq\mfk{X}_F$ not equal to $\mfk{X}_x$, then this is again a contradiction since the orbit cone of such $\mfk{X}_{x'}$ is a proper face of $C_x$ and in particular $\mfk{X}_{x'}$ cannot be mapped to $S_x$ under the quotient map. Note that proper faces of an orbit cone of $\mfk{X}$ always correspond to nontrivial contractions, in contrast to faces of orbit cones in $\Theta$ such as walls of {\it Type 0} introduced below. The last claim also follows from this observation.
\end{proof}

If a wall $W=\overline{C}\cap\overline{C'}$ is not contained in the pullback under $\psi_C: \Theta\to \mrm{Pic}(X)_\R$ of a GIT wall for the action of $T_X$ on $\mfk{X}$, then the morphism $\mfk{M}_C\to\bar{\mfk{M}}_{W}$ is an isomorphism by Lemma \ref{descent}, and we call such a wall of {\it Type 0} following \cite{CI}. Conversely, for any wall $W$ whose image under $\psi_C$ is not of full dimension, the morphism $\mfk{M}_C\to\bar{\mfk{M}}_{W}$ is a non-trivial contraction of either a divisor or a curve by Lemma \ref{lem:vertical}. For a wall $W$ of Type 0, the two moduli spaces $\mfk{M}_C$ and $\mfk{M}_{C'}$ are isomorphic but the tautological bundle changes. In fact $\mcl{R}_{C'}$ is a modification of $\mcl{R}_C$ by tensoring by a line bundle to a subbundle of $\mcl{R}_C$ similar to the abelian cases \cite[Corollary 4.2]{CI}, as we will see below.

In order to study what happens when we cross a wall $W$ of Type 0, we first see that the {\it unstable locus}
$$D_W=\{x\in X\mid \text{any lift }\bar{x}\in\mcl{W}\text{ of }x\text{ is not }W\text{-stable}\}\subset X$$
with respect to $W$ is a compact divisor as a topological space. Later we will see that $D_W$ has a natural scheme structure as a certain moduli space and show that $D_W$ is reduced with respect to this scheme structure in Lemma \ref{lem:tautological}. We will also show that $D_W$ is connected in Proposition \ref{prop:twist}.

\begin{lem}\label{divisor}
For any given GIT-chamber $C\subset\Theta$ and its wall $W$ of Type 0, the  unstable locus $D_W\subset X$ associated to $W$ is a compact divisor.
\end{lem}

\begin{proof}
We first show that $D_W$ is a divisor. We fix an embedding $\mcl{W}\hookrightarrow Y_\mcl{W}$ again, and let $\Sigma_\theta$ be the subfan of $\sigma_\mcl{W}$ defining $Y_\mcl{W}^{\theta\text{-ss}}$ for $\theta\in \Theta$ (see Subsection \ref{3.4}). We also consider the ambient toric variety $Y$ of $X$ and let $\sigma_x$ be the cone of the fan $\Sigma_X$ of $Y$ corresponding to a point $x\in X\subset Y$.

For any point $x\in D_W$, there exists a unique cone $\tilde{\sigma}_x\in \Sigma_\theta$ for $\theta\in C$ such that $p(\tilde{\sigma}_x)=\sigma_x$. Since $\mfk{M}_C\to\bar{\mfk{M}}_W$ is an isomorphism, we have a cone $\bar{\sigma}_x\in\Sigma_{\theta_0}$ strictly containing $\tilde{\sigma}_x$ (as a facet) such that it satisfies $p(\bar{\sigma}_x)=\sigma_x$ where $\theta_0$ is a general point of $W$. Recall that the maximal cones of the fan of $\mcl{W}^{W\text{-ss}}$ are projected under $p$ to the maximal cones of $\Sigma_X$ (see (\ref{condition*})). Since the image of any ray of $\sigma_\mcl{W}$ under $p$ is again a ray of $\Sigma_X$, there must be a ray $\tau$ of $\tilde{\sigma}_x$ and a face $\bar{\tau}\subset \bar{\sigma}_x$ such that $\bar{\tau}\supsetneq\tau$ and $p(\tau)=p(\bar{\tau})$. This implies that $x$ is contained in the unstable divisor of $X$ corresponding to $p(\bar{\tau})$. Thus, $D_W$ is of pure codimension one.

For the compactness of $D_W$, one may apply exactly the same argument in \cite[Proposition 4.4]{CI}.
\end{proof}

Now we describe how $G$-constellations over unstable points change when we cross the wall $W$. We first observe that there is a unique maximal $\C G$-submodule $R_+$ of the regular representation $R$ of $G$ which satisfies $\theta_0(R^+)=0$ and $\theta(R^+)>0$ for a general $\theta_0\in W$ and $\theta\in C$. For a $G$-constellation $F$ corresponding to any point $x\in D_W \subset X=\mfk{M}_C$, we define $S\subsetneq F$ as the unique nonzero minimal $\C[V]$-submodule such that $\theta_0(S)=0$ if $R_+$ is primitive in $R(G)=\bigoplus_{\rho\in \mrm{Irr}(G)} \Z\rho$ and otherwise we define $S\subsetneq F$ as the unique (nonzero) maximal $\C[V]$-submodule such that $\theta_0(S)=0$. We also let $Q$ be the quotient of $F$ by $S$. Note that $R$ is primitive as an element of $R(G)$ and thus at least $R_+$ or its complement in $R$ is primitive in $R(G)$.

In abelian cases every irreducible representation of $G$ has multiplicity one in $R$ and thus the modules $S$ and $Q$ are clearly indecomposable by the $\theta_0$-semistability. When $G$ is non-abelian and either $S$ or $Q$, say $S$, is non-primitive in $R(G)$, the module $S$ may be decomposable in general. Even if $S$ is indecomposable, the $\theta_0$-polystable module $\bar{F}$ of $Q_G$ (i.e. the direct sum of $\theta_0$-stable modules) corresponding to the image of $x$ under $\mfk{M}_C\to\mfk{M}_W$ may not be isomorphic to the direct sum $S\oplus Q$. $S$ may degenerate in $\bar{F}$ to decompose nontrivially into a direct sum $\bigoplus_i \bar{S}_i$ satisfying $\theta_0 (\bar{S}_i)=0$.

We will identify $D_W$ with a certain moduli space $Z_Q$ which plays the same role as the moduli space $Z$ introduced in the proof of \cite[Lemma 3.10]{CI}. $Z_Q$ is constructed as follows. We first choose an irreducible component $D$ of $D_W$, and take its general point and the corresponding $\C[V]$-modules $S, Q$ as above. We then fix a direct summand $R_2$ of $R$ such that $R_2\cong Q$ as a $\C G$-module. Similarly to the construction of $\mcl{N}$ in Subsection \ref{2.1}, we define $\mcl{N}_Q$ to be the space of $G$-equivariant $\C[V]$-module structures on $R$ such that $R_2$ is a quotient $\C[V]$-module. More precisely, if we regard $\mcl{N}$ as the space of quiver representations (see Subsection \ref{2.3}) and fix isotypic decompositions of $R_2$ and its complement $R_1\cong S$ in $R$, then $\mcl{N}_Q$ is defined as the closed subscheme of $\mcl{N}$ given by vanishing of the coordinates (i.e. the entries of matrices) corresponding to arrows which send vectors in $R_1$ to those in $R_2$. Then $Z_Q$ is constructed as the (geometric) quotient of the set $\mcl{N}_Q^{C\text{-ss}}=\mcl{N}_Q \cap \mcl{N}^{C\text{-ss}}$ of $C$-(semi)stable points by the action of $GL_{R_1} \times GL_{R_2}$, similar to the construction of $\mfk{M}_C$. Note that $Z_Q$ admits a natural map $Z_Q\to \mfk{M}_C$, which will turn out to be a closed immersion.

\begin{lem}\label{lem:reduced}
The scheme $Z_Q$ is reduced.
\end{lem}

\begin{proof}
It suffices to show that $\mcl{N}_Q^{C\text{-ss}}$ is reduced. Let
\[
I_\mcl{N}\subset A_G:=\C[\mrm{Hom}_{\C[G]}(V^*\otimes_\C R, R)]
\]
be the defining ideal of $\mcl{N}$ (see (\ref{eq:B wedge B})). Note that $A_G$ is identified with a polynomial ring whose variables are entries of matrices corresponding to arrows of $Q_G$. As mentioned above, the variables of $A_G$ corresponding to the arrows from $R_1$ to $R_2$ generate an ideal $I_S$ such that the defining ideal of $\mcl{N}_Q$ is equal to $I_\mcl{N}+I_S$. Then, by noticing that the defining equations of $\mcl{N}$ come from the commutativity of the actions of the variables $x,y,z\in V^*$, one can check that $I_\mcl{N}$ is generated by polynomials $\{p_i\}_i$ such that every term of $p_i$ is in $I_S$ as long as $p_i$ is in $I_S$. (For example, if the action of $x\otimes y$ sends a vector in $R_1$ to one in $R_2$, so does the action of $y\otimes x$.) This implies that a power $f^m\,(m>0)$ of a polynomial $f\in A_G$ lies in $I_S+I_\mcl{N}$ only if the power $f_0^m$ of the sum $f_0$ of terms of $f$ not in $I_S$ lies in $I_\mcl{N}$. Since $\mcl{N}^{C\text{-ss}}$ is smooth and in particular reduced, it follows that $\mcl{N}_Q^{C\text{-ss}}$ is reduced as well.
\end{proof}

By the choice of $S$, there is no nonzero $G$-equivariant $\C[V]$-linear map $S\to Q$ (i.e. $G\text{-}\Hom_{\C[V]}(S,Q)=0$) and thus there exists a Zariski open subset of $Z_Q$ which admits an immersion into $D$ (cf. the proof of \cite[Lemma 3.10]{CI}). By performing the same operation on each irreducible component of $D_W$, we get a (topologically) dense Zariski open subset $U$ of $D_W$ together with a universal quotient map $\mcl{R}_C|_U\to \bar{\mcl{Q}}$ whose restriction to $U\cap D$ equals the universal quotient map associated to $U\subset Z_Q$. We choose $U$ as the maximal one satisfying this condition, and  denote by $\mfk{Z}$ the complement of $U$ in $D_W$ as a topological space. Then $\mfk{Z}$ parametrizes $G$-constellations whose associated quotient $\C[V]$-modules are strictly larger than $R_2$ as $\C G$-modules. Note that this may a priori happen for non-abelian $G$. Also, the rank of $\bar{\mcl{Q}}$ a priori depends on irreducible components of $D_W$ but it will turn out that $\mfk{Z}$ is empty and that $R_2$ is uniquely determined from $W$ (see Remark \ref{rem:singleQ}).

We take a 1-PS $\lambda: \C^*\to GL_R$ (for each irreducible component of $D_W$) defined by assigning weight one to the generators of the ideal $I_S\subset A_G$ in the proof of Lemma \ref{lem:reduced} and weight zero to the other variables of $A_G$. Then $\lambda$ is positive on $C$ and zero on $W$ as a function on $\Theta\cong\chi(PGL_R)_\R$. $\lambda$ induces a filtration of $F$ by its submodules (cf. \cite[\S3]{Kin}) such that it has two associated graded components one of which is isomorphic to $Q$. This works for a family of $G$-constellations over the open subset $U$ of $D_W$ and thus we obtain an exact sequence 
\begin{equation}\label{eq:exact1}
0\to\bar{\mcl{S}}\to \mcl{R}_C|_U\to \bar{\mcl{Q}}\to0
\end{equation}
of flat families of $\C[V]$-modules. By construction, the quotient map $\mcl{R}_C|_U\to \bar{\mcl{Q}}$ in (\ref{eq:exact1}) coincides with the one induced by the schemes $Z_Q$. Note also that the sequence (\ref{eq:exact1}) is a nontrivial extension as $\C[V]$-modules even if it splits as vector bundles on $U$. It will be shown that this sequence extends to whole $D_W$ in Lemma \ref{lem:tautological}.

Let $\bar{\mcl{K}}$ be the kernel of the natural surjection $\mcl{R}_C|_{X\setminus \mfk{Z}}\to \bar{\mcl{Q}}$. Then $\bar{\mcl{K}}$ is a vector bundle on $X\setminus \mfk{Z}$ (with a $\C[V]$-action) called the {\it elementary transformation} of $\mcl{R}_C|_{X\setminus \mfk{Z}}$ (associated to $\mcl{R}_C|_{X\setminus \mfk{Z}}\to \bar{\mcl{Q}}$) in the sense of \cite{M}. Applying the inverse of the elementary transformation to $\bar{\mcl{K}}$, we obtain an exact sequence
\begin{equation}\label{eq:exact2}
0\to\bar{\mcl{Q}}\otimes \mcl{O}_{X\setminus \mfk{Z}}(-U) \to \bar{\mcl{K}}|_U\to \bar{\mcl{S}}\to0
\end{equation}
of flat families of $\C[V]$-modules. Note that $\bar{\mcl{K}}$ parametrizes the same $G$-constellations as $\mcl{R}_C|_{X\setminus \mfk{Z}}$ outside $U$ while sub and quotient modules are switched over $U$. 

\begin{lem}\label{lem:tautological}
The family $\bar{\mcl{K}}$ over $X\setminus \mfk{Z}$ extends to a family $\mcl{K}$ over the whole $X$ which coincides with the tautological bundle $\mcl{R}_{C'}$ for the GIT-chamber $C'$ up to tensoring by a line bundle on $X$.
\end{lem}

\begin{proof}
First recall that the bundle $\mcl{R}_{C'}=\bigoplus_{\rho\in\mrm{Irr}(G)} \mcl{R}'_\rho$ is normalized by tensoring by a line bundle so that the component $\mcl{R}'_{\rho_0}$ becomes the trivial line bundle. Since $\mfk{Z}\subset X$ is of codimension at least two and $X$ is smooth, it suffices to show that $\bar{\mcl{K}}$ parametrizes $\theta'$-(semi)stable $G$-constellations for $\theta'\in C'$. To this end, we only have to show that the fibers of $\bar{\mcl{K}}$ are indecomposable as $\C[V]$-modules.

Let $F$ be the $G$-constellations corresponding to a (general) point of $x\in U\subset\mfk{M}_C$. If the fibers $S\subset F$ and $Q=F/S$ of $\bar{\mcl{S}}$ and $\bar{\mcl{Q}}$ respectively are primitive as elements in $R$, then they are indecomposable and the same argument as in the proof of \cite[Proposition 4.1]{CI} shows that the sequence (\ref{eq:exact2}) is a nontrivial extension of $\C[V]$-modules, which implies that the fibers of $\bar{\mcl{K}}$ are indecomposable as well.

If one of $S$ and $Q$ is not primitive in $R$, we may assume that $S$ is not primitive and $Q$ is primitive by switching the roles of $C$ and $C'$, if necessary. Then we have a filtration
$$0= \bar{\mcl{S}}_0\subset \bar{\mcl{S}}_1\subset \cdots \subset \bar{\mcl{S}}_{k-1}= \bar{\mcl{S}}\subset \bar{\mcl{S}}_k=\mcl{R}_C|_U$$
consisting of $\C[V]$-equivariant subbundles of $\mcl{R}_C|_U$ such that each fiber of the associated graded bundle $\bigoplus_{i=1}^k \bar{\mcl{S}}_i /\bar{\mcl{S}}_{i-1}$ represents the $W$-polystable module of the corresponding point in $\mfk{M}_W$. This filtration can be obtained by an appropriate 1-PS, similarly to (\ref{eq:exact1}). Since fibers of each component $\bar{\mcl{S}}_i /\bar{\mcl{S}}_{i-1}$ are indecomposable, we can reduce the claim to the previous case by applying elementary transformations repeatedly to quotients by $\bar{\mcl{S}}_i /\bar{\mcl{S}}_{i-1}$ for $i=1,\dots, k-1$.
\end{proof}

\begin{lem}\label{lem:unstable}
$D_W$ is isomorphic to a disjoint union of (the underlying topological spaces of) $Z_Q$'s for all possible $Q$.
\end{lem}

\begin{proof}
By the previous lemma, the extended bundle $\mcl{K}$ is isomorphic to a vector bundle of the form $\mcl{R}_{C'}\otimes L$ for some $L\in \mrm{Pic}(X)$. Let $\overline{U}$ be the scheme-theoretic closure so that it a reduced Cartier divisor whose support is equal to $D_W$. Then there is a quotient morphism over $\overline{U}$ from $\mcl{R}_C|_{\overline{U}}$ to the cokernel $\mcl{Q}$ of $\mcl{K}\to \mcl{R}_C$, whose restriction to $U$ coincides with the universal quotient morphism on $U$. We show that each connected component of $\overline{U}$ is identified with some $Z_Q$ and that $\mcl{R}_C|_{\overline{U}}\to \mcl{Q}$ is identified with the universal quotient morphism for $Z_Q$. For this, it suffices to show that every fiber of $\mcl{Q}$ on each connected component of $D_W$ has the same dimension. Indeed, this condition implies that every point of each connected component defines the same $Q$ as a $\C G$-module and hence the whole $Z_Q$ admits a closed immersion onto $\overline{U}$.

Let $\mcl{S}$ be the kernel of $\mcl{R}_C|_{\overline{U}}\to \mcl{Q}$. Then every fiber of $\mcl{S}$ has the dimension equal to or less than $\dim_\C S=\mrm{rank}\, \mcl{R}_C-\dim_\C Q$. However, by considering the inverse of elementary transformation for the same wall $W$ with the roles of $C$ and $C'$ switched, we also see that every fiber of $\mcl{S}$ has the dimension at least $\dim_\C S$. Therefore, every fiber of $\mcl{S}$ has the same dimension, and the same holds for $\mcl{Q}$ as well.
\end{proof}

\begin{rem}\label{rem:singleQ}
The previous lemma implies that the whole $Z_Q$ is identified with the scheme $Z$ introduced in \cite[Lemma 3.10]{CI}. As we will see in Proposition \ref{prop:twist} that $D_W$ is connected and hence $D_W$ is in fact identified with a single reduced scheme $Z_Q$.
\end{rem}

Next we describe the behaviour of the Grothendieck group $K(X)$ of locally free sheaves on $X$ under the wall-crossing. As we will see below, it is basically the same as the abelian cases and this enables us to proceed the proof of the main result in a similar way to \cite{CI}. To this end we briefly review the necessary notions. See \cite[\S2.4, \S5, and \S7.1]{CI} and also \cite[\S9.2]{BKR} for details.

The result of \cite{BKR} shows that the Fourier-Mukai transform by the universal family $\mcl{U}_C$ gives an equivalence $\Phi_C: D(X)\cong D^G(\C^3)$ whose inverse $\Phi^{-1}_C: D^G(\C^3)\to D(X)$ is given by $\mcl{O}_{\C^3}\otimes\rho\mapsto \mcl{R}^\vee_\rho$. $\Phi_C$ induces an isomorphism between $K(X)$ and the Grothendieck group $K^G(\C^3)$ of $D^G(\C^3)$. Let $K_0(X)$ be the Grothendieck group of the subcategory of $D(X)$ of complexes whose supports are contained in the fiber $\pi^{-1}(0)$ for $\pi:X \to \C^3/G$. Then $K_0(X)$ is identified with the dual of $K(X)$ via the perfect pairing $\chi(-, -): K(X)\times K_0(X)\to\Z$ defined by 
$$(\alpha, \beta)\mapsto \sum_i (-1)^i \dim \mrm{Ext}_{\mcl{O}_X}^i(\alpha, \beta),$$
which is skew-symmetric on $K_0(X)$. Via the equivalence $\Phi_C$ we also have a perfect pairing $\chi^G(-,-)$ between $K^G(\C^3)$ and the Grothendieck group $K^G_0(\C^3)$ of $G$-equivariant sheaves on $\C^3$ supported at the origin. With respect to this pairing, the classes of $\{\mcl{O}_{\C^3}\otimes \rho\}_{\rho\in\mrm{Irr}(G)}$ and $\{\mcl{O}_0 \otimes \rho\}_{\rho\in\mrm{Irr}(G)}$ give mutually dual basis of $K(\C^3)$ and $K_0(\C^3)$ respectively where $\mcl{O}_0$ is the skyscraper sheaf at the origin of $\C^3$.

The equivalence $\Phi_C$ also induces an isomorphism $\phi_C: K_0(X)\to K^G_0(\C^3)$ and its pullback $\phi_C^*: K^G(\C^3)\to K(X)$. The vector space $\Theta$ is regarded as a subspace of $K^G(\C^3)_\R$, and the restriction of $\phi_C^*$ to $\Theta$ gives an isomorphism $\Theta\to F^1$ where we have a filtration
$$K(X)_\R=F^0\supset F^1\supset F^2\supset F^3=0$$
with $F^i$ being the subspace generated by sheaves whose supports are of codimension at least $i$. The homomorphism $\phi_C^*|_\Theta:\Theta \to F^1$ is explicitly given as 
$$\theta=(\theta_\rho)_\rho\mapsto \sum_{\rho\in\mrm{Irr}(G)} \theta_\rho[\mcl{R}_\rho^\vee]$$
where $[\mcl{R}_\rho^\vee]$ denotes the class of the dual sheaf $\mcl{R}_\rho^\vee$ in $F^1$.
We then obtain a homomorphism $L_C:\Theta \to \mrm{Pic}(X)_\R$ such that
the following diagram of vector spaces is commutative:
\begin{equation}\label{Pic}
\begin{CD}
\Theta @> \phi_C^*|_\Theta >> F^1\\
@V L_C VV @V \mrm{pr} VV\\
\mrm{Pic}(X)_\R @<< \det^{-1} < F^1/F^2\\ 
\end{CD}
\end{equation}
where $\mrm{pr}$ is the projection and the lower horizontal arrow denotes the isomorphism $F^1/F^2\cong \mrm{Pic}(X)_\R$ given by sending the class $[\mcl{F}]\in F^1$ to $\det(\mcl{F})^{-1}$.

\begin{lem}\label{lem:commute}
$L_C$ is equal to the homomorphism $\psi_C:\Theta\to\mrm{Pic}(X)_\R$ in (\ref{psi}).
\end{lem}

\begin{proof}
By definition $L_C$ is given as 
$$\theta=(\theta_\rho)_\rho\mapsto \det\left(\bigoplus_{\rho\in\mrm{Irr}(G)} \mcl{R}_\rho^{\oplus \theta_\rho}\right).$$

For a determinantal semi-invariant $f$ having $\theta=(\theta_\rho)_\rho$ as its weight, $\vphi_C(f)$ is regarded as a section of a line bundle $L_C(\theta)$ by noticing that each arrow $a_{\rho,\rho'}$ in $Q_G$ is regarded as a section of $\mcl{R}_\rho^\vee\otimes \mcl{R}_{\rho'}$. Thus, $L_C(\theta)$ is nothing but the line bundle $\psi_C(\theta)$. Since $\Theta$ is spanned by the weights of determinantal semi-invariants by Theorem \ref{thm:semi-invariants}, $L_C$ and $\psi_C$ are the same on the whole $\Theta$. 
\end{proof}

When we cross a wall from a GIT-chamber $C$ to an adjacent one $C'$, the corresponding chambers $\phi_C^*(C)$ and $\phi_{C'}^*(C')$ in $F^1$ may not be adjacent in general. For example, in dimension two, the automorphism $\phi_{C'}^*\circ (\phi_{C}^*)^{-1}$ of $F^1\cong \mrm{Pic}(X)$ is a reflection along the hyperplane containing $\phi_{C}^*(W)$ for the wall $W=\overline{C}\cap\overline{C'}$. In the proof of the main result of \cite{CI}, namely the abelian case of Theorem \ref{main2}, the fact that, for a wall-crossing of Type 0, the chambers $\phi_C^*(C)$ and $\phi_{C'}^*(C')$ in $F^1$ are adjacent up to tensoring by a line bundle plays a crucial role. We will show that the same property holds for non-abelian groups as well (cf. Proposition \ref{prop:twist}). 

Let us recall that a wall $W=\overline{C}\cap\overline{C'}$ of Type 0 induces an exact sequence
\begin{equation}\label{eq:exact3}
0\to\mcl{S}\to \mcl{R}_C|_{D_W}\to \mcl{Q}\to0
\end{equation}
of flat family of $G$-equivariant $\C[V]$-modules on the unstable locus $D_W$ which extends (\ref{eq:exact1}) (see the proof of Lemma \ref{lem:unstable}).

For each connected component $D\subset D_W$, we will show that either (1) one of $\mcl{S}|_D$ and $\mcl{Q}|_D$ parametrizes rigid $\C[V]$-modules (i.e. $\C[V]$-modules $F$ with $G\text{-}\mrm{Ext}^1_{\C[V]}(F, F)=0$), or (2) $D$ is a product of (possibly reducible) curves (Lemma \ref{lem:structure}). To do this, we first take the direct summands $R_1,\, R_2$ of $R$ corresponding to $x\in D$ as in the construction of $Z_Q$. We then consider the closed subscheme $\mcl{V}_{\mcl{S}}$ (resp. $\mcl{V}_{\mcl{Q}}$) of $\mcl{N}$ given by vanishing of coordinates corresponding to all arrows but ones from $R_1$ (resp. $R_2$) to itself. Now we define coarse moduli space $\bar{\mfk{M}}_{\mcl{S}}$ (resp. $\bar{\mfk{M}}_{\mcl{Q}}$) as the GIT-quotient
\[
\mcl{V}_{\mcl{S}}/\!/_{\theta_1}GL_{R_1}\text{ (resp. }\mcl{V}_{\mcl{Q}}/\!/_{\theta_2} GL_{R_2}) 
\]
where $\theta_i$ is the stability condition induced from $\theta_0\in W$ via the natural inclusion $GL_{R_i}\to GL_R$. We denote by $\mfk{M}_{\mcl{S}}$ (resp. $\mfk{M}_{\mcl{Q}}$) the reduced induced subscheme of $\bar{\mfk{M}}_{\mcl{S}}$ (resp. $\bar{\mfk{M}}_{\mcl{Q}}$). Let $\bar{D}$ be the image of $D$ with the reduced structure under the morphism $\mfk{M}_C\to \mfk{M}_W$. Note that $\bar{D}$ is possibly non-normal.

\begin{lem}\label{lem:product}
With the notation as above, there are natural projections from $D_W$ to $\bar{\mfk{M}}_{\mcl{S}}$ and $\bar{\mfk{M}}_{\mcl{Q}}$. Moreover, the bijective map $D\to \bar{D}$ factors through the induced map $D\to \mfk{M}_{\mcl{S}}\times \mfk{M}_{\mcl{Q}}$ which is also bijective.
\end{lem}

\begin{proof}
First note that the projections are well-defined since fibers of $\mcl{S}|_D$ and $\mcl{Q}|_D$ are semistable with respect to $\theta_1$ and $\theta_2$ respectively. We also have a morphism $\mfk{M}_{\mcl{S}}\times \mfk{M}_{\mcl{Q}}\to \bar{D}$ which comes form the assignment $(F_1, F_2)\mapsto F_1\oplus F_2$. This is bijective by the choice of $\theta_1, \theta_2$, and clearly the composition of this map with $D\to \mfk{M}_{\mcl{S}}\times \mfk{M}_{\mcl{Q}}$ is equal to the map $D\to \bar{D}$.
\end{proof}

The following lemma indicates that there can be two essentially different situations for $D$.

\begin{lem}\label{lem:structure}
For each connected component $D\subset D_W$, exactly one of the following two situations occurs:\\
(1) Either $\mcl{S}|_D$ or $\mcl{Q}|_D$ parametrizes rigid $\C[V]$-modules.\\
(2) $D$ is a product of two curves $C_1, C_2$ where each $C_i$ is a tree of rational curves of Dynkin type, that is, $C_i$ is isomorphic to the (reduced) exceptional locus of the minimal resolution of an ADE surface singularity.
\end{lem}

\begin{proof}
Assume that there exists a general point $x\in D$ such that its corresponding $\C[V]$-modules $S\subset F$ and $Q=F/S$ are not rigid. The 2-dimensional tangent space $T_x D$ of $D$ (as a relative Quot scheme) is isomorphic to the kernel of the natural map
$$G\text{-}\mrm{Ext}^1_{\C[V]}(F, F)\to G\text{-}\mrm{Ext}^1_{\C[V]}(S, Q)$$
by the $G$-equivariant version of \cite[Proposition 4.4.4]{S}. Since we have
$$G\text{-}\Hom(S, Q)=0 \quad \text{and} \quad G\text{-}\mrm{Ext}^1(S, Q)\cong G\text{-}\mrm{Ext}^1(Q, S)\cong \C,$$
one sees that $T_x D$ is isomorphic to $G\text{-}\mrm{Ext}^1_{\C[V]}(S, S)\oplus G\text{-}\mrm{Ext}^1_{\C[V]}(Q, Q)$ and hence
$$\dim G\text{-}\mrm{Ext}^1_{\C[V]}(S, S)=\dim G\text{-}\mrm{Ext}^1_{\C[V]}(Q, Q)=1.$$
Therefore, $D$ is a product of two curves. Note that even if $S$ and $Q$ are non-rigid, the moduli space $\bar{\mfk{M}}_{\mcl{S}}$ or $\bar{\mfk{M}}_{\mcl{Q}}$ might be one point with nonreduced structure.

By Lemma \ref{lem:product}, we only have to show that $\mfk{M}_{\mcl{S}}$ and $\mfk{M}_{\mcl{Q}}$ are curves of Dynkin type when they are one-dimensional. Indeed, this particularly implies that each irreducible component of $\mfk{M}_{\mcl{S}} \times \mfk{M}_{\mcl{Q}}$ is smooth and thus $D\to \mfk{M}_{\mcl{S}} \times \mfk{M}_{\mcl{Q}}$ is an isomorphism.

By moving from the relative interior of $W$ to its smaller faces corresponding submodules of the fiber $\mcl{S}_x$ for general $x\in D$, we obtain $\theta\in\Theta$ and its associated birational map $\alpha_\theta: \mfk{M}_C\dashrightarrow \mfk{M}_\theta$ such that $\alpha_\theta$ is defined on $D$ and that $\alpha_\theta$ only collapses the first factor $\mfk{M}_{\mcl{S}}$ of $D$ and in particular the image $\alpha_\theta(D)$ is the curve $\mfk{M}_{\mcl{Q}}$. This is indeed possible since we just have to avoid faces contained in facets of $W$ which induce divisorial contractions corresponding to curves in $\mfk{M}_{\mcl{Q}}$. Note that such a facet is unique for each contraction and that moving from a face to its codimension-one face gives an isomorphism or a contraction of relative Picard number one between the corresponding (normalized) moduli spaces, and in particular each process does not contract the two factors at the same time.

Therefore, the curve $\mfk{M}_{\mcl{S}}$ is identified with the fiber, under $\alpha_\theta$, over a general point of the one-dimensional singular locus of $\bar{\mfk{M}}_\theta$. Then it must be a tree of rational curves of Dynkin type since $\bar{\mfk{M}}_\theta$ has only canonical singularities and $\alpha_\theta$ is crepant (cf. \cite[Corollary 1.14]{Re1}). The same argument works for $\mfk{M}_{\mcl{Q}}$ as well.
\end{proof}

For given $\mcl{E}\in K_0(X)$, we consider the automorphism $\mcl{T}_\mcl{E}$ of $K(X)$ defined by
$$\mcl{T}_\mcl{E}(\xi)=\xi-\chi(\xi, \mcl{E})\cdot \mcl{E}.$$
Note that $\mcl{T}_\mcl{E}$ has positive determinant and fixes the hyperplane
$$\mcl{E}^\perp:=\{\xi\in K(X) \mid \chi(\xi, \mcl{E})=0\}$$
pointwise. We denote by $\mcl{T}'_\mcl{E}$ the inverse of $\mcl{T}_\mcl{E}$. We will show that the automorphism $\phi_{C'}^*\circ (\phi_{C}^*)^{-1}$ of $F^1$ coincides with $\mcl{T}_\mcl{E}$ or $\mcl{T}'_\mcl{E}$ for some $\mcl{E}\in K_0(X)$ up to an automorphism obtained by tensoring by a line bundle of $X$ (Proposition \ref{prop:twist}).

For the families $\mcl{S}$ and $\mcl{Q}$ in (\ref{eq:exact3}) associated to the wall $W$, let $\mcl{S}_\rho$ and $\mcl{Q}_\rho$ be their $\rho$-isotypic components respectively.

\begin{lem}\label{lem:cohomology}
For any $\rho, \sigma\in\mrm{Irr}(G)$, we have $H^i(D_W, \mcl{S}_\rho^\vee\otimes \mcl{Q}_\sigma)=0$ for all $i$.
\end{lem}

\begin{proof}
We first note that the $G$-equivariant bundle $\bigoplus_{\rho\in \mrm{Irr}(G)} \mcl{R}^\vee_\rho \otimes \rho^*$ on $X\times \C^3$ can be regarded as the universal family of the moduli space $\mfk{M}_{\theta^*}$ for $\theta^*\in \Theta$ defined by $\theta^*(\rho)=-\theta(\rho^*)$ with $\theta\in C$ by \cite[Lemma 2.6(ii)]{CI}. Note that $\mfk{M}_{\theta^*}\cong \mfk{M}_C$ and that the associated sub and quotient bundles $\mcl{S}^*,\,\mcl{Q}^*$ of $\mcl{R}_{\theta^*}|_{D_W}\cong \mcl{R}_C^\vee|_{D_W}$ are isomorphic to $\mcl{Q}^\vee$ and $\mcl{S}^\vee$ respectively.

By dualizing the inverse
\begin{equation}\label{eq:K}
0\to \mcl{R}_C(-D_W)\to \mcl{K}\to \mcl{S}\to 0
\end{equation}
of the elementary transformation
\begin{equation}\label{eq:R}
0\to \mcl{K}\to \mcl{R}_C\to \mcl{Q}\to 0
\end{equation}
and then tensoring with $\mcl{O}_X(-D_W)$, we obtain a sequence
$$0\to \mcl{K}^\vee (-D_W) \to \mcl{R}_C^\vee \to \mcl{S}^\vee\to 0$$
which is identified with the tautological sequence for $\mfk{M}_{\theta^*}$. Also, tensoring with $\mcl{O}_X(D_W)$ to the inverse of the elementary transformation for this sequence, we obtain
\begin{equation}\label{eq:dual}
0\to \mcl{R}_C^\vee\to \mcl{K}^\vee\to \mcl{Q}^\vee\otimes \mcl{O}_X(D_W)\to 0.
\end{equation}
Therefore, $\mcl{R}_C^\vee|_{D_W}$ admits $\mcl{Q}^\vee$ (resp. $\mcl{S}^\vee$) as the associated sub (resp. quotient) bundle while $\mcl{K}^\vee|_{D_W}$ admits $\mcl{S}^\vee$ (resp. $\mcl{Q}^\vee(D_W)$) as the associated sub (resp. quotient) bundle.

From (\ref{eq:dual}) we obtain a long exact sequence
\begin{equation}\label{eq:long}
\begin{aligned}
\cdots &\to \mrm{Ext}^i_{\mcl{O}_X}(\mcl{R}_C^\vee, \mcl{R}_C^\vee)\to \mrm{Ext}^i_{\mcl{O}_X}(\mcl{R}_C^\vee, \mcl{K}^\vee)\stackrel{f^i}{\to} \mrm{Ext}^i_{\mcl{O}_X}(\mcl{R}_C^\vee, \mcl{Q}^\vee(D_W))\\
&\to \mrm{Ext}^{i+1}_{\mcl{O}_X}(\mcl{R}_C^\vee, \mcl{R}_C^\vee)\to \cdots
\end{aligned}
\end{equation}
since $\mcl{R}_C^\vee$ is locally free. For all $i$, the map $f^i$ in (\ref{eq:long}) is identified with the composition of the restriction map
\begin{equation}\label{eq:restriction}
\mrm{Ext}^i_{\mcl{O}_X}(\mcl{R}_C^\vee, \mcl{K}^\vee)\to \mrm{Ext}^i_{\mcl{O}_X}(\mcl{R}^\vee_C, \mcl{K}^\vee|_{D_W})
\end{equation}
and the natural projection
$$\mrm{Ext}^i_{\mcl{O}_X}(\mcl{R}^\vee_C, \mcl{K}^\vee|_{D_W})\cong H^i(D_W, \mcl{R}_C\otimes \mcl{K}^\vee|_{D_W})\to H^i(D_W, \mcl{R}_C\otimes \mcl{Q}^\vee(D_W)).$$
By using the long exact sequence obtained from (\ref{eq:R}), we similarly have a map $g^i: \mrm{Ext}^i_{\mcl{O}_X}(\mcl{R}_C^\vee, \mcl{K}^\vee)\to H^i(D_W, \mcl{K}^\vee \otimes \mcl{Q})$ which is the composition of (\ref{eq:restriction}) and the projection
$$\mrm{Ext}^i_{\mcl{O}_X}(\mcl{R}^\vee_C, \mcl{K}^\vee|_{D_W})\to H^i(D_W, \mcl{K}^\vee\otimes \mcl{Q}).$$
The kernels of the two projections are the same since we have
\begin{equation}
\mrm{Ext}^i_{\mcl{O}_X}(\mcl{R}_C^\vee, \mcl{R}_C^\vee)=\mrm{Ext}^i_{\mcl{O}_X}(\mcl{K}^\vee, \mcl{K}^\vee)=0
\end{equation}
for $i>0$ and
\begin{equation}
\Hom_{\mcl{O}_X}(\mcl{R}_C^\vee, \mcl{R}_C^\vee)=\Hom_{\mcl{O}_X}(\mcl{K}^\vee, \mcl{K}^\vee)
\end{equation}
(cf. \cite[Lemma 5.4]{CI}). It follows that $H^i(D_W, \mcl{S}^\vee\otimes \mcl{Q})=0$ (and $H^i(D_W, \mcl{S}\otimes \mcl{Q}^\vee(D_W))=0$) for all $i$.
\end{proof}

Let $W=\overline{C}\cap\overline{C'}$ be a wall of Type 0 with the associated sequence (\ref{eq:exact3}) on the unstable locus $D_W$. Take the decomposition $D_W=\bigcup_i D_i$ into connected components and let $R_{1,i}$ (resp. $R_{2,i}$) be the $\C G$-modules for fibers of $\mcl{S}|_{D_i}$ (resp. $\mcl{Q}|_{D_i}$). Note that all $R_{1,i}$ (resp. $R_{2,i}$) are proportional in $R(G)$ since $D_i$'s are defined by the same wall $W$.

\begin{lem}\label{lem:proportional}
For each $i$, there exists a class $\mcl{E}_i\in K_0(X)$ of a line bundle on $D_i$ such that the image $\phi_C(\mcl{E}_i)\in K_0(\C^3)\cong R(G)$ is proportional to either $R_{1,i}$ or $R_{2,i}$.
\end{lem}

\begin{proof}
We consider the two cases (1), (2) in Lemma \ref{lem:structure} separately. For the rigid case (1), we first assume that $\mcl{Q}|_{D_i}$ parametrizes rigid modules which are indecomposable so that there exists a line bundle $L$ on $D_i$ such that $\mcl{Q}_\rho|_{D_i}\cong L^{\oplus \mrm{rank}(\mcl{Q}_\rho|_{D_i})}$ for all $\rho$. Note that $R_{2,i}=\bigoplus_{\rho\in\mrm{Irr}(G)} \rho^{\oplus\mrm{rank}(\mcl{Q}_\rho|_{D_i})}$ in $R(G)$.

We set $\mcl{E}_i:=[L^{-1}\otimes \omega_{D_i}]\in K_0(X)$. The cohomology sheaves of $\Phi_C(L^{-1}\otimes \omega_{D_i})$, which are $G$-equivariant $\C[V]$-modules supported at the origin, are computed as
$$\begin{aligned}
\Phi_C^j(L^{-1}\otimes \omega_{D_i})&\cong\bigoplus_{\rho\in\mrm{Irr}(G)}H^j(X, \mcl{R}_\rho\otimes L^{-1}\otimes \omega_{D_i})\otimes_\C \rho\\
&\cong\bigoplus_{\rho\in\mrm{Irr}(G)} H^{2-j}(D_i, \mcl{R}_\rho^\vee \otimes L)^\vee\otimes_\C \rho\\
&\cong\bigoplus_{\rho\in\mrm{Irr}(G)} H^{2-j}\left(D_i, ((L^{-1})^{\oplus\mrm{rank}(\mcl{Q}_\rho|_{D_i})}\oplus \mcl{S}_\rho^\vee) \otimes L \right)^\vee \otimes_\C \rho \\
&\cong\bigoplus_{\rho\in\mrm{Irr}(G)} H^{2-j}\left(D_i, \mcl{O}_{D_i}^{\oplus \mrm{rank}(\mcl{Q}_\rho|_{D_i})} \right)^\vee \otimes_\C \rho
\end{aligned}$$
where the last isomorphism holds by Lemma \ref{lem:cohomology}. Therefore, the class $\phi_C^* (\mcl{E}_i)$ is equal to $\chi(\mcl{O}_{D_i})\cdot R_{2,i}\in R(G)$.

If fibers of $\mcl{Q}|_{D_i}$ are rigid but not indecomposable, we can write $\mcl{Q}|_{D_i}$ as a direct sum $\bigoplus_j L_j^{\oplus m_j}$ of line bundles. Then the similar computation as above shows that $\phi_C^* (\mcl{E}_i)$ is again proportional to $R_{2,i}\in R(G)$ when we use any $L_j$ instead of $L$. Note that each direct summand of fibers of $\mcl{Q}|_{D_i}$ is proportional to $R_{2,i}$ in $R(G)$. This implies that the line bundles $L_j$ on $D_i$ are all proportional in $K(X)$. Then $L_j$'s are mutually isomorphic since they define the same class $\mcl{O}_X(-D_i)$ in $\mrm{Pic}(X)$. Therefore, we again have that $\phi_C^* (\mcl{E}_i)$ is equal to $\chi(\mcl{O}_{D_i})\cdot R_{2,i}\in R(G)$ for $\mcl{E}_i:=[L^{-1}\otimes \omega_{D_i}]$ with a line bundle $L$ on $D_i$.

When fibers of $\mcl{S}|_{D_i}$ are rigid, one can similarly show that there exists a line bundle $L'$ on $D_i$ such that $\mcl{S}_\rho|_{D_i}\cong L'^{\bigoplus \mrm{rank}(\mcl{S}_\rho|_{D_i})}$ and that $\phi_C^*([L'^{-1}])$ is equal to $\chi(\mcl{O}_{D_i})\cdot R_{1,i}\in R(G)$.

Next we consider the case (2) in Lemma \ref{lem:structure}. In this case each $D_i$ is the product of curves $C_1$ and $C_2$ which are regarded as moduli spaces parametrizing fibers of $\mcl{S}|_{D_i}$ and $\mcl{Q}|_{D_i}$ respectively. Since each $C_j$ is a union of $\mathbb{P}^1$'s, the bundle $\mcl{S}|_{D_i}$ (resp. $\mcl{Q}|_{D_i}$) on $D_i$ is a direct sum of line bundles such that the restrictions of these line bundles to fibers of the projection $D_i\to C_1$ (resp. $D_i\to C_2$) are mutually isomorphic. We take any line bundles $L_1$ and $L_2$ on $D_i$ appearing as direct summands of $\mcl{S}_\rho|_{D_i}$ and $\mcl{Q}_{\rho'}|_{D_i}$ respectively for any $\rho, \rho'$. By Lemma \ref{lem:cohomology}, we have $H^j(D_i, L_1^{-1}\otimes L_2)=0$ for all $j$. Then either
$$H^j(C_1, (L_1^{-1}\otimes L_2)|_{C_1})=0\;(\forall j)\quad \text{or} \quad H^j(C_2, (L_1^{-1}\otimes L_2)|_{C_2})=0\;(\forall j)$$
holds where $C_1$ and $C_2$ are regarded as fibers of $D_i\to C_2$ and $D_i\to C_1$ respectively. From this we see that one has either $H^j(C_2, (L_1^{-1}\otimes L'_2)|_{C_2})=0\;(\forall j)$ for all line bundles $L'_2$ appearing as direct summands of $\mcl{Q}|_{D_i}$ or $H^j(C_1, (L'^{-1}_1\otimes L_2)|_{C_1})=0\;(\forall j)$ for all line bundles $L'_1$ appearing as direct summands of $\mcl{S}|_{D_i}$. In the former case one has $\phi_C([L_2^{-1}])=R_{2,i}\in R(G)$, and in the latter case $\phi_C([L_1^{-1}]\otimes \omega _{D_i})=R_{1,i}$. This can be checked by noticing that the Euler characteristics $\chi(L'_2|_{C_2})$ (resp. $\chi(L'_1|_{C_1})$) are the same for all the direct summands $L'_2$ of $\mcl{Q}|_{D_i}$ in the former case (resp. $L'_1$ of $\mcl{S}|_{D_i}$ in the latter case). Note that the Euler characteristic of a line bundle $L$ on $C_i$ is equal to $1+\sum_k \deg(L|_{C_{i,k}})$ where $C_{i,k}$ are the irreducible components of $C_i$.
\end{proof}

\begin{prop}\label{prop:twist}
The unstable locus $D_W$ is connected. Moreover, the automorphism of $K(X)$ induced by the autoequivalence $\Phi_{C'}^{-1}\circ \Phi_{C}$ of $D(X)$ is equal to one of the four automorphisms: (a) $\mcl{T}_\mcl{E}$, (b) $\mcl{T}'_\mcl{E}$, (c) $\xi\mapsto \mcl{T}_\mcl{E}(\xi)\otimes \mcl{O}(-D_W)$, and (d) $\xi\mapsto \mcl{T}'_\mcl{E}(\xi\otimes \mcl{O}(D_W))$.\\
In particular, $\phi^*_{C'}(C')$ is adjacent to $\phi^*_C(C)\otimes L$ for some $L$ in $\mrm{Pic}^c(X)$, the subgroup of $\mrm{Pic}(X)$ generated by the classes of compact divisors.
\end{prop}

\begin{proof}
We have shown in the previous lemma that the image in $R(G)$ of a line bundle on $D_i$ is proportional to $R_{1,i}$ (and to $R_{2,i}$) modulo $\Z R$ in any case. Since $D_i$'s are linearly independent in $\mrm{Pic}(X)_\R\cong F^1/F^2$, the classes of the line bundles on $D_i$'s in $K(X)$ are also linearly independent (cf. the proof of \cite[Proposition 5.5(ii)]{CI}). Therefore, $D_W$ must be connected. The similar argument in the proof of \cite[Proposition 5.5(iii)]{CI} also shows that we have $H^1(\mcl{O}_{D_W})=H^2(\mcl{O}_{D_W})=0$ in the rigid case (1) in Lemma \ref{lem:structure}.

To summarize, the following two cases may occur:\\
(I) we have $\chi(\mcl{R}_\rho^\vee, \mcl{E})=\mrm{rank}(\mcl{Q}_\rho)$ for the class $\mcl{E}=[L^{-1}\otimes \omega_{D_W}]\in K_0(X)$ with an arbitrary line bundle $L$ on $D_W$ which is a direct summand of $\mcl{Q}$, and\\
(II) we have $\chi(\mcl{R}_\rho^\vee, \mcl{E})=\mrm{rank}(\mcl{S}_\rho)$ for the class $\mcl{E}=[L^{-1}]\in K_0(X)$ with an arbitrary line bundle $L$ on $D_W$ which is a direct summand of $\mcl{S}$.\\
In the former case (I), one has
$$\mcl{T}_\mcl{E}([\mcl{R}_\rho^\vee])=[\mcl{R}_\rho^\vee]-\mrm{rank}(\mcl{Q}_\rho)[L^{-1}\otimes \omega_{D_W}]=[\mcl{R}_\rho^\vee]-[\mcl{Q}_\rho^\vee(D)]$$
in $K(X)$. By the sequence (\ref{eq:dual}) and the fact that $\{\mcl{R}_\rho^\vee\}_\rho$ forms a basis of $K(X)$, we can conclude that $\mcl{T}_\mcl{E}$ coincides with the automorphism induced by $\Phi_{C'}^{-1}\circ \Phi_{C}$ or its composition with tensoring by $\mcl{O}_X(-D_W)$. Which of the cases (a) and (c) happens depends on whether $\mcl{Q}_{\rho_0}=0$ or $\ne0$. The latter case (II) can be treated similarly by switching the roles of $C$ and $C'$, so that we fall into the cases (b) and (d). This completes the proof.
\end{proof}

\begin{rem}
The above proposition shows that the wall-crossing gives a combination of a spherical twist and tensoring by a line bundle at least at the level of the Grothendieck group. In the rigid case (1) in Lemma \ref{lem:structure}, one can use the same argument as in the proof of \cite[Proposition 7.3]{CI} to prove the similar statement at the level of the derived category. In the case (2), however, this argument does not work since the line bundles in $\mcl{S}$ or $\mcl{Q}$ are not mutually isomorphic (while they have the same Euler characteristics) in contrast to the rigid case. By analogy with the abelian case \cite[Corollary 4.6]{CI}, it is expected that the case (2) does not occur.
\end{rem}

\noindent{\it Proof of Theorem \ref{main2}}

Since all projective crepant resolutions of $\C^3/G$ are connected with each other by a sequence of flops \cite{K}, we only have to show that every flop of $X=\mfk{M}_C$ can be realized by (not necessarily single) wall-crossings. Let $B\subset \overline{\mrm{Amp}}(X)$ be a facet corresponding to a flop $X\to X_B\gets X'$. If we can reach a GIT-chamber $C''$ through walls of Type 0 such that $C''$ admits a wall $\tilde{B}\subset \overline{C}''$ satisfying $\psi_{C''}(\tilde{B})\subset B$, then crossing $\tilde{B}$ induces the desired flop. Indeed, otherwise $\tilde{B}$ would be contained in a facet of the orbit cone of a general point of a divisor since the homomorphism $\psi_{C''}$ must change when we cross $\tilde{B}$. In such a case, however, applying Lemma \ref{lem:vertical} to $\tilde{B}$ shows that $\tilde{B}$ would induce a divisorial contraction, which is a contradiction.

As shown in \cite[\S8]{CI}, each compact exceptional divisor $E\subset X$ acts on $F^1$ by
$$\xi\mapsto \xi\otimes \mcl{O}_X(E)=\xi+[\xi \otimes \mcl{O}_E]$$
with $[\xi \otimes \mcl{O}_E]\in F^2$. Moreover, for the set of all irreducible compact exceptional divisors $S_1,\dots,S_d$ and for any $\xi\in \mrm{pr}^{-1}(\mrm{Amp}(X)\cup \mrm{relint}(B))$, the set $[\xi\otimes\mcl{O}_{S_1}],\dots,[\xi\otimes\mcl{O}_{S_d}]$ forms a basis of $F^2$ where $\mrm{Amp}(X)$ is identified with a cone in $F^1/F^2$ (see (\ref{Pic})).

The existence of the desired chamber $C''$ can be shown similarly to the proof of \cite[Proposition 8.2]{CI}. More precisely, one can show from Proposition \ref{prop:twist} that we have
\begin{equation}\label{inclusion}
\mrm{pr}^{-1}(\mrm{Amp}(X)\cup \mrm{relint}(B))\subset \bigcup_{L\in \mrm{Pic}^c(X)}\bigcup_{\mfk{M}_C\cong X} L\otimes \phi_C^*(\overline{C})
\end{equation}
noticing that GIT-chambers in $\Theta$ are finite. Thus we obtain a chamber $C''$ and its wall $\tilde{B}$ such that $\mrm{pr}(\phi_C^*(\tilde{B}))\subset B$.
\qed

\vspace{3mm}

Theorem \ref{main2} together with (\ref{inclusion}) implies that $\{\psi_C^*(\overline{C})\}_{C\subset\Theta\text{: chamber}}$ covers the movable cone of $X$. Thus, every rational contraction $X\dashrightarrow X'$ (see Subsection \ref{3.2}) is realized as the birational map $\mfk{M}_C\dashrightarrow\bar{\mfk{M}}_\theta$ induced by variation of GIT for some $\theta$. This fact has the following application:

\begin{prop}\label{prop:generate}
Let $h_1, \dots, h_\ell$ be homogeneous generators of the normalization of $\C[\mcl{W}]$. Then the associated $\ell+m$ elements
\[
\tilde{\vphi(h_1)}, \dots, \tilde{\vphi(h_\ell)}, t_1^{-r_1},\dots, t_m^{-r_m}
\]
to $\vphi(h_1), \dots, \vphi(h_\ell)\in\C[V]^{[G,G]}$ (see (\ref{gens})) generate the Cox ring $\Cox(X)\subset \C[V]^{[G,G]}[t_1^{\pm1},\dots,t_m^{\pm1}]$ of any crepant resolution $X\to\C^3/G$.
\end{prop}

\begin{proof} Let $f_1,\dots, f_{\ell'}\in\C[V]^{[G,G]}$ be homogeneous elements such that the associated elements in (\ref{gens}) generate $\Cox(X)$. We consider the ambient toric varieties $Y$ of $X$ and $\bar{Y}_\mcl{W}$ of $\bar{\mcl{W}}$ associated to $f_i$'s (see Subsection \ref{3.3}). Note that each $f_i$ corresponds to a non-exceptional torus-invariant divisor $D_i\subset Y$ and that the image of the corresponding section $f_{D_i}\in H^0(Y, \mcl{O}(-D_i))$ (i.e. $D_i=\{f_{D_i}=0\}$) in $\Cox(X)$ is equal to $\tilde{f_i}$ (cf. Remark \ref{rem:cox(V/G)}). Regarding the weight of $f_{D_i}$ as a character of $T_X$, we consider its associated GIT-quotient $Y_i$ of $\mathbb{A}^{\ell+m}$. In terms of a fan, $Y_i$ is characterized as a birational model $Y\dashrightarrow Y_i$ with $\dim\mrm{Pic}(Y_i)_\R=1$ such that the toric fan of $Y_i$ contains the cone generated by the rays corresponding to the torus-invariant divisors of $Y$ except $D_i$.

For each $i$, let $X_i\subset Y_i$ be the induced GIT-quotient of $\mrm{Spec}(\Cox(X))$. As observed above, there exists $\theta\in\Theta$ such that $X_i\cong \bar{\mfk{M}}_\theta$. Since $\bar{Y}_\mcl{W}$ is normal,  there is a homogeneous semi-invariant $f$ in $\C[\bar{\mcl{W}}]$ such that its (scheme-theoretic) zero locus in $\bar{\mfk{M}}_\theta\cong X_i$ is equal to the strict transform of $D_i\cap X$ under $X\dashrightarrow X_i$. If we take a GIT-chamber $C$ so that $\overline{C}$ contains $\theta$, this implies $\vphi_C(f)$ is equal to $\tilde{f_i}$ (up to constant multiplication). This particularly implies that each $\tilde{f_i}$ can be obtained as the associated element of $\vphi(h)$ for some homogeneous semi-invariant $h$. Thus, the claim follows.
\end{proof}

\appendix
\section{List of Symbols}
\addcontentsline{toc}{chapter}{List of Symbols}

\noindent
\begin{longtable}{l l l }
  $R$ & the regular representation of $G$ & \S\ref{2.1} \\
  $\mrm{Irr}(G)$ & the set of irreducible representations of $G$ & \S\ref{2.1} \\
  $\mcl{N}$ & the affine scheme parametrizing all $G$-constellations & \S\ref{2.1} \\
  $GL_R$ & $:=\mrm{Aut}_{\C G}(R)=\prod_{\rho\in\mrm{Irr}(G)}GL_{\dim \rho}(\C)$ & \S\ref{2.1}  \\
  $PGL_R$ & $:=GL_R/\C^*$ & \S\ref{2.1} \\
  $\bar{\Theta}$ & $:=\{\theta\in\mrm{Hom}_\mathbb{Z}(R(G),\Z)\mid \theta(R)=0\}$ & \S\ref{2.1} \\
 $R(G)$ & $:=\bigoplus_{\rho\in\mrm{Irr}(G)} \Z \rho$ : the representation ring of $G$ & \S\ref{2.1} \\
 $\Theta$ & $:=\bar{\Theta}\otimes_\Z \R$ : the space of stablity conditions & \S\ref{2.1}\\
 $\mcl{M}_\theta$ & the moduli space of $\theta$-semistable $G$-constellations & \S\ref{2.1}\\
 $\mcl{V}\subset \mcl{N}$ & the main irreducible component & \S\ref{2.1}\\
 $\mfk{M}_\theta$ & the coherent component of $\mcl{M}_\theta$ & \S\ref{2.1}\\
 $\mfk{M}_C$ & $:=\mfk{M}_\theta$ for $\theta$ in a GIT-chamber $C\subset \Theta$ & \S\ref{2.1}\\
 $SL_R$ & $:=\prod_{\rho\in\mrm{Irr}(G)}SL_{\dim \rho}(\C)\subset GL_R$ & \S\ref{2.2}  \\
 $Q_G$ & the McKay quiver of $G$ & \S\ref{2.3}  \\
 $\mcl{R}_C$ & the tautological bundle of $\mfk{M}_C$ & \S\ref{2.3}  \\
 $\mcl{R}_\rho$ & the isotypic component of $\mcl{R}_C$ for $\rho\in\mrm{Irr}(G)$  & \S\ref{2.3} \\
 $\Cox(X)$ & the Cox ring of a crepant resolution $X\to \C^n/G$ & \S\ref{2.3} \\
 $Ab(G)$ & $:=G/[G,G]$: the abelianization of $G$ & \S\ref{3.1}\\
 $T_X$ & the torus acting on $\mfk{X}:=\mrm{Spec}\,\Cox(X)$ & \S\ref{3.2} \\
 $\iota_\mfk{X}: \mfk{X}\to \mathbb{A}^{\ell+m}$ & the embedding given by $f_1,\dots, f_\ell\in \C[V]^{G, G}$  & \S\ref{3.3} \\
 $\iota_X: X\to Y$ & the embedding into a toric variety associated to $\iota_\mfk{X}$ & \S\ref{3.3} \\
 $\mcl{W}$ & $:=\mcl{V}/\!/SL_R$ & \S\ref{3.3} \\
 $T_\Theta$ & the acting torus on $\mcl{W}$ & \S\ref{3.3} \\
 $\vphi$ & the ring map $\C[\mcl{N}]^{SL_R}\to \Cox(V/G)\cong \C[V]^{[G,G]}$ & \S\ref{3.3} \\
 $\vphi_C$ & the ring map $\C[\mcl{N}]^{SL_R}\to \Cox(X)$ & \S\ref{3.3} \\
 $\psi_C$ & the $\R$-linear map $\Theta\to \mrm{Pic}(X)_\R$ & \S\ref{3.3} \\
 $N_\mcl{W}$ & $:=\Z^\ell\oplus\bar{\Theta}^\vee$ & \S\ref{3.3} \\
 $M_\mcl{W}$ & $:=\Hom_\Z(N_\mcl{W}, \Z)$ & \S\ref{3.3} \\
 $\mcl{W}\hookrightarrow Y_\mcl{W}$ & the embedding into a toric variety associated to $\iota_\mfk{X}$ & \S\ref{3.3} \\
 $\bar{\mcl{W}}\hookrightarrow \bar{Y}_\mcl{W}$ & the normalization of $\mcl{W}\hookrightarrow Y_\mcl{W}$ & \S\ref{3.3} \\
 $W\subset C$ & a wall (of Type 0, mostly) & \S\ref{4} \\
 $\mfk{M}_W$ & $:=\mfk{M}_{\theta_0}$ for general $\theta_0\in W$ & \S\ref{4}\\
 $\bar{\mfk{M}}_W$ & the normalization of $\mfk{M}_W$ & \S\ref{4}\\
 $D_W\subset \mfk{M}_C$ & the unstable locus associated to $W$ & \S\ref{4} \\
 $\mcl{S},\,\mcl{Q}$ & the associated sub and quotient bundles of $\mcl{R}_C|_{D_W}$ & \S\ref{4} \\
 $K(X)$ & the Grothendieck group of $X$ & \S\ref{4} \\
 $K_0(X)$ & the Grothendieck group of $X$ with compact support & \S\ref{4} \\
 $\mcl{T}_\mcl{E}$ & the automorphism of $K(X)$ defined by $\mcl{E}\in K_0(X)$ & \S\ref{4} \\
 $\mcl{T}'_\mcl{E}$ & the inverse of $\mcl{T}_\mcl{E}$ & \S\ref{4} \\
 $\mcl{S}_\rho$ & the isotypic component of $\mcl{S}$ for $\rho\in\mrm{Irr}(G)$ & \S\ref{4} \\
 $\mcl{Q}_\rho$ & the isotypic component of $\mcl{Q}$ for $\rho\in\mrm{Irr}(G)$ & \S\ref{4}
\end{longtable}

\noindent
Department of Mathematical Sciences, University of Bath, Claverton Down, Bath, BA2 7AY, UK.\\
Email address: \texttt{ry488@bath.ac.uk}


\begin{thebibliography}{99}

\bibitem[AH]{AH}
A. A'Campo-Neuen\ and\ J. Hausen, Quotients of toric varieties by the action of a subtorus, Tohoku Math. J. (2) {\bf 51} (1999), no.~1, 1--12.

\bibitem[ADHL]{ADHL}
I. Arzhantsev, U. Derenthal, J. Hausen and A. Laface, {\it Cox rings}, Cambridge Studies in Advanced Mathematics, 144, Cambridge University Press, Cambridge, 2015.

\bibitem[AG]{AG}
I. V. Arzhantsev\ and\ S. A. Ga\u{\i}fullin, Sb. Math. {\bf 201} (2010), no.~1-2, 1--21; translated from Mat. Sb. {\bf 201} (2010), no. 1, 3--24.

\bibitem[BH]{BH}
F. Berchtold\ and\ J. Hausen, GIT equivalence beyond the ample cone, Michigan Math. J. {\bf 54} (2006), no.~3, 483--515.

\bibitem[BCHM]{BCHM}
C. Birkar\ et al., Existence of minimal models for varieties of log general type, J. Amer. Math. Soc. {\bf 23} (2010), no.~2, 405--468.

\bibitem[BCRSW]{BCRSW}
G. Bellamy, A. Craw, S. Rayan, T. Schedler\ and\ H. Weiss, All 81 crepant resolutions of a finite quotient singularity are hyperpolygon spaces, to appear in Journal of Algebraic Geometry.

\bibitem[BSW]{BSW}
R. Bocklandt, T. Schedler\ and\ M. Wemyss, Superpotentials and higher order derivations, J. Pure Appl. Algebra {\bf 214} (2010), no.~9, 1501--1522.

\bibitem[BKR]{BKR}
T. Bridgeland, A. King\ and\ M. Reid, The McKay correspondence as an equivalence of derived categories, J. Amer. Math. Soc. {\bf 14} (2001), no.~3, 535--554.

\bibitem[CCL]{CCL}
S. Cautis, A. Craw\ and\ T. Logvinenko, Derived Reid's recipe for abelian subgroups of ${\rm SL}_3(\Bbb{C})$, J. Reine Angew. Math. {\bf 727} (2017), 1--48.

\bibitem[Co]{Co}
D. A. Cox, The homogeneous coordinate ring of a toric variety, J. Algebraic Geom. {\bf 4} (1995), no.~1, 17--50.

\bibitem[CLS]{CLS}
D. A. Cox, J. B. Little\ and\ H. K. Schenck, {\it Toric varieties}, Graduate Studies in Mathematics, 124, American Mathematical Society, Providence, RI, 2011.

\bibitem[Cr]{Cr}
A. Craw, An explicit construction of the McKay correspondence for $A$-Hilb $\Bbb C^3$, J. Algebra {\bf 285} (2005), no.~2, 682--705.

\bibitem[CI]{CI}
A. Craw\ and\ A. Ishii, Flops of $G$-Hilb and equivalences of derived categories by variation of GIT quotient, Duke Math. J. {\bf 124} (2004), no.~2, 259--307. 

\bibitem[CM]{CM}
A. Craw\ and\ D. Maclagan, Fiber fans and toric quotients, Discrete Comput. Geom. {\bf 37} (2007), no.~2, 251--266.

\bibitem[CMT]{CMT}
A. Craw, D. Maclagan\ and\ R. R. Thomas, Moduli of McKay quiver representations. I. The coherent component, Proc. Lond. Math. Soc. (3) {\bf 95} (2007), no.~1, 179--198.

\bibitem[CR]{CR}
A. Craw\ and\ M. Reid, How to calculate $A$-Hilb $\Bbb C^3$, in {\it Geometry of toric varieties}, 129--154, S\'{e}min. Congr., 6, Soc. Math. France, Paris.

\bibitem[DeK]{DeK}
H. Derksen\ and\ G. Kemper, {\it Computational invariant theory}, second enlarged edition, Encyclopaedia of Mathematical Sciences, 130, Springer, Heidelberg, 2015.

\bibitem[DK]{DK}
M. Donten-Bury\ and\ S. Keicher, Computing resolutions of quotient singularities, J. Algebra {\bf 472} (2017), 546--572.

\bibitem[DeW1]{DeW1}
H. Derksen\ and\ J. Weyman, Semi-invariants of quivers and saturation for Littlewood-Richardson coefficients, J. Amer. Math. Soc. {\bf 13} (2000), no.~3, 467--479.

\bibitem[DeW2]{DeW2}
H. Derksen\ and\ J. Weyman, {\it An introduction to quiver representations}, Graduate Studies in Mathematics, 184, American Mathematical Society, Providence, RI, 2017.

\bibitem[DZ]{DZ}
M. Domokos\ and\ A. N. Zubkov, Semi-invariants of quivers as determinants, Transform. Groups {\bf 6} (2001), no.~1, 9--24.

\bibitem[D]{D}
M. Donten-Bury, Cox rings of minimal resolutions of surface quotient singularities, Glasg. Math. J. {\bf 58} (2016), no.~2, 325--355.

\bibitem[DG]{DG}
M. Donten-Bury\ and\ M. Grab, Crepant resolutions of 3-dimensional quotient singularities via Cox rings, Exp. Math. {\bf 28} (2019), no.~2, 161--180.

\bibitem[DW]{DW}
M. Donten-Bury\ and\ J. A. Wi\'{s}niewski, On 81 symplectic resolutions of a 4-dimensional quotient by a group of order 32, Kyoto J. Math. {\bf 57} (2017), no.~2, 395--434.

\bibitem[G]{G}
M. Grab, Cox rings and symplectic quotient singularities with torus action, arXiv:1807.11438v2. 

\bibitem[GS]{GS}
D. Grayson and M. Stillman, Macaulay 2: A software system for research in algebraic
geometry; available at http://www.math.uiuc.edu/Macaulay2.


\bibitem[Hu]{Hu}
Y. Hu, Combinatorics and quotients of toric varieties, Discrete Comput. Geom. {\bf 28} (2002), no.~2, 151--174.

\bibitem[HK]{HK}
Y. Hu\ and\ S. Keel, Mori dream spaces and GIT, Michigan Math. J. {\bf 48} (2000), 331--348.

\bibitem[IINdC]{IINdC}
A. Ishii, Y. Ito\ and\ \'{A}. Nolla de Celis, On $G/N$-Hilb of $N$-Hilb, Kyoto J. Math. {\bf 53} (2013), no.~1, 91--130.

\bibitem[IN]{IN}
Y. Ito\ and\ H. Nakajima, McKay correspondence and Hilbert schemes in dimension three, Topology {\bf 39} (2000), no.~6, 1155--1191.

\bibitem[IR]{IR}
Y. Ito\ and\ M. Reid, The McKay correspondence for finite subgroups of ${\rm SL}(3,\bold C)$, in {\it Higher-dimensional complex varieties (Trento, 1994)}, 221-240, de Gruyter, Berlin.

\bibitem[K]{K}
Y. Kawamata, Flops connect minimal models, Publ. Res. Inst. Math. Sci. {\bf 44} (2008), no.~2, 419--423.

\bibitem[Kin]{Kin}
A. D. King, Moduli of representations of finite-dimensional algebras, Quart. J. Math. Oxford Ser. (2) {\bf 45} (1994), no.~180, 515--530. 

\bibitem[Kir]{Kir}
A. Kirillov, Jr., {\it Quiver representations and quiver varieties}, Graduate Studies in Mathematics, 174, American Mathematical Society, Providence, RI, 2016.

\bibitem[KSZ]{KSZ}
M. M. Kapranov, B. Sturmfels\ and\ A. V. Zelevinsky, Quotients of toric varieties, Math. Ann. {\bf 290} (1991), no.~4, 643--655.


\bibitem[Le]{Le}
B. Leng. The Mckay correspondence and orbifold Riemann-Roch. PhD thesis, University of Warwick, 2002.

\bibitem[Lo]{Lo}
Logvinenko, Timothy Natural G-constellation families. Doc. Math. 13 (2008), 803-823.

\bibitem[M]{M}
M. Maruyama, On a family of algebraic vector bundles, in {\it Number theory, algebraic geometry and commutative algebra, in honor of Yasuo Akizuki}, 95--146, Kinokuniya, Tokyo.

\bibitem[MFK]{MFK}
D. Mumford, J. Fogarty\ and\ F. Kirwan, {\it Geometric invariant theory}, third edition, Ergebnisse der Mathematik und ihrer Grenzgebiete (2), 34, Springer-Verlag, Berlin, 1994.


\bibitem[NdCS]{NdCS}
\'{A}. Nolla de Celis\ and\ Y. Sekiya, Flops and mutations for crepant resolutions of polyhedral singularities, Asian J. Math. {\bf 21} (2017), no.~1, 1--45.

\bibitem[O]{O}
R. Ohta, On the relative version of Mori dream spaces, Eur. J. Math. {\bf 8} (2022), suppl. 1, S147--S181.

\bibitem[Re1]{Re1}
M. Reid, Canonical $3$-folds, in {\it Journ\'{e}es de G\'{e}ometrie Alg\'{e}brique d'Angers, Juillet 1979/Algebraic Geometry, Angers, 1979}, 273--310, Sijthoff \& Noordhoff, Alphen aan den Rijn.

\bibitem[Re2]{Re2}
M. Reid. Mckay correspondence. Proc. of Algebraic Geom. Symposium (Kinosaki, November 1996), pages 14-41, 1997.

\bibitem[Ro]{Ro}
M. Rossi, Embedding non-projective Mori dream space, Geom. Dedicata {\bf 207} (2020), 355--393.

\bibitem[S]{S}
E. Sernesi, {\it Deformations of algebraic schemes}, Grundlehren der mathematischen Wissenschaften, 334, Springer, Berlin, 2006.

\bibitem[SV]{SV}
A. Schofield\ and\ M. van den Bergh, Semi-invariants of quivers for arbitrary dimension vectors, Indag. Math. (N.S.) {\bf 12} (2001), no.~1, 125--138.

\bibitem[W]{W}
M. Wemyss, Flops and clusters in the homological minimal model programme. Invent. Math. {\bf 211} (2018), no.~2, 435--521.

\bibitem[Y1]{Y1}
R. Yamagishi, On smoothness of minimal models of quotient singularities by finite subgroups of $SL_n(\Bbb C)$, Glasg. Math. J. {\bf 60} (2018), no.~3, 603--634.

\bibitem[Y2]{Y2}
R. Yamagishi, Moduli of $G$-constellations and crepant resolutions I: the abelian case, in {\it McKay correspondence, mutation and related topics}, 159--193, Adv. Stud. Pure Math., 88, Math. Soc. Japan, Tokyo.

\end{thebibliography}
\end{document}